# Efficient Finite Difference WENO Scheme for Hyperbolic Systems with Non-Conservative Products

**By**

**Dinshaw S. Balsara[1,2], Deepak Bhoriya[1], Chi-Wang Shu[3] and Harish Kumar[4]**

**[1]Physics Department, [2]ACMS Department, University of Notre Dame**

**[3]Division of Applied Mathematics, Brown University**

**[4]Department of Mathematics, Indian Institute of Technology, Delhi**

**Abstract**

Higher order finite difference Weighted Essentially Non-Oscillatory (WENO) schemes have been constructed for conservation laws. For multidimensional problems, they offer high order accuracy at a fraction of the cost of a finite volume WENO or DG scheme of comparable accuracy. This makes them quite attractive for several science and engineering applications. But, to the best of our knowledge, such schemes have not been extended to non-linear hyperbolic systems with non-conservative products. In this paper, we perform such an extension which improves the domain of applicability of such schemes. The extension is carried out by writing the scheme in fluctuation form. We use the HLLI Riemann solver of Dumbser and Balsara (2016) as a building block for carrying out this extension. Because of the use of an HLL building block, the resulting scheme has a proper supersonic limit. The use of anti-diffusive fluxes ensures that stationary discontinuities can be preserved by the scheme, thus expanding its domain of applicability. Our new finite difference WENO formulation uses the same WENO reconstruction that was used in classical versions, making it very easy for users to transition over to the present formulation.

For conservation laws, the new finite difference WENO is shown to perform as well as the classical version of finite difference WENO, with two major advantages:- 1) It can capture jumps in stationary linearly degenerate wave families exactly. 2) It only requires the reconstruction to be applied once. Several examples from hyperbolic PDE systems with non-conservative products are shown which indicate that the scheme works and achieves its design order of accuracy for smooth multidimensional flows. Stringent Riemann problems and several novel multidimensional problems that are drawn from compressible Baer-Nunziato multiphase flow, multiphase debris flow and two-layer shallow water equations are also shown to document the robustness of the method. For some test problems that require well-balancing we have even been able to apply the scheme without any modification and obtain good results. Many useful PDEs may have stiff relaxation source terms for which the finite difference formulation of WENO is shown to provide some genuine advantages.

## I) Introduction

Finite volume Essentially Non-Oscillatory (ENO) schemes were first invented in the pathbreaking paper by Harten et al [32]. The invention of finite difference Essentially Non-Oscillatory (FD-ENO) schemes by Shu and Osher [55], [56] followed soon on the heels of this important invention. The FD-ENO schemes offered high order accuracies at a price point that is much lower than their finite volume counterparts of comparable order for multidimensional problems. Early ENO schemes suffered from the deficiency that certain problems could cause very rapid switching of the stencil, causing a loss of accuracy. Weighted Essentially Non-Oscillatory (WENO) schemes were invented to overcome this deficiency (Liu, Osher & Chan [43], Jiang & Shu [40]). The methods were extended to eleventh order by Balsara & Shu [3] and much later to seventeenth order by Gerolymos, Sénéchal & Vallet [28]. It was found that WENO schemes could suffer a loss of accuracy at critical points, and a way out of this problem was presented in Henrick, Aslam & Powers [34], Borges *et al*. [12] and Castro *et al*. [13]. For a comprehensive review of WENO schemes, see Shu [57] [58].

While the above-mentioned FD-WENO schemes emphasize accuracy, the recent trend in scheme design has also begun to emphasize stability. The central WENO (CWENO) schemes (Levy, Puppo & Russo [42], Cravero & Semplice [18], Semplice, Coco & Russo [54]), which have also been extended to unstructured meshes by Friedrichs [27], Käser and Iske [38] and Dumbser and Käser [20], fall in this category of WENO schemes. Such schemes emphasize the central, and most stable, stencil over and above all the other stencils. The WENO-ZQ scheme of Zhu and Qiu [63] non-linearly hybridized a central stencil with an extrema-clipping Van Albada limiter. Multiresolution WENO, which permits the use of a hierarchy of central stencils has also been developed by Zhu and Shu [64]. In Balsara, Garain and Shu [3] and Balsara *et al*. [4] we presented closed form expressions for the one-dimensional stencil operations which also produce smoothness indicators that can be expressed as the sum of perfect squares. Since the resulting smoothness indicators can be evaluated with minimal operation count, it greatly enhances the stability and speed of high order FD-WENO schemes. Another advance in Balsara, Garain and Shu [6] was the invention of WENO-AO where "AO" stands for adaptive order. This was done by non-linearly hybridizing large high order central stencils with an r=3 WENO central scheme for stability. WENO-AO has the advantage over WENO-ZQ in that it preserves local extrema. WENO-AO limiting has been extended to unstructured meshes in Balsara et al. [7]; furthermore, Boscheri and Balsara [11] were able to also use it as a limiter for Discontinuous Galerkin (DG) schemes. In the Supplement to the paper of Balsara, Samantaray and Subramanian [8] we were able to show that the smoothness indicators for finite volume WENO could also be expressed as the sum of perfect squares, thus greatly enhancing the result from Balsara, Garain and Shu [6]. Despite all the fine accomplishments of FD-WENO schemes mentioned above, all the above-mentioned FD-WENO schemes were restricted to the treatment of conservation laws.

To the best of our knowledge, FD-WENO schemes have not been extended to hyperbolic PDE systems with non-conservative products. In recent years, hyperbolic PDE systems have been presented in the literature that have non-conservative products. As a result, FD-WENO schemes have not been able to serve those communities that need to treat hyperbolic PDE systems with non-conservative products. The *first goal* of this paper is to rectify this situation. As a result, we seek a solution methodology that can handle non-conservative products in a properly upwinded fashion when they are present. These classes of PDE systems include Baer-Nunziato compressible multiphase flow, multiphase debris flow and two layer shallow water equations. A study of these hyperbolic PDE systems indicates that they all apply to physical systems that have the property that given enough time they can tend to an equilibrium steady state. For example, given enough time, a lake of water should subside to a non-moving steady state equilibrium. Extensive study of such hyperbolic PDEs has shown that their ability to retrieve this "lake-at-rest" property is strongly linked to their ability to preserve stationary discontinuities exactly. Therefore, the *second goal* of this paper is to arrive at FD-WENO schemes that can preserve stationary discontinuities exactly. To that end, we are helped along by the relatively recent invention by Dumbser and Balsara [26] of an HLLI Riemann solver that can handle non-conservative products and retrieve the "lake-at-rest" property. Indeed, the FD-WENO scheme presented here can be viewed as a higher order extension of the work in Dumbser and Balsara [26]. What is indeed novel here is that we explicitly show the steps that allow us to justify that the scheme presented here is indeed high order accurate while satisfying the two goals set forth in this paragraph. Many of the PDEs discussed here can be endowed with stiff relaxation source terms. For problems with stiff source terms FD-WENO shows some especially good advantages, and that facet of the scheme is also documented here.

The plan of the rest of this paper is as follows. In Section II we present a discussion of the linearized formulation for FD-WENO; in this Section we show its shortcomings and build perspective. In Section III we present a quick introduction to the HLLI Riemann solver by Dumbser and Balsara [26]. Section IV extends this to higher order, presenting us with the higher order FD-WENO scheme that can handle hyperbolic PDE systems with non-conservative products. Section V explicitly justifies why our formulation is higher order accurate. Section VI presents a pointwise description of the FD-WENO scheme and its implementation. Section VII provides an accuracy study, showing that the scheme meets is design accuracy for various hyperbolic PDEs. Section VIII presents one-dimensional Riemann problems for hyperbolic PDE systems with non-conservative products. The same section also shows the ability of our FD-WENO scheme to capture stationary, linearly degenerate discontinuities and the ability to handle problems that require well-balancing. Section IX presents several multidimensional test problems, many of which are new to the literature. Section X shows that our FD-WENO scheme takes especially well to the inclusion of stiff source terms. Section XI offers conclusions.

## II) Preliminary Discussion of Linearized Formulation for FD-WENO

Finite Difference WENO schemes are extremely efficient compared to the finite volume alternatives for solving hyperbolic PDE systems with high accuracy in multiple dimensions. FD-WENO schemes have already been formulated for hyperbolic systems in conservation form (Shu and Osher [55], [56]); i.e. systems of the form $\partial_t \mathbf{U} + \partial_x \mathbf{F} = 0$ in one dimension. Up to now, and to the best of our knowledge, there have been no formulations of FD-WENO for hyperbolic systems with non-conservative products; i.e. systems of the form $\partial_t \mathbf{U} + \partial_x \mathbf{F} + \mathbf{C} \partial_x \mathbf{U} = 0$ in one dimension. Overcoming this deficiency is the *most important goal* of this paper. The vector "**U**" is the entity that is evolved in time on a computational mesh and is usually referred to the vector of primal variables for the PDE system.

In order to intuitively understand the FD-WENO schemes that we arrive at for treating hyperbolic systems with non-conservative products, it helps to start with a formulation of FD-WENO for linear systems. This process will lead us to the design of FD-WENO for hyperbolic systems with non-conservative products. In this section we will only focus on spatial discretization and assume that temporal accuracy is achieved with a suitably accurate SSP Runge-Kutta time discretization. This section is only intended to provide some stage-setting for the subsequent sections.

Let us establish some notation. Please see Fig. 1. It shows a small sub-section of the mesh function in a few adjacent zones. The zones are labeled by "$i-1, i, i+1$" etc. and their zone centers are denoted by "$x_{i-1}, x_i, x_{i+1}$" etc. The associated mesh functions are specified in pointwise fashion at the zone centers and are labeled by "$\mathbf{U}_{i-1}, \mathbf{U}_i, \mathbf{U}_{i+1}$" etc. Here "$\mathbf{U}$" is a vector of primal variables for the hyperbolic PDE that we are considering. The zones are assumed to be uniform with a size $\Delta x$; and for the moment we consider a one-dimensional mesh. Using our familiar FD-WENO reconstruction strategy, as applied to the point values of the mesh function, we can obtain a suitably high order reconstruction within each zone. (Sections 2 and 3 of Balsara, Garain and Shu [6] provide full implementation-specific details on the FD-WENO reconstruction strategy that will be primarily used in this paper.) The reconstruction within zone "*i*" gives us reconstructed $\hat{\mathbf{U}}_{i+1/2}^{-}$ and $\hat{\mathbf{U}}_{i-1/2}^{+}$ at the right and left boundaries of the zone being considered; please also see Fig. 1. We will use a caret to denote such reconstructed variables. Note that $\hat{\mathbf{U}}_{i+1/2}^{-}$ is available at the left side of the zone boundary $x_{i+1/2}$ and $\hat{\mathbf{U}}_{i-1/2}^{+}$ is available at the right side of the zone boundary $x_{i-1/2}$. Similarly, using our FD-WENO reconstruction within the zone "*i+1*", we obtain $\hat{\mathbf{U}}_{i+3/2}^{-}$ and $\hat{\mathbf{U}}_{i+1/2}^{+}$ at the right and left boundaries of that zone. Note that $\hat{\mathbf{U}}_{i+1/2}^{+}$ is available at the right side of the zone boundary $x_{i+1/2}$. Likewise, using our FD-WENO reconstruction within the zone "*i-1*", we obtain $\hat{\mathbf{U}}_{i-1/2}^{-}$ and $\hat{\mathbf{U}}_{i-3/2}^{+}$ at the right and left boundaries of that zone. Consequently, $\hat{\mathbf{U}}_{i-1/2}^{-}$ is available at the left side of the zone boundary $x_{i-1/2}$. Because we are dealing with non-conservative products, we will eventually want to work with fluctuations. To that end, we assume that a Riemann solver with left and right states given by $\hat{\mathbf{U}}_{i+1/2}^{-}$ and $\hat{\mathbf{U}}_{i+1/2}^{+}$ is applied at zone boundary $x_{i+1/2}$ and it produces a resolved state within the Riemann fan given by $\mathbf{U}_{i+1/2}^{*}$. Likewise, we assume that a Riemann solver with left and right states given by $\hat{\mathbf{U}}_{i-1/2}^{-}$ and $\hat{\mathbf{U}}_{i-1/2}^{+}$ is applied at zone boundary $x_{i-1/2}$ and it produces a resolved state within the Riemann fan given by $\mathbf{U}_{i-1/2}^{*}$. Since we are hoping to produce a very light-weight scheme, we assume that the Riemann solver will be something very simple like HLL or LLF with a single resolved state within the Riemann fan. This completes our description of Fig. 1. In the next section we will provide all the details for obtaining $\mathbf{U}^{*}_{i+1/2}$ from $\hat{\mathbf{U}}^{-}_{i+1/2}$ and $\hat{\mathbf{U}}^{+}_{i+1/2}$ at any given zone boundary $x_{i+1/2}$ for the HLL or LLF Riemann solver.

Now say that we are considering a one-dimensional linear system so that the flux is given by $\mathbf{F} = \mathbf{A}\mathbf{U}$, with "$\mathbf{A}$" being a constant square matrix. We can write the linear system as

$$\partial_t \mathbf{U} + \partial_x \mathbf{F} = 0 \qquad \Leftrightarrow \qquad \partial_t \mathbf{U} + \mathbf{A}\partial_x \mathbf{U} = 0 \tag{2.1}$$

For this linear system, the simplest 1st order accurate upwind Lax-Friedrichs (LF) spatial discretization that is also formally conservative, is given by

$$\partial_t \mathbf{U}_i = -\mathbf{A}^{+} \frac{\left(\mathbf{U}_i - \mathbf{U}_{i-1}\right)}{\Delta x} - \mathbf{A}^{-} \frac{\left(\mathbf{U}_{i+1} - \mathbf{U}_i\right)}{\Delta x} \quad \text{with } \mathbf{A}^{+} \equiv \frac{1}{2}\left(\mathbf{A} + S\ \mathbf{I}\right) \quad \text{and } \mathbf{A}^{-} \equiv \frac{1}{2}\left(\mathbf{A} - S\ \mathbf{I}\right) \tag{2.2}$$

Here "**I**" is the identity matrix. Furthermore, "*S*" is a signal speed that is given by the largest eigenvalue of $\left|\mathbf{A}\right|$. As a result, $\mathbf{A}^{+}$ only has non-negative eigenvalues and $\mathbf{A}^{-}$ only has non-positive eigenvalues. It is easy to write eqn. (2.2) in a parabolized form which immediately shows us that the largest CFL for the scheme in eqn. (2.2) is given by $\Delta t = \Delta x / S$.

The higher order accurate FD-WENO version of eqn. (2.2) which is properly upwinded is given by

$$\partial_t \mathbf{U}_i = -\mathbf{A}^{+} \frac{\left(\hat{\mathbf{U}}^{-}_{i+1/2} - \hat{\mathbf{U}}^{-}_{i-1/2}\right)}{\Delta x} - \mathbf{A}^{-} \frac{\left(\hat{\mathbf{U}}^{+}_{i+1/2} - \hat{\mathbf{U}}^{+}_{i-1/2}\right)}{\Delta x} \tag{2.3}$$

For a linear system with a constant matrix "A", the above scheme is a conservative and high order accurate FD-WENO scheme. It reduces identically to the LLF scheme from Shu and Osher [55] or Jiang and Shu [40] when the flux is linear in the conserved variables. If it is applied in a dimension-by-dimension fashion to a multidimensional hyperbolic system with constant characteristic matrices, it will even preserve pointwise accuracy for multidimensional problems. Therefore, if it can be properly upgraded, it would be a good building block for non-linear systems with non-conservative products. Note though that it will not preserve a stationary contact discontinuity exactly, so it would be inadequate for any PDE system that needs to preserve stationary discontinuities.

Let us make the simplest, and most natural, upgrade of the scheme shown in eqn. (2.3) for a genuinely non-linear one-dimensional PDE with non-conservative products. Let us, therefore, consider

$$\partial_t \mathbf{U} + \partial_x \mathbf{F} + \mathbf{C} \partial_x \mathbf{U} = 0 \tag{2.4}$$

Where "**F**" can be a non-linear flux that depends on "**U**", and "**C**" can be a solution-dependent square matrix associated with the non-conservative products. The above equation can be rewritten as

$$\partial_t \mathbf{U} + \mathbf{A} \partial_x \mathbf{U} = 0 \quad \text{with} \quad \mathbf{A} \equiv \frac{\partial \mathbf{F}}{\partial \mathbf{U}} + \mathbf{C} \equiv \mathbf{B} + \mathbf{C} \tag{2.5}$$

Therefore, we can say that for the special case where the square matrices "**B**" and "**C**" are constant matrices, we just have a formal split between $\partial_x \mathbf{F}$ and $\mathbf{C} \partial_x \mathbf{U}$. Therefore, in that special case, the formulation from the previous paragraph will hold exactly. To accommodate a non-linear PDE with non-conservative products, let us denote the solution-dependent characteristic matrix by $\mathbf{A}_i \equiv \mathbf{A}\left(\mathbf{U}_i\right)$. One might be instinctively tempted to upgrade the scheme shown in eqn. (2.3) as

$$\partial_t \mathbf{U}_i = -\mathbf{A}_i^+ \frac{\left(\hat{\mathbf{U}}_{i+1/2}^- - \hat{\mathbf{U}}_{i-1/2}^-\right)}{\Delta x} - \mathbf{A}_i^- \frac{\left(\hat{\mathbf{U}}_{i+1/2}^+ - \hat{\mathbf{U}}_{i-1/2}^+\right)}{\Delta x} \quad \text{with } \mathbf{A}_i^+ \equiv \frac{1}{2}\left(\mathbf{A}_i + S_i \ \mathbf{I}\right) \text{ and } \mathbf{A}_i^- \equiv \frac{1}{2}\left(\mathbf{A}_i - S_i \ \mathbf{I}\right) \tag{2.6}$$

In the above equation it is beneficial to set $S_i$ to be the maximum wave speed obtained not just from $\left|\mathbf{A}_i\right|$ at the zone "*i*" but also from all the neighboring zones that might contribute to the stencils used for the update in eqn. (2.6). Eqn. (2.6) might even retain higher order accuracy and work well for systems with mild non-linearities. However, as with eqn. (2.3), the scheme in eqn. (2.6) cannot preserve stationary contacts exactly. For that reason, eqn. (2.6) is only helpful for sharpening our discussion; but it is otherwise unacceptable.

## III) Preserving Stationary Discontinuities – Anti-Diffusive Terms for the HLL and LLF Riemann Solvers

There are many application domains that use hyperbolic PDEs with non-conservative products. In many of these fields, the non-conservative products are sub-dominant compared to the flux terms. However, it is easy to see from the structure of those PDEs that the non-conservative products usually couple to linearly degenerate characteristic fields. Wave structures associated with such characteristic fields do not self-steepen. Indeed, if the simulation is run with a plain-vanilla scheme, wave structures associated with such wave families will diffuse ever so slowly even when a very high order scheme is used. This tendency becomes even more pernicious when, say for instance, discontinuous bathymetry is present in any variant of a shallow water equation. In that case, if isolated stationary contact discontinuities are not exactly preserved, the simulation might never achieve steady state. This is the "lake-at-rest" property. It has undergone extensive study and it has been realized that the ability to preserve isolated stationary contact discontinuities is crucial for achieving the "lake-at-rest" property. The ability to preserve isolated stationary contact discontinuities is also crucial for arriving at a well-balanced treatment of PDEs that are required to approach and hold a steady state. (We don't consider well-balancing in this paper. However, please see Käppeli and Mishra [35], Berberich *et al*. [10], Grosheintz-Laval and Käppeli [31] and the recent review by Käppeli [36], all of which reinforce the need for preserving isolated stationary contact discontinuities. All those authors identify the need for designing Riemann solvers that can preserve such contact discontinuities.)

The formulation presented in the previous Section will not, by itself, preserve isolated stationary contact discontinuities. However, it is based on the LLF Riemann solver. The LLF Riemann solver itself is a special case of the HLL Riemann solver. For an HLL Riemann solver, Dumbser and Balsara ([26]; DB16 henceforth) have presented a strategy for restoring isolated stationary contact discontinuities via adding an anti-diffusive term to the flux or the fluctuations. Therefore, we adapt that formulation of the anti-diffusive flux term to this work. It should be noted too, that the anti-diffusive contribution is always conservative.

Notice from Fig. 1 that at each zone boundary, say $x_{i+1/2}$, we have two states on the left and right of that boundary, given by $\hat{\mathbf{U}}^{-}_{i+1/2}$ and $\hat{\mathbf{U}}^{+}_{i+1/2}$ respectively. Associated with those states, we can evaluate a single resolved state $\mathbf{U}^*_{i+1/2}$ associated with the Riemann fan at that zone boundary. Let the left-most wave emanating from the Riemann fan be denoted by $S_{L;i+1/2}$ and let the right-most wave emanating from the Riemann fan be denoted by $S_{R;i+1/2}$ . The Riemann solver of DB16 has three very nice properties that are very useful to the subsequent development of this paper:- First, it can be written in a fluctuation form which is general enough to accommodate conservative and non-conservative formulations. When non-conservative products are absent, it reduces to a conservation form. Second, it allows for $S_{L;i+1/2}$ to be minimized and $S_{R;i+1/2}$ to maximized within some reasonable physical limits and, even with that adjustment, it will still produce good results. This will be very useful for the order property of the FD-WENO scheme that we will design. Third, it allows for anti-diffusive terms to be designed which will exactly preserve a stationary discontinuity in any wave family. This will be very useful for incorporating the "lake-at-rest" property in the final scheme that we will design.

While the resolved $\mathbf{U}^*_{i+1/2}$ for an HLL or LLF Riemann solver is given in DB16, we provide it explicitly in this paragraph for the sake of completeness.

$$\begin{aligned}\mathbf{U}^*_{i+1/2} = &\frac{\left(S_{R;i+1/2}\hat{\mathbf{U}}^+_{i+1/2} - S_{L;i+1/2}\hat{\mathbf{U}}^-_{i+1/2}\right) - \left(\mathbf{F}\left(\hat{\mathbf{U}}^+_{i+1/2}\right) - \mathbf{F}\left(\hat{\mathbf{U}}^-_{i+1/2}\right)\right)}{S_{R;i+1/2} - S_{L;i+1/2}} \\ &- \frac{1}{S_{R;i+1/2} - S_{L;i+1/2}}\left[\tilde{\mathbf{C}}\left(\hat{\mathbf{U}}^-_{i+1/2}, \mathbf{U}^*_{i+1/2}\right)\left(\mathbf{U}^*_{i+1/2} - \hat{\mathbf{U}}^-_{i+1/2}\right) + \tilde{\mathbf{C}}\left(\mathbf{U}^*_{i+1/2}, \hat{\mathbf{U}}^+_{i+1/2}\right)\left(\hat{\mathbf{U}}^+_{i+1/2} - \mathbf{U}^*_{i+1/2}\right)\right]\end{aligned} \tag{3.1}$$

It is easy to see that when $\mathbf{C} = 0$ the above equation reduces to the conventional formula for the HLL Riemann solver; which also makes it very computationally efficient. The matrices $\tilde{\mathbf{C}}\left(\hat{\mathbf{U}}^-_{i+1/2}, \mathbf{U}^*_{i+1/2}\right)$ and $\tilde{\mathbf{C}}\left(\mathbf{U}^*_{i+1/2}, \hat{\mathbf{U}}^+_{i+1/2}\right)$ are obtained via path integration across the jumps in the Riemann fan. Theory does not provide us any further guidance as to the optimal choice of path in phase space. For PDEs where the matrix "C" is linear in the flow variables, this choice makes physical sense. However, in general, a straight line path integrated with four point Gauss-Lobatto quadrature is often sufficient. Sometimes, one might indeed find additional support from entropy considerations (Peshkov and Romenski [47], Peshkov *et al*. [48]) but such methods have not yet been developed into an easily-implementable technology.) Therefore, to keep things simple, we define the matrices by

$$\tilde{\mathbf{C}}\left(\hat{\mathbf{U}}_{i+1/2}^{-}, \mathbf{U}_{i+1/2}^{*}\right) \equiv \int_{p=0}^{p=1} \mathbf{C}\left(\mathbf{\Phi}(p)\right) dp \quad \text{with the path in phase space given by} \quad \mathbf{\Phi}(p) \equiv \hat{\mathbf{U}}_{i+1/2}^{-} + p\left(\mathbf{U}_{i+1/2}^{*} - \hat{\mathbf{U}}_{i+1/2}^{-}\right) \tag{3.2}$$

and

$$\tilde{\mathbf{C}}\left(\mathbf{U}_{i+1/2}^{*}, \hat{\mathbf{U}}_{i+1/2}^{+}\right) \equiv \int_{p=0}^{p=1} \mathbf{C}\left(\mathbf{\Phi}(p)\right) dp \quad \text{with the path in phase space given by} \quad \mathbf{\Phi}(p) \equiv \mathbf{U}_{i+1/2}^{*} + p\left(\hat{\mathbf{U}}_{i+1/2}^{+} - \mathbf{U}_{i+1/2}^{*}\right) \tag{3.3}$$

Notice that when non-conservative products are involved, the eqn. (3.1) is implicit in the resolved state, $\mathbf{U}_{i+1/2}^{*}$. While an inexact Newton iteration strategy has been designed for its solution in DB16, our experience has been that a few fixed point iterations are sufficient to drive the left hand side in eqn. (3.1) to convergence. The convergence is indeed quite rapid and those who feel the need for exact convergence can indeed drive the process to convergence with a few more iterations.

The starting point of the iteration can be obtained by integrating over a straight line path that goes from $\hat{\mathbf{U}}_{i+1/2}^{-}$ to $\hat{\mathbf{U}}_{i+1/2}^{+}$ without the introduction of an intermediate state. Therefore, we have for the first iterate:-

$$\begin{aligned} \mathbf{U}_{i+1/2}^{*} = & \frac{\left(S_{R;i+1/2}\hat{\mathbf{U}}_{i+1/2}^{+} - S_{L;i+1/2}\hat{\mathbf{U}}_{i+1/2}^{-}\right) - \left(\mathbf{F}\left(\hat{\mathbf{U}}_{i+1/2}^{+}\right) - \mathbf{F}\left(\hat{\mathbf{U}}_{i+1/2}^{-}\right)\right)}{S_{R;i+1/2} - S_{L;i+1/2}} \\ & - \frac{1}{S_{R;i+1/2} - S_{L;i+1/2}} \tilde{\mathbf{C}}\left(\hat{\mathbf{U}}_{i+1/2}^{-}, \hat{\mathbf{U}}_{i+1/2}^{+}\right)\left(\hat{\mathbf{U}}_{i+1/2}^{+} - \hat{\mathbf{U}}_{i+1/2}^{-}\right) \end{aligned} \tag{3.4}$$

Here we define the matrix $\tilde{\mathbf{C}}\left(\hat{\mathbf{U}}_{i+1/2}^{-}, \hat{\mathbf{U}}_{i+1/2}^{+}\right)$ as:-

$$\tilde{\mathbf{C}}\left(\hat{\mathbf{U}}_{i+1/2}^{-}, \hat{\mathbf{U}}_{i+1/2}^{+}\right) \equiv \int_{p=0}^{p=1} \mathbf{C}\left(\mathbf{\Phi}(p)\right) dp \quad \text{with the path in phase space given by} \quad \mathbf{\Phi}(p) \equiv \hat{\mathbf{U}}_{i+1/2}^{-} + p\left(\hat{\mathbf{U}}_{i+1/2}^{+} - \hat{\mathbf{U}}_{i+1/2}^{-}\right) \tag{3.5}$$

We have found that the convergence of the fixed point iteration in the above paragraph was always very rapid for all the problems that we tested when we started with the first iterate presented above.

Eqns. (3.1) to (3.5) complete the discussion for obtaining $\mathbf{U}^{*}_{i+1/2}$ when an HLL Riemann solver is used with non-conservative products. The fluctuations, as well as the fluxes, from the HLL Riemann solver can then be written in all the supersonic and subsonic limits as

If $\left(S_{L;i+1/2} \geq 0\right)$ then

$$\mathbf{D}^{-}_{HLL;i+1/2} = 0 \quad ; \quad \mathbf{D}^{+}_{HLL;i+1/2} = S_{L;i+1/2}\left(\mathbf{U}^{*}_{i+1/2} - \hat{\mathbf{U}}^{-}_{i+1/2}\right) + S_{R;i+1/2}\left(\hat{\mathbf{U}}^{+}_{i+1/2} - \mathbf{U}^{*}_{i+1/2}\right) \quad ;$$

$$\mathbf{F}^{*}_{HLL;i+1/2} = \mathbf{F}\left(\hat{\mathbf{U}}^{-}_{i+1/2}\right)$$

Else if $\left(S_{R;i+1/2} \leq 0\right)$ then

$$\mathbf{D}^{-}_{HLL;i+1/2} = S_{L;i+1/2}\left(\mathbf{U}^{*}_{i+1/2} - \hat{\mathbf{U}}^{-}_{i+1/2}\right) + S_{R;i+1/2}\left(\hat{\mathbf{U}}^{+}_{i+1/2} - \mathbf{U}^{*}_{i+1/2}\right) \quad ; \quad \mathbf{D}^{+}_{HLL;i+1/2} = 0 \quad ;$$

$$\mathbf{F}^{*}_{HLL;i+1/2} = \mathbf{F}\left(\hat{\mathbf{U}}^{+}_{i+1/2}\right)$$

Else

$$\mathbf{D}^{-}_{HLL;i+1/2} = S_{L;i+1/2}\left(\mathbf{U}^{*}_{i+1/2} - \hat{\mathbf{U}}^{-}_{i+1/2}\right) \quad ; \quad \mathbf{D}^{+}_{HLL;i+1/2} = S_{R;i+1/2}\left(\hat{\mathbf{U}}^{+}_{i+1/2} - \mathbf{U}^{*}_{i+1/2}\right) \quad ;$$

$$\mathbf{F}^{*}_{HLL;i+1/2} = \frac{S_{R;i+1/2}\mathbf{F}\left(\hat{\mathbf{U}}^{-}_{i+1/2}\right) - S_{L;i+1/2}\mathbf{F}\left(\hat{\mathbf{U}}^{+}_{i+1/2}\right) + S_{R;i+1/2}S_{L;i+1/2}\left(\hat{\mathbf{U}}^{+}_{i+1/2} - \hat{\mathbf{U}}^{-}_{i+1/2}\right)}{S_{R;i+1/2} - S_{L;i+1/2}}$$

$$-\frac{1}{S_{R;i+1/2} - S_{L;i+1/2}}\left[S_{R;i+1/2}\tilde{\mathbf{C}}\left(\hat{\mathbf{U}}^{-}_{i+1/2}, \mathbf{U}^{*}_{i+1/2}\right)\left(\mathbf{U}^{*}_{i+1/2} - \hat{\mathbf{U}}^{-}_{i+1/2}\right) + S_{L;i+1/2}\tilde{\mathbf{C}}\left(\mathbf{U}^{*}_{i+1/2}, \hat{\mathbf{U}}^{+}_{i+1/2}\right)\left(\hat{\mathbf{U}}^{+}_{i+1/2} - \mathbf{U}^{*}_{i+1/2}\right)\right]$$

End if

(3.6)

Here, $\mathbf{D}^{-}_{HLL;i+1/2}$ and $\mathbf{D}^{+}_{HLL;i+1/2}$ are the left- and right-going fluctuations from the HLL Riemann solver at the zone boundary $x_{i+1/2}$. The corresponding HLL flux is $\mathbf{F}^{*}_{HLL;i+1/2}$ and is given here for the sake of conceptual completeness. This paragraph completes the discussion for obtaining $\mathbf{U}^{*}_{i+1/2}$ when an HLL or LLF Riemann solver is used with non-conservative products. It also provides all the details for obtaining the fluctuations.

Once $\mathbf{U}^{*}_{i+1/2}$ is available, we can evaluate a diagonal matrix for its eigenvalues, $\mathbf{\Lambda}\left(\mathbf{U}^{*}_{i+1/2}\right)$; we can also evaluate the matrices of left and right eigenvectors $\mathbf{L}\left(\mathbf{U}^{*}_{i+1/2}\right)$ and $\mathbf{R}\left(\mathbf{U}^{*}_{i+1/2}\right)$. We can now project the jump across the Riemann fan into a space of eigenweights, so that we write a vector of eigenweights as $\left[\mathbf{L}\left(\mathbf{U}^{*}_{i+1/2}\right)\cdot\left(\hat{\mathbf{U}}^{+}_{i+1/2} - \hat{\mathbf{U}}^{-}_{i+1/2}\right)\right]$. Depending on its wave speed relative to $S_{L;i+1/2}$ and $S_{R;i+1/2}$, each wave family can take on differing amounts of anti-diffusion, with stationary waves taking on the most anti-diffusive contribution. The amount of this anti-diffusion is calibrated to give the maximum amount of anti-diffusion consistent with maintaining stability of the resulting numerical scheme. The diagonal matrix that governs the anti-diffusion is given by

$$\boldsymbol{\delta}\left(\mathbf{U}^*_{i+1/2}\right)=\mathbf{I}+\frac{\boldsymbol{\Lambda}^-\left(\mathbf{U}^*_{i+1/2}\right)}{S_{i+1/2}}-\frac{\boldsymbol{\Lambda}^+\left(\mathbf{U}^*_{i+1/2}\right)}{S_{i+1/2}} \tag{3.7}$$

Here $\mathbf{I}$ is the identity matrix; $\boldsymbol{\Lambda}^-\left(\mathbf{U}^*_{i+1/2}\right)$ is obtained from $\boldsymbol{\Lambda}\left(\mathbf{U}^*_{i+1/2}\right)$ by zeroing out all its positive eigenvalues ; $\boldsymbol{\Lambda}^+\left(\mathbf{U}^*_{i+1/2}\right)$ is obtained from $\boldsymbol{\Lambda}\left(\mathbf{U}^*_{i+1/2}\right)$ by zeroing out all its negative eigenvalues. With all the entities defined in this paragraph, we can now write out the anti-diffusive flux, $\Phi_{AD;i+1/2}$, that can be added to the numerical flux or the numerical fluctuations. The anti-diffusive flux contribution is given by

$$\Phi_{AD;i+1/2}=-\psi_{i+1/2}\frac{S_{R;i+1/2}S_{L;i+1/2}}{S_{R;i+1/2}-S_{L;i+1/2}}\mathbf{R}\left(\mathbf{U}^*_{i+1/2}\right)\boldsymbol{\delta}\left(\mathbf{U}^*_{i+1/2}\right)\left[\mathbf{L}\left(\mathbf{U}^*_{i+1/2}\right)\cdot\left(\hat{\mathbf{U}}^+_{i+1/2}-\hat{\mathbf{U}}^-_{i+1/2}\right)\right] \tag{3.8}$$

Note that $\psi_{i+1/2}$ is a shock detector which is 0 in the presence of shocks and 1 when the solution is smooth (see Balsara [5]). This ensures that there is no over-steepening from the anti-diffusive terms when a shock is present. Notice that the anti-diffusive term in eqn. (3.8) is such that if the speeds $S_{L;i+1/2}$ and $S_{R;i+1/2}$ are slightly modified, the anti-diffusive contribution adjusts accordingly. This ensures that the contribution from the anti-diffusive terms is exactly right whether we choose an HLL formulation or an LLF formulation. (The HLL Riemann solver, with the anti-diffusive contribution, is referred to as the HLLI Riemann solver; where the "I" stands for intermediate waves.) The anti-diffusive contribution is only needed when the Riemann fan is subsonic. The subsonic fluctuations, as well as the subsonic HLL flux, in eqn. (3.6) then get modified as

$$\mathbf{D}^-_{HLLI;i+1/2}=S_{L;i+1/2}\left(\mathbf{U}^*_{i+1/2}-\hat{\mathbf{U}}^-_{i+1/2}\right)+\Phi_{AD;i+1/2} \quad ; \quad \mathbf{D}^+_{HLLI;i+1/2}=S_{R;i+1/2}\left(\hat{\mathbf{U}}^+_{i+1/2}-\mathbf{U}^*_{i+1/2}\right)-\Phi_{AD;i+1/2} \quad ;$$
$$\mathbf{F}^*_{HLLI;i+1/2}=\mathbf{F}^*_{HLL;i+1/2}+\Phi_{AD;i+1/2} \tag{3.9}$$

Observe that the anti-diffusive contribution is fully conservative in the sense that the contribution that is added to $\mathbf{D}^{-}_{HLLI;i+1/2}$ is subtracted from $\mathbf{D}^{+}_{HLLI;i+1/2}$. This completes our description of the anti-diffusive flux in eqns. (3.8) and (3.9). It should be added to the numerical fluxes in order to preserve stationary contact discontinuities, or stationary wave structures from linearly degenerate wave fields in general.

It is worthwhile making another observation about eqn. (3.7) that contributes to its usefulness. For many PDE systems the linearly degenerate right and left eigenvectors are quite easy to evaluate, whereas the genuinely non-linear right and left eigenvectors are very difficult to evaluate. For some systems (and the debris flow model of Pitman and Le is one such prominent example) the genuinely non-linear eigenvectors can only be evaluated numerically and there is no known analytical way to evaluate them. Because our goal is to introduce anti-diffusion only for the linearly degenerate wave fields, it is acceptable to only retain those eigenvectors in $\mathbf{L}\left(\mathbf{U}^{*}_{i+1/2}\right)$ and $\mathbf{R}\left(\mathbf{U}^{*}_{i+1/2}\right)$; the other rows of $\mathbf{L}\left(\mathbf{U}^{*}_{i+1/2}\right)$ and the other columns of $\mathbf{R}\left(\mathbf{U}^{*}_{i+1/2}\right)$ can indeed be set to zero. Indeed, for all the results shown in this work, we altogether avoided the evaluation and use of eigenvectors corresponding to genuinely non-linear wave fields in eqn. (3.7).

**IV) An Efficient FD-WENO for Conservative and Non-Conservative Systems**

We are now ready to specify our FD-WENO formulation for conservative and non-conservative systems. For the one-dimensional PDE shown in eqns. (2.4) and (2.5) we retain the same WENO construction that was shown in Fig. 1 and described in detail in the second paragraph of Section II. The only difference is that now we invoke the Riemann solver at each zone boundary. Consequently, at each zone boundary we can have two wave speeds – a left-going wave speed $S_{L;i+1/2}$ and a right-going wave speed $S_{R;i+1/2}$. Our strategy is to evaluate these wave speeds at each and every zone boundary first. Using this information, we reset $S_{L;i+1/2}$ to be the minimum value from all the stencils that could contribute to the zone boundary $x_{i+1/2}$. Similarly, we reset $S_{R;i+1/2}$ to be the maximum value from all the stencils that could contribute to the zone boundary $x_{i+1/2}$. This is needed in order to obtain the order property, which would be damaged if the extremal signal speeds varied too much over the stencils of interest. This choice of extremizing the left and right wave speeds is acceptable because the Riemann solver that is described in Section III is very forgiving to small modifications of the bounding wave speeds. Furthermore, it is also acceptable because the anti-diffusive terms can adjust as $S_{L;i+1/2}$ and $S_{R;i+1/2}$ are slightly modified.

With this very slight modification to the extremal signal speeds, at each zone boundary, say $x_{i+1/2}$, the higher order WENO reconstruction gives us a left state $\hat{\mathbf{U}}^{-}_{i+1/2}$ and a right state $\hat{\mathbf{U}}^{+}_{i+1/2}$. These left and right states are then used in eqn. (3.1) to obtain the single resolved state within the Riemann fan, given by $\mathbf{U}^{*}_{i+1/2}$. Eqn. (3.6) can then be used to obtain the left-going fluctuation $\mathbf{D}^{-}_{HLL;i+1/2}$ and the right-going fluctuation $\mathbf{D}^{+}_{HLL;i+1/2}$ if we want to use the HLL Riemann solver without anti-diffusion. Alternatively, if we wish to include anti-diffusion, eqn. (3.9) can then be used to obtain the left-going fluctuation $\mathbf{D}^{-}_{HLLI;i+1/2}$ and the right-going fluctuation $\mathbf{D}^{+}_{HLLI;i+1/2}$. With those fluctuations in hand, the discrete in space but continuous in time scheme can be written very compactly as

$$\partial_t \mathbf{U}_i = -\mathbf{A}_i \frac{\left(\hat{\mathbf{U}}^{-}_{i+1/2} - \hat{\mathbf{U}}^{+}_{i-1/2}\right)}{\Delta x} - \frac{1}{\Delta x}\left(\mathbf{D}^{-}_{HLLI;i+1/2} + \mathbf{D}^{+}_{HLLI;i-1/2}\right) \tag{4.1}$$

This completes our description of the higher order FD-WENO scheme. Like all FD-WENO formulations, the extension to multidimensions is trivial because the scheme is a finite difference scheme. In Section V we will show that it reduces transparently to eqn. (2.3) in the limit where the characteristic matrix is a constant. As a result, we may expect eqn. (4.1) to also retain high order accuracy on multidimensional problems. Because of our use of anti-diffusive fluctuations, we also expect the scheme to be able to preserve stationary discontinuities on the computational mesh. In Section VI we give a pointwise description of its implementation.

## V) Equivalence of Our FD-WENO Scheme to the Linear FD-WENO Scheme in Section II

It is easy to retrieve the LLF limit for the formulation in the previous Section because one has simply to define $S \equiv \max\left(\left|S_{L;i+1/2}\right|, \left|S_{R;i+1/2}\right|\right)$ and then reset $S_{R;i+1/2} \to S$ and $S_{L;i+1/2} \to -S$. In this Section we demonstrate the equivalence of eqn. (4.1) to eqn. (2.3) when the characteristic matrix is a constant and the LLF limit of the HLL Riemann solver is used.

We start with eqn. (2.3) and write it as

$$\partial_t \mathbf{U}_i = -\mathbf{A}\frac{\left(\hat{\mathbf{U}}^-_{i+1/2} - \hat{\mathbf{U}}^+_{i-1/2}\right)}{\Delta x} - \mathbf{A}^+\frac{\left(\hat{\mathbf{U}}^-_{i+1/2} - \hat{\mathbf{U}}^-_{i-1/2}\right)}{\Delta x} - \mathbf{A}^-\frac{\left(\hat{\mathbf{U}}^+_{i+1/2} - \hat{\mathbf{U}}^+_{i-1/2}\right)}{\Delta x} + \mathbf{A}\frac{\left(\hat{\mathbf{U}}^-_{i+1/2} - \hat{\mathbf{U}}^+_{i-1/2}\right)}{\Delta x}, \quad (5.1)$$

where we realize that the same term has been added and subtracted to eqn. (2.3). We can now write eqn. (5.1) as

$$\begin{aligned}\partial_t \mathbf{U}_i = &-\mathbf{A}\frac{\left(\hat{\mathbf{U}}^-_{i+1/2} - \hat{\mathbf{U}}^+_{i-1/2}\right)}{\Delta x} \\ &-\frac{1}{\Delta x}\left\{\left[\mathbf{A}^+\hat{\mathbf{U}}^-_{i+1/2} + \mathbf{A}^-\hat{\mathbf{U}}^+_{i+1/2} - \mathbf{A}\hat{\mathbf{U}}^-_{i+1/2}\right] + \left[-\mathbf{A}^+\hat{\mathbf{U}}^-_{i-1/2} - \mathbf{A}^-\hat{\mathbf{U}}^+_{i-1/2} + \mathbf{A}\hat{\mathbf{U}}^+_{i-1/2}\right]\right\}\end{aligned} \quad (5.2)$$

In the LLF limit, where $S_{R;i+1/2} \to S$ and $S_{L;i+1/2} \to -S$, it is easy to use the equations from Section III to show that

$$\begin{aligned}&\mathbf{U}^*_{i+1/2} = \frac{1}{S}\left[\mathbf{A}^+\hat{\mathbf{U}}^-_{i+1/2} - \mathbf{A}^-\hat{\mathbf{U}}^+_{i+1/2}\right] \quad ; \quad \mathbf{U}^*_{i-1/2} = \frac{1}{S}\left[\mathbf{A}^+\hat{\mathbf{U}}^-_{i-1/2} - \mathbf{A}^-\hat{\mathbf{U}}^+_{i-1/2}\right] \\ &\mathbf{D}^-_{LLF;i+1/2} = S\ \hat{\mathbf{U}}^-_{i+1/2} - \mathbf{A}^+\hat{\mathbf{U}}^-_{i+1/2} + \mathbf{A}^-\hat{\mathbf{U}}^+_{i+1/2} \quad ; \quad \mathbf{D}^+_{LLF;i-1/2} = S\ \hat{\mathbf{U}}^+_{i-1/2} - \mathbf{A}^+\hat{\mathbf{U}}^-_{i-1/2} + \mathbf{A}^-\hat{\mathbf{U}}^+_{i-1/2}\end{aligned} \quad (5.3)$$

Now we can see quite easily that

$$\begin{aligned}\mathbf{D}^-_{LLF;i+1/2} &= S\ \hat{\mathbf{U}}^-_{i+1/2} - \mathbf{A}^+\hat{\mathbf{U}}^-_{i+1/2} + \mathbf{A}^-\hat{\mathbf{U}}^+_{i+1/2} = S\ \mathbf{I}\ \hat{\mathbf{U}}^-_{i+1/2} - \frac{1}{2}\left(\mathbf{A} + S\ \mathbf{I}\right)\hat{\mathbf{U}}^-_{i+1/2} + \mathbf{A}^-\hat{\mathbf{U}}^+_{i+1/2} \\ &= -\mathbf{A}\hat{\mathbf{U}}^-_{i+1/2} + \frac{1}{2}\left(\mathbf{A} + S\ \mathbf{I}\right)\hat{\mathbf{U}}^-_{i+1/2} + \mathbf{A}^-\hat{\mathbf{U}}^+_{i+1/2} = \mathbf{A}^+\hat{\mathbf{U}}^-_{i+1/2} + \mathbf{A}^-\hat{\mathbf{U}}^+_{i+1/2} - \mathbf{A}\hat{\mathbf{U}}^-_{i+1/2}\end{aligned} \quad (5.4)$$

This shows that the first square bracket in eqn. (5.2) is indeed $\mathbf{D}^-_{LLF;i+1/2}$. We can also see that

$$\begin{aligned}\mathbf{D}^+_{LLF;i-1/2} &= S\ \hat{\mathbf{U}}^+_{i-1/2} - \mathbf{A}^+\hat{\mathbf{U}}^-_{i-1/2} + \mathbf{A}^-\hat{\mathbf{U}}^+_{i-1/2} = S\ \mathbf{I}\ \hat{\mathbf{U}}^+_{i-1/2} + \frac{1}{2}\left(\mathbf{A} - S\ \mathbf{I}\right)\hat{\mathbf{U}}^+_{i-1/2} - \mathbf{A}^+\hat{\mathbf{U}}^-_{i-1/2} \\ &= \mathbf{A}\hat{\mathbf{U}}^+_{i-1/2} - \frac{1}{2}\left(\mathbf{A} - S\ \mathbf{I}\right)\hat{\mathbf{U}}^+_{i-1/2} - \mathbf{A}^+\hat{\mathbf{U}}^-_{i-1/2} = -\mathbf{A}^+\hat{\mathbf{U}}^-_{i-1/2} - \mathbf{A}^-\hat{\mathbf{U}}^+_{i-1/2} + \mathbf{A}\hat{\mathbf{U}}^+_{i-1/2}\end{aligned} \quad (5.5)$$

This shows that the second square bracket in eqn. (5.2) is indeed $\mathbf{D}^+_{LLF;i-1/2}$. Consequently, by examining eqn. (4.2) along with eqns. (5.2) to (5.5), we see that the equivalence that we sought to establish has indeed been established!

Note, however, that the above-mentioned equivalence is only true for a linear PDE in the limit when an LLF Riemann solver is used. Eqn. (4.1) is based on including all the non-linearities of the PDE in the HLL (or LLF) Riemann solver. As a result, it is accurate in the presence of non-linearities; at least when the flow is smooth. When used with an HLL Riemann solver, eqn. (4.1) can reach the correct supersonic limit. By virtue of being in fluctuation form, eqn. (4.1) also includes the contribution from non-conservative products. By virtue of including anti-diffusive contributions, eqn. (4.1) can also preserve isolated stationary discontinuities. We see, therefore, that our build-up to eqn. (4.1) has given us a very useful formulation of FD-WENO that significantly improves on eqns. (2.3) or (2.6).

We would like to make a further observation about this Section. Eqn. (4.1) is not novel. Indeed, it occurs, with slight modification, as eqn. (32) in DB16 where it is presented as a second order scheme! So in what sense is eqn. (4.1) novel? It is novel in the sense that the analysis in this section has shown that it is also legitimately higher order accurate (in a finite difference sense) for smooth flow! It is also novel because we use eqn. (4.1) to demonstrate the higher order accuracy of our FD-WENO scheme for a range of attractive applications areas.

## VI) Pointwise Implementation of Our FD-WENO Scheme

Now that the above discussions are understood, we provide a pointwise implementation of our FD-WENO scheme for treating non-linear PDEs with non-conservative products. The scheme also works when non-conservative products are absent. It indeed presents one significant improvement over all past FD-WENO formulations in that it can represent stationary discontinuities exactly on the computational mesh. The pointwise implementation of our FD-WENO scheme goes according to the following steps:-

**1)** We start with the mesh function as shown in Fig. 1. This means that at each zone center $x_i$ we have a pointwise value for the primal variable $\mathbf{U}_i$. FD-WENO schemes are always implemented in dimension-by-dimension fashion, so we only describe one of the dimensional updates here.

**2)** From the conserved variables in each zone, obtain the primitive variables. Use them both to obtain the normalized right and left eigenvectors in the conserved variables.

**3)** As shown in Fig. 1, we use the WENO-AO algorithm from Balsara, Garain and Shu [6]. That paper includes all closed form expressions that are needed for WENO reconstruction in one dimension. This consists of making a non-linear hybridization between a large high order accurate stencil and smaller lower order accurate stencils. We use WENO-AO reconstruction in characteristic variables. As a result, the neighboring zones around zone "*i*" are projected into the characteristic space of zone "*i*". The extent of these neighboring zones depends on the desired order of the scheme. The fifth order case is explicitly shown in Fig. 1. Expanding the large stencil by one zone on either side adds two further orders of accuracy. Once the variables in the neighboring zones around zone "*i*" are projected into the characteristic space of zone "*i*", WENO reconstruction is carried out in the characteristic space. Projecting the reconstructed characteristic variables back into the space of right eigenvectors gives us high order accurate $\hat{\mathbf{U}}^{-}_{i+1/2}$ and $\hat{\mathbf{U}}^{+}_{i-1/2}$ within each zone "*i*", as shown in Fig. 1.

**4)** At each zone boundary $x_{i+1/2}$, use the left and right states $\hat{\mathbf{U}}^{-}_{i+1/2}$ and $\hat{\mathbf{U}}^{+}_{i+1/2}$ to obtain the left-most and right-most going speeds of the Riemann fan; these are denoted by $S_{L;i+1/2}$ and $S_{R;i+1/2}$. Please note that at this point in the game we are not yet seeking the resolved state within the Riemann fan; that will come later.

**5)** Once the left-most and right-most going speeds have been obtained for any zone boundary, say $x_{i+1/2}$, we examine all the stencils that will contribute to that zone boundary. We use the extremal speeds from all the zone boundaries that are present in those stencils to extend $S_{L;i+1/2}$ to the left and $S_{R;i+1/2}$ to the right. This gives us the largest possible extremal speeds that will be needed at the zone boundary in question.

**6)** Now, at each zone boundary $x_{i+1/2}$, we hand in the speeds $S_{L;i+1/2}$ and $S_{R;i+1/2}$ as well as the states $\hat{\mathbf{U}}^{-}_{i+1/2}$ and $\hat{\mathbf{U}}^{+}_{i+1/2}$ to the Riemann solver. This gives us the resolved state $\mathbf{U}^{*}_{i+1/2}$ and the left and right-going fluctuations $\mathbf{D}^{-}_{HLLI;i+1/2}$ and $\mathbf{D}^{+}_{HLLI;i+1/2}$ respectively. Section III provides all the formulae for doing this. (While this is optional, if one can design a flattener function for the PDE (see Balsara [5] for an example), one can avoid the introduction of anti-diffusive terms at shocks.)

**7)** At each zone center, obtain $\mathbf{A}_i = \mathbf{A}\left(\mathbf{U}_i\right)$. We now have everything we need for evaluating the right hand side of eqn. (4.1). For multidimensional problems, sum the contributions from each dimension. This gives us a single stage for a multi-stage SSP Runge-Kutta scheme. Imposing the timestep gives us the full scheme.

Notice that this scheme entails only one WENO reconstruction, whereas the LLF scheme in Shu and Osher [55] or Jiang and Shu [40] entails two WENO reconstructions. If the WENO reconstruction should dominate the costs in an application, then the FD-WENO scheme presented here is half as computationally expensive as prior generations of FD-WENO schemes. The above points have only shown one stage of the scheme. It can be coupled with an SSP-RK update strategy, say from Shu and Osher [55] or Spiteri and Ruuth [59], [60], to achieve higher order in time.

Many of the PDEs with non-conservative products also have stiff source terms; these are usually relaxation terms that enable the system to relax to several useful physical limits. The FD-WENO method makes it very simple to treat stiff source terms because the source terms are treated pointwise and are collocated at the exact same location as the primal variables. For this reason, when stiff source terms are present, we recommend using the Runge-Kutta IMEX methods from Pareschi and Russo [45]; see also Kupka *et al*. [39]

## VII) Accuracy Study

We present several accuracy studies for Euler flow, to show that the method even works for conservative systems. We also present accuracy studies of Baer-Nunziato compressible multi-phase flow, multiphase debris flow and two-layer shallow water flow, as an examples of non-conservative systems. This shows the reader the generality of our approach.

For the timestepping we used a third order SSP-RK scheme from Shu and Osher [46] and a fourth order scheme from Spiteri and Ruuth [59], [60] with the result that for the fifth and seventh order schemes we had to reduce the timestep as the mesh was refined. This allows the temporal build up of error to be contained at levels that remain below the error build up from the spatial accuracy of the scheme. The base level mesh in all of these accuracy tests was run with a CFL of 0.3 for all our two-dimensional tests. When 3rd order timestepping is used with an $r^{th}$ order spatial scheme then one should have $\left(\Delta t_{fine}/\Delta t_{coarse}\right)=\left(\Delta x_{fine}/\Delta x_{coarse}\right)^{r/3}$ . When 4th order timestepping is used with an $r^{th}$ order spatial scheme then one should have $\left(\Delta t_{fine}/\Delta t_{coarse}\right)=\left(\Delta x_{fine}/\Delta x_{coarse}\right)^{r/4}$ . Here $\Delta t_{fine}$ and $\Delta t_{coarse}$ are timesteps on the fine and coarse mesh respectively; and $\Delta x_{fine}$ and $\Delta x_{coarse}$ are zone sizes on the fine and coarse mesh respectively.

**VII.a) Accuracy Study for Two Dimensional Euler Flow**

As our first test, we show the accuracy study of a hydrodynamic vortex when it is propagated diagonally on a square mesh. Our purpose in doing this is to show the reader that the same formulation works just as well for conservation laws as for PDEs with non-conservative products. The set-up of the vortex is described in Pao and Salas [44] and this test problem has been described in Balsara and Shu [3]. Table I shows that the present FD-WENO scheme works well for multidimensional conservation laws. Since seventh and ninth order schemes reach machine precision very rapidly, the domain was doubled to $\left[-10,10\right]\times\left[-10,10\right]$ and the stopping time was also doubled, but the rest of the problem was kept the same. We show the results from using WENO-AO reconstruction (Balsara, Garain and Shu [6]) as well as Multiresolution WENO reconstruction (Zhu and Shu [64]). For both variants of WENO reconstruction, the method does achieve the desired order of accuracy. For WENO-AO we had to change the "$\tau$" term that is used in the smoothness indicators so that it scales with the order of accuracy of the scheme. For Multiresolution WENO, we got less than second order accuracy with the published third order accurate version and had to modify it in order to restore third order of accuracy. The requisite modifications are documented in Appendix D.

**Table I shows the accuracy of the 2D Euler vortex problem using the LLF Riemann solver; the density variable is shown. The first half of the table shows the results from WENO-AO with non-linear limiting. The second half shows the same results when Multiresolution WENO reconstruction is used. Here the lowest order polynomial in the Multiresolution WENO was modified to have piecewise linear variation with the $MC_\beta$ limiter.**

| WENO-AO-3 | $L_1$ Error | $L_1$ Accuracy | $L_{inf}$ Error | $L_{inf}$ Accuracy |
|---|---|---|---|---|
| $64^2$ | 3.77899E-03 | | 6.79759E-02 | |
| $128^2$ | 4.83119E-04 | 2.97 | 1.17871E-02 | 2.53 |
| $256^2$ | 6.37053E-05 | 2.92 | 2.12058E-03 | 2.47 |
| $512^2$ | 1.05943E-05 | 2.59 | 4.31110E-04 | 2.30 |
| WENO-AO-(5,3) | | | | |
| $64^2$ | 4.03476E-04 | | 2.73665E-02 | |
| $128^2$ | 4.02117E-05 | 3.33 | 3.87474E-03 | 2.82 |
| $256^2$ | 1.67993E-06 | 4.58 | 1.05470E-04 | 5.20 |
| $512^2$ | 5.37313E-08 | 4.97 | 3.54283E-06 | 4.90 |
| WENO-AO-(7,3) | | | | |
| $64^2$ | 2.09416E-04 | | 2.95412E-02 | |
| $128^2$ | 5.71233E-06 | 5.20 | 5.41615E-04 | 5.77 |
| $256^2$ | 6.31543E-08 | 6.50 | 7.56964E-06 | 6.16 |
| $512^2$ | 5.18406E-10 | 6.93 | 5.93836E-08 | 6.99 |
| WENO-AO-(9,3) | | | | |
| $64^2$ | 9.17405E-05 | | 9.93215E-03 | |
| $128^2$ | 1.23503E-06 | 6.21 | 2.23517E-04 | 5.47 |
| $256^2$ | 3.56071E-09 | 8.44 | 3.64506E-07 | 9.26 |
| $384^2$ | 1.01355E-10 | 8.78 | 1.03245E-08 | 8.79 |
| | | | | |
| Order 3 Multires WENO | | | | |
| $64^2$ | 2.31330E-03 | | 4.11024E-02 | |
| $128^2$ | 2.82763E-04 | 3.03 | 6.74873E-03 | 2.61 |
| $256^2$ | 3.37308E-05 | 3.07 | 1.11187E-03 | 2.60 |
| $512^2$ | 5.06473E-06 | 2.74 | 1.96406E-04 | 2.50 |
| Order 5 Multires WENO | | | | |
| $64^2$ | 5.76628E-04 | | 3.36392E-02 | |
| $128^2$ | 5.34567E-05 | 3.43 | 4.02233E-03 | 3.06 |
| $256^2$ | 1.94283E-06 | 4.78 | 1.31879E-04 | 4.93 |
| $512^2$ | 7.56482E-08 | 4.68 | 4.53360E-06 | 4.86 |
| Order 7 Multires WENO | | | | |
| $64^2$ | 2.19228E-04 | | 2.87452E-02 | |
| $128^2$ | 5.87231E-06 | 5.22 | 5.44809E-04 | 5.72 |
| $256^2$ | 6.31596E-08 | 6.54 | 7.56856E-06 | 6.17 |
| $512^2$ | 5.18407E-10 | 6.93 | 5.93836E-08 | 6.99 |
| Order 9 Multires WENO | | | | |
| $64^2$ | 7.89925E-05 | | 9.67640E-03 | |
| $128^2$ | 1.23503E-06 | 6.00 | 2.23518E-04 | 5.44 |
| $256^2$ | 3.56071E-09 | 8.44 | 3.64506E-07 | 9.26 |
| $384^2$ | 1.01352E-10 | 8.78 | 1.03249E-08 | 8.79 |

**VII.b) Accuracy Study for Two Dimensional Baer-Nunziato Model for Compressible Multi-Phase Flows**

Now let us move on to a PDE that has non-conservative products. We focus on the Baer-Nunziato compressible multiphase flow proposed by Baer and Nunziato [2] and extensively studied at the numerical level by Saurel and Abgrall [52], Adrianov and Warnecke [1], Schwendeman *et al*. [53], Dumbser *et al*. [25], Tokareva and Toro [61], Coquel et al. [17] and Chiochetti and Müller [16]. A very useful set of eigenvectors have been presented in Tokareva and Toro [61].

The PDE system assumes two phases, a solid phase denoted by density $\rho_1$, volume fraction $\phi_1$, velocity $\mathbf{v}_1 = (u_1, v_1, w_1)$ and a pressure $p_1$ and a gas phase denoted by density $\rho_2$, volume fraction $\phi_2$, velocity $\mathbf{v}_2 = (u_2, v_2, w_2)$ and a pressure $p_2$. The two phases have an interfacial pressure $P_I$ and an interfacial velocity $\mathbf{V}_I$, but it is suggested in Baer and Nunziato [2] to set $P_I = p_2$ and $\mathbf{V}_I = \mathbf{v}_1$. The total energy density for phase "*j*" is related to the specific internal energy $e_j$ by $\rho_j E_j = \rho_j e_j + \rho_j \mathbf{v}_j^2 / 2$.

$$
\begin{aligned}
&\partial_t(\phi_1\rho_1) + \nabla\cdot(\phi_1\rho_1\mathbf{v}_1) = 0\\
&\partial_t(\phi_1\rho_1\mathbf{v}_1) + \nabla\cdot(\phi_1(\rho_1\mathbf{v}_1\otimes\mathbf{v}_1 + \mathbf{I}\,p_1)) - P_I\nabla\phi_1 = 0\\
&\partial_t(\phi_1\rho_1E_1) + \nabla\cdot(\phi_1\mathbf{v}_1(\rho_1E_1 + p_1)) + P_I\partial_t\phi_1 = 0\\
&\partial_t(\phi_2\rho_2) + \nabla\cdot(\phi_2\rho_2\mathbf{v}_2) = 0\\
&\partial_t(\phi_2\rho_2\mathbf{v}_2) + \nabla\cdot(\phi_2(\rho_2\mathbf{v}_2\otimes\mathbf{v}_2 + \mathbf{I}\,p_2)) - P_I\nabla\phi_2 = 0\\
&\partial_t(\phi_2\rho_2E_2) + \nabla\cdot(\phi_2\mathbf{v}_2(\rho_2E_2 + p_2)) + P_I\partial_t\phi_2 = 0\\
&\partial_t\phi_1 + \mathbf{V}_I\cdot\nabla\phi_1 = 0
\end{aligned}
$$

The system requires that the phases volume fractions add up to unity, $\phi_1 + \phi_2 = 1$. The closure relations for each phase are also given by

$$
\rho_j e_j = \frac{p_j + \gamma_j \pi_j}{\gamma_j - 1}
$$

Here $\gamma_j$ is the ratio of specific heats and $\pi_j$ is a constant. For the above EOS, the sound speed $c_j$ in each phase is given by

$$c_j = \sqrt{\gamma_j \frac{p_j + \pi_j}{\rho_j}}$$

In Section 7.1 of Coquel *et al*. [17] a 2D smooth test problem that is useful for accuracy analysis was presented for Baer-Nunziato flow. The problem consists of a diagonal flow of sinusoidal perturbations in the density. We show the accuracies from that test problems here. Table II shows the 2D test from that paper, again showing that we reach the design accuracies.

In Section 4 of Dumbser *et al*. [22] a 2D smooth vortex test problem was designed for Baer-Nunziato flow. The problem is an analogue of the hydrodynamic vortex from the previous sub-section. We would simply like to fix a typographical error in their eqn. (43) which should actually read as

$$u_\theta^1 = \frac{r}{2s_1 D}\sqrt{D\left[p_{10}\left(4\sqrt{2\pi}F_1 + 6H_1 - 12Gs_1^2 + 3H_1 s_1^2\right) + 3p_{20}s_1^2\left(4G - H_2\right)\right]}$$

$$u_\theta^2 = \frac{r\sqrt{2}}{2\rho_2 s_2}\sqrt{\rho_2 p_{20} F_2}$$

For fifth and seventh orders, we also had to double size of the computational domain in order to obtain results that were independent of the exponential fall off in the velocity of the vortex. Appendix A of this paper fully documents this vortex problem for the Baer-Nunziato equations for the sake of completeness. Table III shows the 2D test from that paper. In Table III we show the accuracy with the non-linear limiter and also the same accuracy when only the central stencil is used. As before, we see that the third and fifth order schemes reach their design accuracies very efficiently in all circumstances. The seventh order scheme shows a slight degradation in accuracy when the non-linear limiting is turned on. But Table III shows us that it too reaches its design accuracy when the contribution from the limiter is suppressed.

**Table II shows the accuracy of the 2D Baer-Nunziato problem using the LLF Riemann solver; the average of the solid and gas densities is shown.**

| WENO-AO-3 | $L_1$ Error | $L_1$ Accuracy | $L_{inf}$ Error | $L_{inf}$ Accuracy |
|---|---|---|---|---|
| $64^2$ | 2.65678E-02 | | 8.47012E-02 | |
| $128^2$ | 9.33135E-03 | 1.51 | 2.32778E-02 | 1.86 |
| $256^2$ | 1.91409E-03 | 2.29 | 8.17852E-03 | 1.51 |
| $512^2$ | 3.05585E-04 | 2.65 | 2.49099E-03 | 1.72 |
| WENO-AO-(5,3) | | | | |
| $64^2$ | 1.40573E-02 | | 3.00121E-02 | |
| $128^2$ | 1.97121E-03 | 2.83 | 5.21534E-03 | 2.52 |
| $256^2$ | 6.47488E-05 | 4.93 | 1.75707E-04 | 4.89 |
| WENO-AO-(7,3) | | | | |
| $64^2$ | 2.11856E-05 | | 5.83654E-05 | |
| $128^2$ | 1.64576E-07 | 7.01 | 4.55047E-07 | 7.00 |
| $256^2$ | 1.27134E-09 | 7.02 | 3.51204E-09 | 7.02 |
| WENO-AO-(9,3) | | | | |
| $64^2$ | 4.29813E-07 | | 1.22794E-06 | |
| $128^2$ | 8.22844E-10 | 9.03 | 2.29503E-09 | 9.06 |
| $192^2$ | 2.36162E-11 | 8.76 | 7.14733E-11 | 8.56 |

**Table III shows the accuracy of the 2D Baer-Nunziato Vortex problem using the LLF Riemann solver; the solid volume fraction is shown. The first half of the table shows the results from WENO-AO with non-linear limiting. The second half shows the same results when Multiresolution WENO reconstruction is used. Here the lowest order polynomial in the Multiresolution WENO was modified to have piecewise linear variation with the $MC_\beta$ limiter.**

| WENO-AO-3 | $L_1$ Error | $L_1$ Accuracy | $L_{inf}$ Error | $L_{inf}$ Accuracy |
|---|---|---|---|---|
| $64^2$ | 3.65555E-04 | | 2.26149E-02 | |
| $128^2$ | 6.74918E-05 | 2.44 | 8.13744E-03 | 1.47 |
| $256^2$ | 1.19388E-05 | 2.50 | 2.80208E-03 | 1.54 |
| $512^2$ | 1.90516E-06 | 2.65 | 9.48405E-04 | 1.56 |
| WENO-AO-(5,3) | | | | |
| $32^2$ | 4.42755E-04 | | 3.54580E-02 | |
| $64^2$ | 5.75719E-05 | 2.94 | 4.72247E-03 | 2.91 |
| $128^2$ | 1.97053E-06 | 4.87 | 2.29502E-04 | 4.36 |
| $256^2$ | 6.20766E-08 | 4.99 | 7.72485E-06 | 4.89 |
| WENO-AO-(7,3) | | | | |
| $32^2$ | 2.46277E-04 | | 8.87359E-03 | |
| $64^2$ | 7.37154E-06 | 5.06 | 6.29709E-04 | 3.82 |
| $128^2$ | 6.43636E-08 | 6.84 | 8.50361E-06 | 6.21 |
| $256^2$ | 5.27973E-10 | 6.93 | 7.07152E-08 | 6.91 |
| WENO-AO-(9,3) | | | | |
| $32^2$ | 3.26754E-04 | | 1.45547E-02 | |
| $64^2$ | 1.18362E-06 | 8.11 | 1.02302E-04 | 7.15 |
| $128^2$ | 3.14880E-09 | 8.55 | 4.13505E-07 | 7.95 |
| $256^2$ | 8.25669E-12 | 8.58 | 8.63329E-10 | 8.90 |
| | | | | |
| Order 3 Multires WENO | | | | |
| $64^2$ | 2.53132E-04 | | 1.60473E-02 | |
| $128^2$ | 4.20954E-05 | 2.59 | 5.39172E-03 | 1.57 |
| $256^2$ | 7.36060E-06 | 2.52 | 1.82224E-03 | 1.57 |
| $512^2$ | 1.13418E-06 | 2.70 | 6.14775E-04 | 1.57 |
| Order 5 Multires WENO | | | | |
| $32^2$ | 8.53261E-04 | | 3.88555E-02 | |
| $64^2$ | 7.84887E-05 | 3.44 | 9.43636E-03 | 2.04 |
| $128^2$ | 4.55467E-06 | 4.11 | 6.38861E-04 | 3.88 |
| $256^2$ | 1.63743E-07 | 4.80 | 5.43007E-05 | 3.56 |
| Order 7 Multires WENO | | | | |
| $32^2$ | 4.21351E-04 | | 1.50146E-02 | |
| $64^2$ | 7.67412E-06 | 5.78 | 6.89414E-04 | 4.44 |
| $128^2$ | 9.44628E-08 | 6.34 | 1.38913E-05 | 5.63 |
| $256^2$ | 5.40514E-10 | 7.45 | 7.20604E-08 | 7.59 |
| Order 9 Multires WENO | | | | |
| $32^2$ | 2.57891E-04 | | 1.00495E-02 | |
| $64^2$ | 1.32054E-05 | 4.29 | 1.19933E-03 | 3.07 |
| $128^2$ | 3.09693E-09 | 12.06 | 3.97746E-07 | 11.56 |
| $256^2$ | 8.25749E-12 | 8.55 | 8.63332E-10 | 8.85 |

To complete our description of the Baer-Nunziato model, we also mention that the model has seen extensive, multifaceted development in the literature, which we cannot address here in a numerical paper. Kapila *et al*. [37], realized the need for viscous regularization, which was addressed further by Delchini *et al*. [19]. Powers [50] mentioned the need for resolving the thin layers down to a physical scale and this was discussed further by Bdzil *et al*. [9] and Hennessey *et al*. [33]. It was also deemed possible to reformulate the non-conservative terms as conservative terms, see Gonthier *et al*. [29] [30].

**VII.c) Accuracy Study for Two Dimensional Multiphase Debris Flow Model of Pitman and Le**

The multiphase debris flow model of Pitman and Le [49] is described here. We use the formulation of Pelanti *et al*. [46] instead of the original formulation by Pitman and Le [49]. The PDE system in 2D is given by defining $h_s$, $u_s$, $v_s$ as the solid height, the solid x-velocity and the solid y-velocity respectively and by defining $h_f$, $u_f$, $v_f$ as the fluid height, the fluid x-velocity and the solid y-velocity respectively. We also have $h_f = (1-\phi)h$ where $\phi$ is the solid volume fraction and $h$ is the total height. The variable "*b*" refers to the bottom topography and is kept constant in all our test problems. The variable "g" refers to the gravitational acceleration and "$\rho \equiv \rho_f / \rho_s$" refers to the density ratio of the fluid and the solid. The time evolutionary equations for this model are:-

$$\frac{\partial}{\partial t}\begin{pmatrix} h_s \\ h_s u_s \\ h_s v_s \\ h_f \\ h_f u_f \\ h_f v_f \\ b \end{pmatrix} + \frac{\partial}{\partial x}\begin{pmatrix} h_s u_s \\ h_s u_s^2 + g h_s^2 \rho / 2 + g(1-\rho) h_s (h_f + h_s)/2 \\ h_s u_s v_s \\ h_f u_f \\ h_f u_f^2 + g h_f^2 / 2 \\ h_f u_f v_f \\ 0 \end{pmatrix}$$

$$+\frac{\partial}{\partial y}\begin{pmatrix} h_s v_s \\ h_s u_s v_s \\ h_s v_s^2 + g h_s^2 \rho / 2 + g(1-\rho) h_s (h_f + h_s)/2 \\ h_f v_f \\ h_f u_f v_f \\ h_f v_f^2 + g h_f^2 / 2 \\ 0 \end{pmatrix} + \begin{pmatrix} 0 \\ g\ h_s \rho\ \partial_x h_f + g\ h_s \partial_x b \\ 0 \\ 0 \\ g\ h_f \partial_x h_s + g\ h_f \partial_x b \\ 0 \\ 0 \end{pmatrix} + \begin{pmatrix} 0 \\ 0 \\ g\ h_s \rho\ \partial_y h_f + g\ h_s \partial_y b \\ 0 \\ 0 \\ g\ h_f \partial_y h_s + g\ h_f \partial_y b \\ 0 \end{pmatrix} = 0$$

All the linearly degenerate eigenvectors in both directions can be evaluated analytically. The non-linear eigenvectors and their eigenvalues have to be found numerically.

In Section 4.2 of Dumbser *et al.* [21] a 2D smooth vortex test problem was designed for the debris flow model. For this accuracy study we used $\rho = 0.9$ and $g = 10$. Appendix B of this paper fully documents this vortex problem for the multiphase debris flow equations for the sake of completeness. Table IV shows the 2D test from that paper. As before, we see that the scheme reaches its design accuracy.

**Table IV shows the accuracy of the 2D Multiphase Debris Flow Model Vortex problem using the LLF Riemann solver; the variable $h_s$ is shown. The first half of the table shows the results**

**from WENO-AO with non-linear limiting. The second half shows the same results when Multiresolution WENO reconstruction is used. Here the lowest order polynomial in the Multiresolution WENO was modified to have piecewise linear variation with the $MC_\beta$ limiter.**

| WENO-AO-3 | $L_1$ Error | $L_1$ Accuracy | $L_{inf}$ Error | $L_{inf}$ Accuracy |
|---|---|---|---|---|
| $128^2$ | 6.89476E-04 | | 2.47006E-02 | |
| $256^2$ | 1.38601E-04 | 2.31 | 6.92758E-03 | 1.83 |
| $512^2$ | 2.65048E-05 | 2.39 | 1.58420E-03 | 2.13 |
| $1024^2$ | 4.25903E-06 | 2.64 | 4.89453E-04 | 1.69 |
| WENO-AO-(5,3) | | | | |
| $64^2$ | 1.56725E-03 | | 1.09193E-01 | |
| $128^2$ | 1.45710E-04 | 3.43 | 2.51346E-02 | 2.12 |
| $256^2$ | 5.19690E-06 | 4.81 | 9.78580E-04 | 4.68 |
| $512^2$ | 1.65696E-07 | 4.97 | 3.19263E-05 | 4.94 |
| WENO-AO-(7,3) | | | | |
| $64^2$ | 1.08266E-03 | | 8.58208E-02 | |
| $128^2$ | 1.48342E-05 | 6.19 | 3.45509E-03 | 4.63 |
| $256^2$ | 1.33724E-07 | 6.79 | 3.40910E-05 | 6.66 |
| $512^2$ | 1.08264E-09 | 6.95 | 2.84816E-07 | 6.90 |
| WENO-AO-(9,3) | | | | |
| $32^2$ | 6.37146E-03 | | 4.44387E-01 | |
| $64^2$ | 7.47838E-04 | 3.09 | 4.17393E-02 | 3.41 |
| $128^2$ | 2.00125E-06 | 8.55 | 5.03900E-04 | 6.37 |
| $256^2$ | 4.88885E-09 | 8.68 | 1.43552E-06 | 8.46 |
| | | | | |
| Order 3 Multires WENO | | | | |
| $128^2$ | 4.47909E-04 | | 1.51502E-02 | |
| $256^2$ | 7.21345E-05 | 2.63 | 3.45687E-03 | 2.13 |
| $512^2$ | 1.16748E-05 | 2.63 | 7.08754E-04 | 2.29 |
| $1024^2$ | 1.78774E-06 | 2.71 | 2.07025E-04 | 1.78 |
| Order 5 Multires WENO | | | | |
| $64^2$ | 1.46645E-03 | | 7.24576E-02 | |
| $128^2$ | 1.81027E-04 | 3.02 | 2.08778E-02 | 1.80 |
| $256^2$ | 6.56540E-06 | 4.79 | 1.04102E-03 | 4.33 |
| $512^2$ | 2.27916E-07 | 4.85 | 3.98982E-05 | 4.71 |
| Order 7 Multires WENO | | | | |
| $64^2$ | 7.18279E-04 | | 4.62518E-02 | |
| $128^2$ | 1.52092E-05 | 5.56 | 2.80551E-03 | 4.04 |
| $256^2$ | 1.39429E-07 | 6.77 | 3.32243E-05 | 6.40 |
| $512^2$ | 1.08301E-09 | 7.01 | 2.82047E-07 | 6.88 |
| Order 9 Multires WENO | | | | |
| $32^2$ | 6.10595E-03 | | 4.23764E-01 | |
| $64^2$ | 6.97989E-04 | 3.13 | 1.03191E-01 | 2.04 |
| $128^2$ | 2.00655E-06 | 8.44 | 4.86319E-04 | 7.73 |
| $256^2$ | 4.88896E-09 | 8.68 | 1.43483E-06 | 8.40 |

**VII.d) Accuracy Study for Two Dimensional Two-layer Shallow Water Equations**

The two layer shallow water equations that we use here obtain from the paper of Castro *et al*. [15]. The PDE system in 2D is given by defining $h_1$, $u_1$, $v_1$ as the height of the upper fluid, its x-velocity and its y-velocity respectively and by defining $h_2$, $u_2$, $v_2$ as the height of the lower fluid, its x-velocity and its y-velocity respectively. The bottom topography is denoted by "*b*" and "*g*" is the gravity. The ratio $\rho \equiv \rho_1/\rho_2$ denotes the ratio of the fluid densities. The total surface height, therefore, becomes $\eta = \eta_1 = b + h_1 + h_2$ and the surface elevation of the interior layer is denoted by $\eta_2 = b + h_2$. The conservation law is given by

$$\frac{\partial}{\partial t}\begin{pmatrix} h_1 \\ h_1u_1 \\ h_1v_1 \\ h_2 \\ h_2u_2 \\ h_2v_2 \\ b \end{pmatrix} + \frac{\partial}{\partial x}\begin{pmatrix} h_1u_1 \\ h_1u_1^2 + gh_1^2/2 \\ h_1u_1v_1 \\ h_2u_2 \\ h_2u_2^2 + gh_2^2/2 \\ h_2u_2v_2 \\ 0 \end{pmatrix} + \frac{\partial}{\partial y}\begin{pmatrix} h_1v_1 \\ h_1u_1v_1 \\ h_1v_1^2 + gh_1^2/2 \\ h_2v_2 \\ h_2u_2v_2 \\ h_2v_2^2 + gh_2^2/2 \\ 0 \end{pmatrix} + \begin{pmatrix} 0 \\ gh_1\partial_x h_2 + g\,h_1\partial_x b \\ 0 \\ 0 \\ \rho g h_2 \partial_x h_1 + gh_2\partial_x b \\ 0 \\ 0 \end{pmatrix} + \begin{pmatrix} 0 \\ 0 \\ gh_1\partial_y h_2 + gh_1\partial_y b \\ 0 \\ 0 \\ \rho g h_2\partial_y h_1 + gh_2\partial_y b \\ 0 \end{pmatrix} = 0$$

All the linearly degenerate eigenvectors in both directions can be evaluated analytically. The non-linear eigenvectors and their eigenvalues have to be found numerically.

In Section 4.1 of Dumbser *et al*. [21] a 2D smooth vortex test problem was designed for the two-layer shallow water flow model. For this accuracy study we used $\rho = 0.9$ and $g = 10$. Appendix C of this paper fully documents this vortex problem for the two-layer shallow water equations for the sake of completeness. Table V shows the 2D test from that paper. As before, we see that the scheme reaches its design accuracy.

**Table V shows the accuracy of the 2D Two-layer Shallow Water Equations Vortex problem using the LLF Riemann solver; the variable $h_1$ is shown. The first half of the table shows the results from WENO-AO with non-linear limiting. The second half shows the same results when Multiresolution WENO reconstruction is used. Here the lowest order polynomial in the Multiresolution WENO was modified to have piecewise linear variation with the $MC_\beta$ limiter.**

| WENO-AO-3 | $L_1$ Error | $L_1$ Accuracy | $L_{inf}$ Error | $L_{inf}$ Accuracy |
|---|---|---|---|---|
| $128^2$ | 3.51577E-03 | | 8.48035E-02 | |
| $256^2$ | 7.52103E-04 | 2.22 | 2.10029E-02 | 2.01 |
| $512^2$ | 1.17975E-04 | 2.67 | 4.39758E-03 | 2.26 |
| $1024^2$ | 1.67899E-05 | 2.81 | 1.14952E-03 | 1.94 |
| WENO-AO-(5,3) | | | | |
| $64^2$ | 4.37768E-04 | | 2.73332E-02 | |
| $128^2$ | 4.87690E-05 | 3.17 | 1.21278E-02 | 1.17 |
| $256^2$ | 1.86615E-06 | 4.71 | 5.55867E-04 | 4.45 |
| $512^2$ | 5.99287E-08 | 4.96 | 1.82648E-05 | 4.93 |
| WENO-AO-(7,3) | | | | |
| $64^2$ | 1.98274E-04 | | 3.63034E-02 | |
| $128^2$ | 5.03312E-06 | 5.30 | 1.96856E-03 | 4.20 |
| $256^2$ | 4.67313E-08 | 6.75 | 1.97459E-05 | 6.64 |
| $512^2$ | 3.81742E-10 | 6.94 | 1.62755E-07 | 6.92 |
| WENO-AO-(9,3) | | | | |
| $32^2$ | 6.12951E-03 | | 4.27426E-01 | |
| $64^2$ | 1.32529E-04 | 5.53 | 1.70932E-02 | 4.64 |
| $128^2$ | 6.42973E-07 | 7.69 | 2.92168E-04 | 5.87 |
| $256^2$ | 1.59016E-09 | 8.66 | 7.54719E-07 | 8.60 |
| | | | | |
| Order 3 Multires WENO | | | | |
| $64^2$ | 2.15283E-03 | | 4.56624E-02 | |
| $128^2$ | 3.74725E-04 | 2.52 | 8.87986E-03 | 2.36 |
| $256^2$ | 5.25003E-05 | 2.84 | 2.05985E-03 | 2.11 |
| $512^2$ | 7.21921E-06 | 2.86 | 5.34445E-04 | 1.95 |
| Order 5 Multires WENO | | | | |
| $64^2$ | 7.91399E-04 | | 8.68716E-02 | |
| $128^2$ | 7.84474E-05 | 3.33 | 1.09385E-02 | 2.99 |
| $256^2$ | 2.74993E-06 | 4.83 | 5.96627E-04 | 4.20 |
| $512^2$ | 6.78261E-08 | 5.34 | 2.20924E-05 | 4.76 |
| Order 7 Multires WENO | | | | |
| $64^2$ | 2.03243E-04 | | 1.84106E-02 | |
| $128^2$ | 5.48277E-06 | 5.21 | 1.50036E-03 | 3.62 |
| $256^2$ | 5.24435E-08 | 6.71 | 2.24150E-05 | 6.06 |
| $512^2$ | 3.83162E-10 | 7.10 | 1.61585E-07 | 7.12 |
| Order 9 Multires WENO | | | | |
| $32^2$ | 5.07733E-03 | | 3.72370E-01 | |
| $64^2$ | 4.32464E-04 | 3.55 | 8.23819E-02 | 2.18 |
| $128^2$ | 6.85031E-07 | 9.30 | 2.81237E-04 | 8.19 |
| $256^2$ | 1.59031E-09 | 8.75 | 7.55231E-07 | 8.54 |

## VIII) One-dimensional Test Problems for Hyperbolic PDE Systems with Non-Conservative Products; Including Well-Balancing

In this Section, we focus on one-dimensional Riemann problems as well as some test problems that show that the algorithm is capable of handling "lake-at-rest" types of problems. We focus on hyperbolic PDE systems with non-conservative products that are of scientific or engineering interest. In Sub-Section VIII.a we focus on the multiphase debris flow model of Pitman and Le. In Sub-section VIII.b we focus on the Baer-Nunziato model of compressible multi-phase flows. In Sub-section VIII.c we present results from the two-layer shallow water equations. Just to show that the method also handles conservation laws well, in Sub-section VIII.d we present a very stringent result involving Euler flow. For all the simulations presented in this Section we used a CFL of 0.8 with a third order SSP-RK scheme.

### VIII.a) One-dimensional Test Problems for Multiphase Debris Flow Model of Pitman and Le

We present one-dimensional test problems for the multiphase debris flow model of Pitman and Le [49] in the formulation of Pelanti *et al*. [46]. The equations for the multiphase debris flow model are described in Sub-section VII.c. The test problems we present are drawn from Pelanti *et al*. [46], Dumbser *et al*. [21] and Rhebergen *et al*. [51] and they have also been concatenated in Dumbser and Balsara [26].

Table VI describes the three Riemann problems that we present here; they are drawn from Dumbser and Balsara [26]. In all three Riemann problems we put $\rho = 0.5$ and $g = 9.8$. In all cases we use a 1D mesh with 200 zones. The computational domain and the final stopping time is also shown in Table VI. Fig. 2, which corresponds to the first Riemann problem in Table VI, corresponds to the superposition of three stationary jump discontinuities in the linearly degenerate fields. The free surface is flat with a bottom jump and a jump in the depths of the solid and the fluid. This is coupled with two shear waves in the solid and fluid respectively. The solid volume fraction is also flat in this problem. The problem was run with the anti-diffusive fluxes turned on for all the linearly degenerate waves. Fig. 2 shows that all the jumps in the linearly degenerate stationary wave families are exactly preserved on the mesh at all orders. Please note that all the results from third, fifth and seventh order FD-WENO that is based on the HLL Riemann solver look identical to Fig. 2. Likewise, all the results from third, fifth and seventh order FD-WENO that is based on the LLF Riemann solver also look identical to Fig. 2. In general, an LLF Riemann solver is considered much more dissipative than an HLL Riemann solver. However, owing to the fact that the anti-diffusive fluxes can also self-adjust, both algorithms give the same excellent quality of solution. For all the simulations that we show in Fig. 2, we also measured that the velocity fluctuations in the x-direction remain less than $1.4\times10^{-14}$. Fig. 3 shows the next two 1D Riemann problems from Table VI. Figs. 3a, 3b and 3c show the results of these two 1D Riemann problems when they were run with third, fifth and seventh order FD-WENO schemes that are based on the HLL Riemann solver. A third order FD-WENO scheme was also run on a mesh with 2000 zones in order to generate a reference solution, and the reference solution is shown with solid lines in Fig. 3. We see that the solution quality is very good at third order and indeed improves at fifth and seventh orders, as expected.

**Table VI gives the left and right initial states, computational domain $\left[x_L, x_R\right]$, location of the discontinuity $(x_c)$ and final times $(t_{end})$ of the three Riemann problems for the Multiphase Debris Flow Model.**

| Test Case: | | $h_s$ | $u_s$ | $v_s$ | $h_f$ | $u_f$ | $v_f$ | $b$ | $x_L$ | $x_R$ | $x_c$ | $t_{end}$ |
|---|---|---|---|---|---|---|---|---|---|---|---|---|
| **RP1** | Left: | 1.5 | 0.0 | 0.2 | 0.5 | 0.0 | -0.5 | 0.0 | -5.0 | 5.0 | 0.0 | 1.0 |
| | Right: | 1.125 | 0.0 | -0.2 | 0.375 | 0.0 | 0.5 | 0.5 | | | | |
| **RP2** | Left: | 2.1 | 0.0 | 0.0 | 0.9 | 0.0 | 0.0 | 0.0 | -5.0 | 5.0 | 0.0 | 0.5 |
| | Right: | 0.8 | 0.0 | 0.0 | 1.2 | 0.0 | 0.0 | 0.0 | | | | |
| **RP3** | Left: | 2.1 | -1.4 | 0.0 | 0.9 | 0.3 | 0.0 | 0.0 | -5.0 | 5.0 | 0.0 | 0.5 |
| | Right: | 0.8 | -0.9 | 0.0 | 1.2 | 0.1 | 0.0 | 0.0 | | | | |

It is also customary and desirable to show that the algorithm produces well-balanced results, and we do the same here. In all our simulations to demonstrate the well-balancing, we used the code without any special modifications for the well-balancing. The well-balanced solution for this problem is always given by

$$\phi = \frac{h_s}{h_s + h_f} = const. \quad ; \quad \eta = h_s + h_f + b = const. \quad ; \quad u_s = v_s = u_f = v_f = 0$$

which is why we plot out "$\eta$" as a way of demonstrating the approach to equilibrium. The test problem for demonstrating well-balancing for the multiphase debris flow model comes from Pelanti *et al*. [46] which builds on the work by LeVeque [41]. It consists of setting the bottom topography as

$$b(x) = \begin{cases} 0.25\left(\cos\left(10\pi\left(x-1/2\right)\right)+1\right) & \text{if} \quad |x-0.5| \le 0.1 \\ 0 & \text{otherwise} \end{cases}$$

All velocities are initially set to zero and a small perturbation is added to the free surface elevation $\eta$ and the solid volume fraction $\phi$ as follows:-

$$\eta(x,0) = \eta_0 + \varepsilon \quad ; \quad \phi(x,0) = \phi_0 - \varepsilon \qquad \text{for} \quad -0.6 \le \text{x} \le -0.5$$

Here we set $\eta_0 = 1$, $\phi_0 = 0.6$, $g = 1$ and $\rho = 0.5$. The problem was run on a 400 zone computational domain spanning $\left[-1.2, 1.2\right]$ with a stopping time of 1.25. We used $\varepsilon = 10^{-3}$. The results are shown in Fig. 4. We see that the surface elevation at early and late times remains bounded, showing that a well-balanced solution has been achieved. Our solutions are competitive with the reference results shown in Dumbser and Balsara [26].

### VIII.b) One-dimensional Test Problems for Baer-Nunziato Model for Compressible Multi-Phase Flows

The equations for the Baer-Nunziato model were already described in Section VII.b. For a 2D problem, each phase can sustain an entropy wave (contact discontinuity) and a shear wave along with a further wave that corresponds to a jump in the volume fraction. The Abgrall condition, see Saurel and Abgrall [52], is based on the idea that a mixture of the two phases moving in 1D with uniform velocity and pressure should be able to continue moving in this fashion without generating any wiggles in the velocity or pressure. We present such a test problem here on a 200 zone mesh that spans the domain $[-0.5, 0.5]$. In the solid phase we set a stiffened EOS with $\gamma_1 = 3$ and $\pi_1 = 100$, whereas in the gas phase we set $\gamma_2 = 1.4$ and $\pi_2 = 0$. The pressures and longitudinal velocities in both phases are set to unity. The problem is initialized with $\rho_{1,L} = 800$, $\rho_{2,L} = 2$, $\phi_{1,L} = 0.99$ on the left of $x = 0$ and $\rho_{1,R} = 1000$, $\rho_{2,R} = 1$, $\phi_{1,L} = 0.01$ to its right. The problem was run to a final time of $t = 0.25$. The results of pressure, velocity and solid volume fraction are shown in Fig. 5a at third, fifth and seventh orders when the HLL Riemann solver is used. The same information in shown in Fig. 5b at third, fifth and seventh orders when the LLF Riemann solver is used. We see that the pressure and velocity profiles are absolutely flat. The pressure in Fig. 5 varies from unity by less than $9\times10^{-13}$; the velocities in the solid and gas vary from unity by less than $2\times10^{-13}$. As expected, we see a clear improvement in the crispness of the volume fraction profile as we go from third to higher orders. While an LLF Riemann solver is considered much more dissipative than an HLL Riemann solver, we see that the anti-diffusive term has compensated for that so that the profiles in Figs. 5a and 5b are entirely comparable.

**Table VII gives the left and right initial states, computational domain $[x_L, x_R] = [-0.5, 0.5]$ and final times $(t_{end})$ of the six Riemann problems for Baer-Nunziato model.**

| Test Case: | $\rho_1$ | $u_1$ | $p_1$ | $\rho_2$ | $u_2$ | $p_2$ | $\phi_1$ | $t_{end}$ |
|---|---|---|---|---|---|---|---|---|
| **RP1** | $\gamma_1 = 1.4, \quad \pi_1 = 0, \quad \gamma_2 = 1.4, \quad \pi_2 = 0$ | | | | | | | |
| Left: | 1.0 | 0.0 | 1.0 | 0.5 | 0.0 | 1.0 | 0.4 | 0.1 |
| Right: | 2.0 | 0.0 | 2.0 | 1.5 | 0.0 | 2.0 | 0.8 | |
| **RP2** | $\gamma_1 = 3.0, \quad \pi_1 = 100, \quad \gamma_2 = 1.4, \quad \pi_2 = 0$ | | | | | | | |
| Left: | 800.0 | 0.0 | 500.0 | 1.5 | 0.0 | 2.0 | 0.4 | 0.1 |
| Right: | 1000.0 | 0.0 | 600.0 | 1.0 | 0.0 | 1.0 | 0.3 | |
| **RP3** | $\gamma_1 = 1.4, \quad \pi_1 = 0, \quad \gamma_2 = 1.4, \quad \pi_2 = 0$ | | | | | | | |
| Left: | 1.0 | 0.9 | 2.5 | 1.0 | 0.0 | 1.0 | 0.9 | 0.1 |
| Right: | 1.0 | 0.0 | 1.0 | 1.2 | 1.0 | 2.0 | 0.2 | |
| **RP4** | $\gamma_1 = 3.0, \quad \pi_1 = 3400, \quad \gamma_2 = 1.35, \quad \pi_2 = 0$ | | | | | | | |
| Left: | 1900.0 | 0.0 | 10.0 | 2.0 | 0.0 | 3.0 | 0.2 | 0.15 |
| Right: | 1950.0 | 0.0 | 1000.0 | 1.0 | 0.0 | 1.0 | 0.9 | |
| **RP5** | $\gamma_1 = 1.4,\quad \pi_1 = 0,\quad \gamma_2 = 1.4,\quad \pi_2 = 0$ | | | | | | | |
| Left:<br>Right: | 1.0<br>1.0 | 0.0<br>0.0 | 1.0<br>1.0 | 0.2<br>1.0 | 0.0<br>0.0 | 0.3<br>1.0 | 0.8<br>0.3 | 0.20 |
| **RP6** | $\gamma_1 = 1.4,\quad \pi_1 = 0,\quad \gamma_2 = 1.4,\quad \pi_2 = 0$ | | | | | | | |
| Left:<br>Right: | 0.2068<br>2.2263 | 1.4166<br>0.9366 | 0.0416<br>6.0 | 0.5806<br>0.4890 | 1.5833<br>-0.70138 | 1.375<br>0.986 | 0.1<br>0.2 | 0.1 |

In Dumbser *et al*. [22] a fine set of six 1D Riemann problems were compiled for the Baer-Nunziato system. Table VII catalogues the parameters for these six Riemann problems. The results for the FD-WENO scheme are shown in Fig. 6 when the HLL Riemann solver was used. Figs. 6a and 6b show the solid and gas densities for RP1 to RP3 and then again for RP4 to RP6 when a third order scheme was used. Figs. 6c and 6d show similar information when a fifth order scheme was used. Figs. 6e and 6f show similar information when a seventh order scheme was used. A third order FD-WENO scheme was also run on a mesh with 2000 zones in order to generate a reference solution, and the reference solution is shown with solid lines in Fig. 6. We see that the higher order schemes represent the interfaces more crisply. But we also see that inevitable wiggles close to shocks are also exaggerated by the higher order schemes. In all of Fig. 6 we switched off any contribution from the anti-diffusive flux. Since a good flattener function, i.e. a good shock detection algorithm, has not been designed in the literature for the Baer-Nunziato we had to set $\psi = 1$ in eqn. (3.8). The availability of advanced algorithms, like the ones given here, for the Baer-Nunziato system, makes a good case for finding a good shock detection algorithm for this system. But such work is out of scope here.

It is also worth noting that the solid and gas mass densities in Baer-Nunziato flow are indeed in flux conservation form; even though the entire system is non-conservative. This has the consequence that when the velocities at the left and right boundaries are zero, the mass fluxes at the two boundaries should also be zero. As a result, the spatially integrated solid and gas mass should be constant as a function of time. Riemann problems 2, 4 and 5 in Table VII have zero velocities at the left and right boundaries. As a result, we can verify that the integrated solid and gas masses are constant as a function of simulation time for these three test problems. This was done by evaluating the total solid and gas masses on the mesh at the initial time (subscript "0" in the legends). We then kept track of the integrated solid and gas masses as a function of time. Fig.

6g shows the ratio of the mass to the initial mass as a function of time. We see that the ratio is extremely close to unity. With the help of Fig. 6g we verify that mass conservation is respected to at least one part in $10^{12}$ by the methods developed here.

**VIII.c) One-dimensional Test Problems for Two-layer Shallow Water Equations**

The two-layer shallow water equations were described in Sub-section VII.d. For all the problems shown in this sub-section we used $\rho = 0.8$ and $g = 9.8$. The three Riemann problems that we would like to show have been concatenated in Dumbser and Balsara [26] but derive from Dumbser *et al*. [21] and Castro *et al*. [14]. Table VIII gives the relevant information for these Riemann problems. All the two-layer shallow water Riemann problems were run on a 200 zone mesh. The first Riemann problem from Table VIII illustrates the preservation of stationary jump discontinuities in the linearly degenerate intermediate fields, and the results are shown in Fig. 7 for third order FD-WENO schemes that use the HLL Riemann solver; the results from higher order FD-WENO schemes are identical. The use of the LLF Riemann solver with the same schemes produces identical results. This shows the ability of our high order schemes to preserve stationary discontinuities in linearly degenerate intermediate fields. Fig. 8 shows the results from the remaining two Riemann test problems in Table VIII. A third order FD-WENO scheme was also run on a mesh with 2000 zones in order to generate a reference solution, and the reference solution is shown with solid lines in Fig. 8. We see that the results in Fig. 8 indeed closely match the reference solution.

**Table VIII gives the left and right initial states, computational domain $[x_L, x_R]$, location of the discontinuity $(x_c)$ and final times $(t_{end})$ of the three Riemann problems for the Two-Layer Shallow Water Model.**

| Test Case: | | $h_1$ | $u_1$ | $v_1$ | $h_2$ | $u_2$ | $v_2$ | $b$ | $x_L$ | $x_R$ | $x_c$ | $t_{end}$ |
|---|---|---|---|---|---|---|---|---|---|---|---|---|
| **RP1** | Left: | 0.5 | 0.0 | 0.5 | 0.8 | 0.0 | -0.2 | 0.2 | -5.0 | 5.0 | 0.0 | 1.0 |
| | Right: | 0.5 | 0.0 | -0.5 | 0.2 | 0.0 | 0.2 | 0.8 | | | | |
| **RP2** | Left: | 0.4 | 0.0 | 0.0 | 0.6 | 0.0 | 0.0 | 0.0 | -5.0 | 5.0 | 0.0 | 1.25 |
| | Right: | 0.6 | 0.0 | 0.0 | 0.4 | 0.0 | 0.0 | 0.0 | | | | |
| **RP3** | Left: | 1.0 | 0.0 | 0.0 | 1.0 | 0.0 | 0.0 | 0.0 | -5.0 | 5.0 | 0.0 | 1.0 |
| | Right: | 0.5 | 0.0 | 0.0 | 0.5 | 0.0 | 0.0 | 0.5 | | | | |

We would also like to illustrate that our algorithm produces well-balanced results. In all our simulations to demonstrate the well-balancing, we used the code without any special modifications for the well-balancing. The well-balanced solution for this problem is always given by

$$h_1 = const. \quad ; \quad \eta = h_1 + h_2 + b = const. \quad ; \quad u_1 = v_1 = u_2 = v_2 = 0$$

We also denote the elevation of the interior layer by $\eta_2 = b + h_2$. The problem was run on a 400 zone mesh that spans $[0,2]$ to a final time of 0.2. The present test problem is a small variation on LeVeque [41] where the initial velocities are zero and the bottom topography is given by

$$b(x) = \begin{cases} 0.25\left(\cos\left(10\pi(x-1.5)\right)+1\right) & \text{if} \quad 1.4 \le x \le 1.6 \\ 0 & \text{otherwise} \end{cases}$$

For this problem we initially set $\eta_2(x,0) = 0.65$ for all values of "*x*". We also initialize

$$\eta(x,0) = \begin{cases} 1+\varepsilon & \text{if} \quad 1.1 \le \text{x} \le 1.2 \\ 1 & \text{otherwise} \end{cases}$$

Two cases were run. The first case had a large perturbation with $\varepsilon = 0.2$. The second case had a small perturbation with $\varepsilon = 10^{-3}$. The left panel in Fig. 9a shows the $\eta(x,t)$, $\eta_2(x,t)$ and the bottom topography at the final time of $t = 0.2$ for the large perturbation; the right panel shows $\eta(x,t)$ and $\eta_2(x,t)$ at the final time of $t = 0.2$ for the small perturbation. Fig. 9a shows the results when a third order FD-WENO scheme was used with an HLL Riemann solver. Figs. 9b and 9c show similar results when fifth and seventh order FD-WENO schemes were used with an HLL Riemann solver. We find a very good agreement with the reference results presented in Dumbser and Balsara [26], showing that our method is well-balanced. As in the debris flow model, we did not use any special variables for the well-balancing; the native scheme, without any modification, was able to achieve well-balanced results.

**VIII.d) Euler Flow**

It is also worthwhile to show that for an extremely stringent problem involving conservation laws, the method performs well. We, therefore, focus on the blast wave interaction problem from Woodward and Colella [62]. The problem uses the Euler system, which is in conservation form. The problem consists of two blast waves that interact with one another. The 1D domain has 1000 zones and spans $x \in [-0.5, 0.5]$. The problem was stopped at a time of $t$ = 0.038. Fig. 10 shows the resulting density profile at the final time. The panels show the density variable for third, fifth and seventh order schemes. A reference calculation was run using a third order scheme on a 10,000 zone mesh and is superposed on all the figures as a solid line. The precise coincidence of all the density profiles with the reference solution gives one confirmation that our method produces convergent shock locations for conservation laws.

## IX) Multidimensional Test Problems for Hyperbolic PDE Systems with Non-Conservative Products

In this Section, we focus on several two-dimensional test problems. Sub-sections IX.a and IX.b show test problems involving Baer-Nunziato compressible multi-phase flow. Sub-sections IX.c and IX.d show test problems involving two-layer shallow water flow. Sub-sections IX.e and IX.f show test problems involving multiphase debris flow. Lastly, Sub-sections IX.g and IX.h show two very popular test problems involving Euler flow. For all the simulations presented in this Section we used a third order SSP-RK scheme.

### IX.a) Shock-Bubble Interaction for Baer-Nunziato Compressible Multi-Phase Flow

This shock-bubble interaction problem for the Baer-Nunziato compressible multi-phase flow was first presented in Dumbser *et al*. [24]. It enables us to demonstrate that the higher order FD-WENO schemes presented here work well for stringent problems with strong shocks for hyperbolic systems that have non-conservative products.

The problem is set up on a domain that spans $[-0.5, 3.0] \times [-0.75, 0.75]$ and is run to a final time of 0.025 with a CFL of 0.4. The problem consists of a planar, right-going shock propagating into an ambient medium which contains a bubble with radius 0.25. The shock that flows in from the left boundary has its parameters denoted by an "L", the ambient medium into which the shock propagates has its parameters denoted by an "R" and the bubble's parameters are denoted with a "B" in Table IX. The left and right boundaries are Dirichlet and the top and bottom boundaries are periodic. The shock is initialized at $x = 0$ and the bubble is centered at $(0.5, 0)$. The parameters for the equation of state are given by $\gamma_1 = 3.0,\ \pi_1 = 100,\ \gamma_2 = 1.4,\ \pi_2 = 0.$

**Table IX gives the left (L), right (R) and bubble (B) state for the shock-bubble interaction.**

| | $\rho_1$ | $u_1$ | $p_1$ | $\rho_2$ | $u_2$ | $p_2$ | $\phi_1$ |
|---|---|---|---|---|---|---|---|
| Left (L) | 1999.939402 | 49.998485 | 4999849.5 | 1.0 | 0.0 | 1.0 | 0.75 |
| Right (R) | 1000.0 | 0.0 | 1.0 | 1.0 | 0.0 | 1.0 | 0.75 |
| Bubble (B) | 1000.0 | 0.0 | 1.0 | 1.0 | 0.0 | 1.0 | 0.25 |

Figs. 11a, 11b and 11c show the resulting solid density for third, fifth and seventh order FD-WENO schemes on a mesh with 700×300 zones. We see that the resulting density structures are not resolved very sharply by the third order scheme. However, the fifth and seventh order schemes resolve all the flow structures very well, thereby highlighting the value of higher order schemes.

**IX.b) Shock-Vortex Interaction for Baer-Nunziato Compressible Multi-Phase Flow**

In Balsara and Shu [3] a very useful test problem was constructed for a vortex interacting with a shock involving Euler flow. Since we have well-packaged vortex flow examples in the Appendices, it is possible to now supply analogous test problems for the Baer-Nunziato compressible multi-phase flow.

The domain for this problem spans $[-0.5,1.5]\times[-0.5,1.5]$ and the problem was run to a final time of 0.84 with a CFL of 0.4. A standing shock was initialized that goes from $(x, y) = (-0.5,1)$ to $(x, y) = (1,-0.5)$. The shock is initialized along the line $x + y = 0.5$ and the same initialization is retained for the duration of the simulation on the ghost zones that bound the computational domain. A vortex that is centered at $(x, y) = (0,0)$ is initialized on the mesh. The unperturbed pre-shock flow to the left and bottom of $x + y = 0.5$ is denoted by "LB" and the unperturbed post-shock flow to the right and top of $x + y = 0.5$ is denoted by "RT". The flow variables for these pre-shock and post-shock flows are given in Table X.

**Table X gives the Left-Bottom (LB) and Right-Top (RT) states of the unperturbed shock flow for the Baer-Nunziato Shock-Vortex interaction problem.**

| | $\rho_1$ | $u_1$ | $v_1$ | $p_1$ | $\rho_2$ | $u_2$ | $v_2$ | $p_2$ | $\phi_1$ |
|---|---|---|---|---|---|---|---|---|---|
| LB | 1.0 | 1.1652078 | 1.1652078 | 1.0 | 2.0 | 0.0 | 0.0 | 1.5 | 1.0/3.0 |
| RT | 1.672601 | 0.696644 | 0.696644 | 2.090198 | 2.0 | 0.0 | 0.0 | 1.5 | 1.0/3.0 |

The vortex itself uses the same style of parametrization as in Appendix A, however, the vortex parameters are given by

$$r_s = 0.075, \quad p_{10} = 1, \quad p_{20} = \frac{3}{2}, \quad s_1 = \frac{3}{2}, \quad s_2 = \frac{7}{5}, \quad \gamma_1 = 1.4, \quad \pi_1 = 0, \quad \gamma_2 = 1.4, \quad \pi_2 = 0.$$

The vortex takes on the mean velocity of the pre-shocked flow, with the result that the vortex propagates into the stationary shock. By an intermediate time of 0.23, the vortex will have propagated halfway through the shock and this time is also an interesting time for visualizing the vortex.

Figs. 12a, 12b and 12c show solid volume fraction at times t=0.0, 0.23, 0.84. The intermediate time of 0.23 corresponds to a time when the vortex has propagated halfway through the shock. Note that the shock does not reveal itself in this variable. Figs. 12d, 12e, 12f show solid x-velocity at times t=0.0, 0.23, 0.84. The shock is clearly visible in this variable. We see that as the vortex propagates through the shock, the structure of the vortex is rearranged by the shock. However, the vortex is not fully destroyed by the shock. As a result, the vortex sheds some of its angular momentum after it has passed through the shock in an attempt to rearrange its structure. This is revealed by the spiral arms that are shed by the vortex in Figs. 12c and 12f. However, the figure shows us that vortices are very robust flow structures because they carry the angular momentum of the fluid. (Recall that the Baer-Nunziato equations are rotationally invariant; therefore, by Noether's theorem they must conserve angular momentum.) Because angular momentum is conserved by Baer-Nunziato flow, the vortices will not be fully destroyed as they pass through shocks of modest strength. This was known for Euler flow, but Fig. 12 reconfirms this insight for Baer-Nunziato flow.

### IX.c) Shock-"Bubble" Interaction for Two-layer Shallow Water Flow

The reason for putting "bubble" within apostrophes in the title of this Sub-section is that this is not truly a bubble but is made to look like one by controlling the bathymetry in the two-layer shallow water equations. Here we consider a domain that spans $[-0.5, 3.0] \times [-0.75, 0.75]$ and the problem was run to a final time of 0.30 with a CFL of 0.4. For this problem we take the gravitational acceleration to be $g = 9.8$ and the density ratio in the two fluids to be $\rho \equiv \rho_1 / \rho_2 = 0.8$ ; i.e., the lighter fluid with a density of $\rho_1$ is on top of the bottom fluid with a density of $\rho_2$. The actual densities of the two fluids do not need to be specified. The problem consists of a right-going shock that is initialized at $x = 0$ and propagates into the unshocked medium that exists for $x > 0$. The post-shock variables for the region with $x \le 0$ are given by $\left(h_{1L}, u_{1L}, v_{1L}, h_{1R}, u_{1R}, v_{1R}, b\right) = \left(2.44516, 3.331868, 0, 3.785719, 2.927431, 0, 0\right)$. The two-layer fluid that this shock propagates into is initially at rest and in hydrodynamical equilibrium. For most of the domain, we have $b = 0$ . However, centered at $\left(x, y\right) = \left(0.5, 0.0\right)$ we have a circular depression in the bottom of the domain with a mean radius of 0.25. Let "*r*" denote the distance from the center of this depression, so that we specify the bottom bathymetry with a tapered profile that is given by

$$b\left(r\right) = -0.25 + 0.25 \tanh\left(\left(r - 0.25\right) / \left(2\Delta x\right)\right).$$

Here $\Delta x$ is the zone size, assuming square zones. To maintain hydrodynamical equilibrium we specify the unshocked medium with $x > 0$ by

$$h_1 = 1.5 \qquad \text{and} \qquad h_2\left(r\right) = 4 - h_1 - b\left(r\right)$$

The left and right boundaries are Dirichlet and the upper and lower boundaries are set periodic. If one wishes, the shock can be initialized with a one or two zone taper on the mesh to eliminate any oscillations associated with the initialization of the shock on the mesh.

Fig. 13 shows the final image of the height of the lower fluid. The shock propagates into the unshocked medium without perturbation till it encounters the depression. When it encounters the "bubble" it sets up a left-going bow shock and a wake around the bubble. The wake is shown in green in Fig. 13. There is also a shock that has propagated three-fourths of the way through the "bubble". A very interesting flow structure develops behind the "bubble" as can be seen in Fig. 13. A conical pair of shocks that meet at a focal point behind the "bubble" can also be seen in the figure. Figs. 13a, 13b and 13c correspond to third, fifth and seventh order simulations on a $700 \times 300$ zone mesh. We see that the fifth and seventh order simulations show flow structures that are very crisp, highlighting the value of higher order schemes.

**IX.d) Shock-Vortex Interaction for Two-layer Shallow Water Flow**

In this Sub-section we analyze the interaction of a vortex flow with a stationary shock for the two-layer shallow water equations. For the Euler flows, analogous test problem has been constructed in Balsara and Shu [3].

The computational domain for this problem spans $\left[-0.5,1.5\right]\times\left[-0.5,1.5\right]$. The test problem was run to a final time of 0.24 with a CFL of 0.4. A standing shock was initialized along the line $x+y=0.5$ that goes from $\left(x,y\right)=\left(-0.5,1\right)$ to $\left(x,y\right)=\left(1,-0.5\right)$. The unperturbed pre-shock flow to the left and bottom of the line $x+y=0.5$ is denoted by "LB" and the unperturbed post-shock flow to the right and top of the line $x+y=0.5$ is denoted by "RT". For the entire simulation, the same initialization is used to define the ghost zones that bound the computational domain. The flow variables for these pre-shock and post-shock flows are given in Table XI.

**Table XI gives the Left-Bottom (LB) and Right-Top (RT) states of the unperturbed shock flow for the Two-layer Shallow Water Shock-Vortex interaction problem.**

| | $h_1$ | $u_1$ | $v_1$ | $h_2$ | $u_2$ | $v_2$ | $b$ |
|---|---|---|---|---|---|---|---|
| LB | 1.746044 | 4.849948 | 4.849948 | 0.330408 | 4.849948 | 4.849948 | 0.0 |
| RT | 3.002690 | 2.820212 | 2.820212 | 0.517654 | 3.091776 | 3.091776 | 0.0 |

A vortex that uses the same style of parametrization as in Appendix C is initialized on the computational mesh. However, the vortex parameters for this test problem are given by

$$r_s=0.1,\ g=10,\ \rho=\frac{9}{10},\ s_1=\frac{1}{2},\ s_2=1,\ v_{10}=\frac{3}{4},\ v_{20}=\frac{1}{10},\ h_{10}=1 \text{ and } h_{20}=1.$$

The vortex propagates diagonally into the stationary shock with the mean velocity of the pre-shocked region. By an intermediate time of 0.06, the vortex will have propagated halfway through the stationary-shock so that the vortex core is strongly interacting with the shock and interesting transition of the vortex is observed in the Figs. 14b and 14e.

Figs. 14a, 14b and 14c show height of the upper fluid at times t=0.0, 0.06, 0.24. The intermediate time of 0.06 corresponds to a time when the vortex has propagated halfway through the shock and a little bend in the in the shock front can be observed at this time. Figs. 14d, 14e,

14f show x-velocity of the upper fluid at times t=0.0, 0.06, 0.24. Similar to the Figs. 14a, 14b and 14c for the height of upper fluid, a stationary shock is clearly visible in the Figs. 14d, 14e, 14f for the x-velocity of the upper fluid. We observe that as the vortex propagates through the standing shock, the vortex is not fully destroyed by the shock flow. Instead, the vortex sheds some of its angular momentum after it has passed through the shock. In doing so, the vortex has rearranged its form. The shock used for this problem is not as strong as the shock used in Sub-section IX.b. Consequently, the arms that are shed by the vortex in Figs. 14c and 14f are not a prominent as the ones in sub-section IX.b.. Analogous to the Baer-Nunziato flow equations, the two-layer shallow water equations preserve angular momentum. As a result, the core of the vortex is not destroyed as it passes through the shock.

**IX.e) Shock-"Bubble" Interaction for Multiphase Debris Flow**

As before, we put "bubble" within apostrophes because we will form the bubble by controlling the bathymetry. Here we consider a domain that spans $\left[-0.5, 3.0\right]\times\left[-1.25, 1.25\right]$ and the problem was run to a final time of 0.30 with a CFL of 0.4. For this problem we take the gravitational acceleration to be $g = 9.8$ and the density ratio in the two fluids to be $\rho \equiv \rho_f / \rho_s = 0.5$; i.e., the fluid with a density of $\rho_f$ is on top of the solid with a density of $\rho_s$. The problem consists of a right-going shock that is initialized at $x = 0$ and propagates into the unshocked medium that exists for $x > 0$. The post-shock variables for the region with $x \le 0$ are given by $\left(h_{sL}, u_{sL}, v_{sL}, h_{fR}, u_{fR}, v_{fR}, b\right) = \left(2.3394, 3.17729, 0, 3.899, 3.17742, 0, 0\right)$. This shock propagates into medium that is initially at rest and in hydrodynamical equilibrium. For most of the domain, we have $b = 0$. However, centered at $\left(x, y\right) = \left(0.5, 0.0\right)$ we have a circular depression in the bottom of the domain with a mean radius of 0.25. Let "*r*" denote the distance from the center of this depression, so that we specify the bottom bathymetry with a tapered profile that is given by

$$b\left(r\right) = -0.5 + 0.5 \tanh\left(\left(r - 0.25\right) / \left(2\Delta x\right)\right).$$

Here $\Delta x$ is the zone size, assuming square zones. To maintain hydrodynamical equilibrium we specify the unshocked medium with $x > 0$ by

$$h_s = 1.5 \qquad \text{and} \qquad h_f\left(r\right) = 4 - h_s - b\left(r\right)$$

The left and right boundaries are Dirichlet and the upper and lower boundaries are set periodic.

Fig. 15 shows the final image of the height of the fluid phase. The shock propagates into the unshocked medium without perturbation till it encounters the depression. When it encounters the "bubble" it sets up a left-going bow shock. However, the rest of the flow structure that is set up here is very different from the flow that developed in Fig. 13. In Fig. 15, a ring-like structure emanates from the location of the depression. Figs. 15a, 15b and 15c correspond to third, fifth and seventh order simulations on a $700 \times 500$ zone mesh. We see that the fifth and seventh order simulations show flow structures that are very crisp, highlighting the value of higher order schemes.

**IX.f) Shock-Vortex Interaction for Multiphase Debris Flow**

Analogous to the Sub-sections IX.b and IX.d, it is possible to design a shock-vortex interaction test problem for the Multiphase Debris Flow equations. We begin by considering a computational domain which spans $[-0.5,1.5]\times[-0.5,1.5]$. As before, a stationary shock was initialized along the line $x+y=0.5$ that goes from $(x,y)=(-0.5,1)$ to $(x,y)=(1,-0.5)$. The unperturbed pre-shock flow to the left and bottom of the line $x+y=0.5$ is denoted by "LB" and the unperturbed post-shock flow to the right and top of the line $x+y=0.5$ is denoted by "RT". The same initialization is retained to define the ghost zones during the entire simulation. The flow variables for these pre-shock and post-shock flows are given in Table XII.

**Table XII gives the Left-Bottom (LB) and Right-Top (RT) states of the unperturbed shock flow for the Debris Flow Shock-Vortex interaction problem.**

| | $h_s$ | $u_s$ | $v_s$ | $h_f$ | $u_f$ | $v_f$ | $b$ |
|---|---|---|---|---|---|---|---|
| LB | 1.746044 | 4.694546 | 4.694546 | 0.330408 | 4.694546 | 4.694546 | 0.0 |
| RT | 2.825618 | 2.900917 | 2.900917 | 0.534696 | 2.900921 | 2.900921 | 0.0 |

Following Appendix B, a vortex centered at $(0,0)$ has been initialized in the pre-shock region. The vortex is considered in the same style as in Appendix B with the particular values of the parameters given by $r_s = 0.1,\ g = 10,\ \rho = \frac{9}{10},\ s_1 = \frac{1}{2},\ s_2 = 1,\ v_{10} = \frac{3}{4},\ v_{20} = \frac{1}{10},\ h_{10} = 1$ and $h_{20} = 1$.

The vortex propagates diagonally into the standing shock with the mean velocity given by the velocity of the pre-shock flow. The test problem was run to a final time of 0.24 with a CFL of 0.4. By an intermediate time of t=0.6, half of the vortex has interacted with the shock. As a result, change in the structure of vortex and the shock can be observed in Figs. 16b and 16e.

Figs. 16a, 16b and 16c show the solid height at times t=0.0, 0.06, 0.24. The intermediate simulation time of 0.06 corresponds to a time when the vortex has propagated halfway through the shock. At this time, a slight reduction in the strength of the solid height can be observed in the vortex. Figs. 16d, 16e, 16f show the solid x-velocity at times t=0.0, 0.06, 0.24. The standing shock is clearly visible in all of Fig. 16. We can again verify that the vortex is very robust flow structure in the simulation shown here. As with the two-layer shallow water example, the shock is not a very strong shock and the rearrangement in the structure of the vortex is minimal.

**IX.g) Forward Facing Step Problem for Euler Flow**

We wish to demonstrate that the FD-WENO scheme presented here also works well for stringent problems involving the Euler system, which is indeed in conservation form. This problem was first presented by Woodward and Colella [62]. It consists of a two-dimensional domain that spans $[0,3]\times[0,1]$. An ideal gas flows in from the left boundary at a speed of Mach 3, a density of 1.4 , a pressure of 1 and a ratio of specific heats of 1.4. A forward-facing step is located with its upper corner at the position $(0.6, 0.2)$. Outflow boundary conditions are applied on the right boundary. The top and bottom walls have reflective boundary conditions, and that includes all parts of the forward facing step. We treat the singularity at the corner of the box in the same manner as Woodward and Colella [62]. The simulation is run to a time of 4.0 with a CFL number of 0.4.

Fig. 17 shows the density variable at the final time from our present seventh order FD-WENO algorithm on a 1440×480 zone mesh. The HLLI formulation of the Riemann solver was used. All the shocks have sharp profiles and the vortex sheet shows little or no spreading along the length of the computational domain. We see a very pronounced roll up of the vortex sheet, which demonstrates the value of a higher order scheme. Our fifth order scheme also shows a similar rollup.

### IX.h) Double Mach Reflection Problem for Euler Flow

We wish to build further confidence that the FD-WENO scheme presented here also works well for stringent problems involving the Euler system, which is indeed in conservation form. This problem was presented in Woodward and Colella [62]; and we use the same parameters as those authors. The problem simulates the similarity solution that develops in multi-dimensions when an angled wedge is put in a supersonic flow. The computational domain of this problem is $[0,4]\times[0,1]$. Initially, a Mach 10 shock is positioned at an angle of 60° to the bottom boundary, meeting that boundary at $x = 1/6$. For values of $x < 1/6$, the post-shock conditions are used at the boundary, which mimics the start of the wedge. For $x > 1/6$, the boundary is reflective, which mimics the windward face of the wedge. The upper boundary is set to exactly track the motion of the oblique shock. The left and right edges have inflow and outflow boundary conditions, respectively. The unshocked material is initialized with a density of 1.4, a pressure of unity and a ratio of specific heats of 1.4. The problem was run to a time of 0.2 using a CFL number of 0.4. It is customary to only image the domain $[0,3]\times[0,1]$ .

Fig. 18 shows the density variable at the final time from a 1920×480 zone simulation using our FD-WENO scheme at seventh order. The smaller panel in the same figure shows a blow-up of the region around the Mach stem. Notice that our scheme resolves all the structures, including the instabilities that develop around the Mach stem. Our fifth order scheme also shows a similar development of instabilities, showing the value of higher order schemes.

## X) Test Problems for Hyperbolic PDE Systems with Stiff Source Terms

FD-WENO schemes show one of their great advantages when stiff source terms are present. This is because these schemes operate pointwise, rather than in a finite volume sense. As a result, the source terms are also enforced at zone centers, which coincides with the location where the primal variables are collocated. Several very popular Runge-Kutta IMEX schemes have been presented in Pareschi and Russo [45] and we reproduce the one that we used in all our examples here. We first write the PDE formally as

$$\partial_t \mathbf{U} = \mathbf{L}(\mathbf{U}) + \mathbf{S}(\mathbf{U}) \tag{10.1}$$

where $\mathbf{L}(\mathbf{U})$ denotes the contribution from the flux terms as well as the non-conservative products and $\mathbf{S}(\mathbf{U})$ denotes the contribution from the stiff source terms. The third order accurate IMEX-SSP3(4,3,3) scheme given by

$$\mathbf{U}^{(1)} = \mathbf{U}^{n} + \alpha\ \Delta t\ \mathbf{S}\left(\mathbf{U}^{(1)}\right) \tag{10.2}$$

$$\mathbf{U}^{(2)} = \mathbf{U}^{n} - \alpha\ \Delta t\ \mathbf{S}\left(\mathbf{U}^{(1)}\right) + \alpha\ \Delta t\ \mathbf{S}\left(\mathbf{U}^{(2)}\right) \tag{10.3}$$

$$\mathbf{U}^{(3)} = \mathbf{U}^{n} + \Delta t\ \mathbf{L}\left(\mathbf{U}^{(2)}\right) + (1-\alpha)\Delta t\ \mathbf{S}\left(\mathbf{U}^{(2)}\right) + \alpha\ \Delta t\ \mathbf{S}\left(\mathbf{U}^{(3)}\right) \tag{10.4}$$

$$\begin{aligned}\mathbf{U}^{(4)} = \mathbf{U}^{n} &+ \frac{\Delta t}{4}\mathbf{L}\left(\mathbf{U}^{(2)}\right) + \frac{\Delta t}{4}\mathbf{L}\left(\mathbf{U}^{(3)}\right) + \beta\ \Delta t\ \mathbf{S}\left(\mathbf{U}^{(1)}\right) + \eta\ \Delta t\ \mathbf{S}\left(\mathbf{U}^{(2)}\right) \\ &+ \left(\frac{1}{2} - \beta - \eta - \alpha\right)\Delta t\ \mathbf{S}\left(\mathbf{U}^{(3)}\right) + \alpha\ \Delta t\ \mathbf{S}\left(\mathbf{U}^{(4)}\right)\end{aligned} \tag{10.5}$$

$$\begin{aligned}\mathbf{U}^{n+1} = \mathbf{U}^{n} &+ \frac{\Delta t}{6}\ \mathbf{L}\left(\mathbf{U}^{(2)}\right) + \frac{\Delta t}{6}\ \mathbf{L}\left(\mathbf{U}^{(3)}\right) + \frac{2\Delta t}{3}\ \mathbf{L}\left(\mathbf{U}^{(4)}\right) \\ &+ \Delta t\left[\frac{1}{6}\ \mathbf{S}\left(\mathbf{U}^{(2)}\right) + \frac{1}{6}\ \mathbf{S}\left(\mathbf{U}^{(3)}\right) + \frac{2}{3}\ \mathbf{S}\left(\mathbf{U}^{(4)}\right)\right]\end{aligned} \tag{10.6}$$

The coefficients in the above equations are given by

$$\alpha = 0.24169426078821 \quad ; \quad \beta = 0.06042356519705 \quad ; \quad \eta = 0.12915286960590 \tag{10.7}$$

This timestepping scheme has three flux evaluations and is, therefore, only third order accurate in time. Without the source terms, it can be rewritten in a format that exactly reproduces the very popular third order SSP-RK3 scheme of Shu and Osher [55].

We show a couple of test problems that are based on the Baer-Nunziato compressible multi-phase flow with stiff source terms, as given in Dumbser and Boscheri [23]. The equations with stiff source terms are given by:-

$$\partial_t(\phi_1\rho_1)+\nabla\cdot(\phi_1\rho_1\mathbf{v}_1)=0$$
$$\partial_t(\phi_1\rho_1\mathbf{v}_1)+\nabla\cdot(\phi_1(\rho_1\mathbf{v}_1\otimes\mathbf{v}_1+\mathbf{I}\,p_1))-P_I\nabla\phi_1=-\lambda(\mathbf{v}_1-\mathbf{v}_2)$$
$$\partial_t(\phi_1\rho_1E_1)+\nabla\cdot(\phi_1\mathbf{v}_1(\rho_1E_1+p_1))+P_I\partial_t\phi_1=-\lambda\mathbf{V}_I\cdot(\mathbf{v}_1-\mathbf{v}_2)$$
$$\partial_t(\phi_2\rho_2)+\nabla\cdot(\phi_2\rho_2\mathbf{v}_2)=0$$
$$\partial_t(\phi_2\rho_2\mathbf{v}_2)+\nabla\cdot(\phi_2(\rho_2\mathbf{v}_2\otimes\mathbf{v}_2+\mathbf{I}\,p_2))-P_I\nabla\phi_2=-\lambda(\mathbf{v}_2-\mathbf{v}_1)$$
$$\partial_t(\phi_2\rho_2E_2)+\nabla\cdot(\phi_2\mathbf{v}_2(\rho_2E_2+p_2))+P_I\partial_t\phi_2=-\lambda\mathbf{V}_I\cdot(\mathbf{v}_2-\mathbf{v}_1)$$
$$\partial_t\phi_1+\mathbf{V}_I\cdot\nabla\phi_1=\mu(p_1-p_2)$$

The stiff terms appear on the right hand sides of the above set of equations and the stiffness is regulated by two parameters $\lambda$ and $\mu$. For the sake of completeness, we present below the two test problems from Dumbser and Boscheri [23].

### X.a) One-Dimensional Riemann Problem with Stiff Source Terms

This first test problem has the range $[-0.5, 0.5]$ and was run on a 400 zone mesh with a CFL of 0.8. The problem uses interphase drag $\lambda=10^3$ and pressure relaxation parameter $\mu=10^2$. These values turn this problem into a test problem with a moderately stiff source term. With these stiff relaxation parameters the problem becomes comparable to a Riemann problem for Euler flow with different values of the ratio of specific heats on either side of the initial discontinuity. It was run to a stopping time of 0.2. The left and right initial conditions were initialized around $x=0$ and are documented in Table XIII. Figs. 19a, 19b and 19c show the solid density, solid x-velocity and solid pressure. We see that the results are comparable to those presented in Dumbser and Boscheri [23]. The reference solution is also shown as the solid line in Fig. 19 and was obtained using a third order scheme with 4000 zones.

**Table XIII gives the left states, right states and final time of the one-dimensional Riemann problem for the Baer-Nunziato compressible multi-phase flow with stiff source terms.**

| | $\rho_1$ | $u_1$ | $p_1$ | $\rho_2$ | $u_2$ | $p_2$ | $\phi_1$ | $t_{end}$ |
|---|---|---|---|---|---|---|---|---|
| | $\gamma_1=1.4,\ \pi_1=0,\ \gamma_2=1.67,\ \pi_2=0$ | | | | | | | |
| **Left:** | 1.0 | 0.0 | 1.0 | 1.0 | 0.0 | 1.0 | 0.99 | 0.2 |
| **Right:** | 0.125 | 0.0 | 0.1 | 0.125 | 0.0 | 0.1 | 0.01 | |

### X.b) Two-Dimensional Riemann Problem with Stiff Source Terms

This second test problem consists of a two dimensional Riemann problem. Such two-dimensional Riemann problems consist of four initial states that cover the four quadrants of a square Cartesian domain. In this problem, the domain spans $[-0.5, 0.5]\times[-0.5, 0.5]$ and uses $400\times 400$ zones. The problem was run with a CFL of 0.4 to a stopping time of 0.15. Here the right-upper quadrant is denoted with a subscript "RU"; the left-upper quadrant is denoted with a subscript "LU"; the left-down quadrant is denoted with a subscript "LD" and the right-down quadrant is denoted with a subscript "RD". The physical values for those states are documented in Table XIV. Here we used $\lambda = 10^5$ and $\mu = 10^2$ which corresponds to rather severely stiff source terms. Fig. 20a shows the solid density, Fig. 20b shows the gas density and Fig. 20c shows the solid volume fraction. We see that the results match the solution in Dumbser and Boscheri [23].

**Table XIV gives the right-upper (RU), the left-upper (LU), left-down (LD) and the right-down (RD) states of the two-dimensional Riemann problem for the Baer-Nunziato compressible multi-phase flow with stiff source terms.**

| | $\rho_1$ | $u_1$ | $v_1$ | $p_1$ | $\rho_2$ | $u_2$ | $v_2$ | $p_2$ | $\phi_1$ |
|---|---|---|---|---|---|---|---|---|---|
| | $\gamma_1 = 1.4,\quad \pi_1 = 0,\quad \gamma_2 = 1.67,\quad \pi_2 = 0$ | | | | | | | | |
| **RU:** | 2.0 | 0.0 | 0.0 | 2.0 | 1.5 | 0.0 | 0.0 | 2.0 | 0.8 |
| **LU:** | 1.0 | 0.0 | 0.0 | 1.0 | 0.5 | 0.0 | 0.0 | 1.0 | 0.4 |
| **LD:** | 2.0 | 0.0 | 0.0 | 2.0 | 1.5 | 0.0 | 0.0 | 2.0 | 0.8 |
| **RD:** | 1.0 | 0.0 | 0.0 | 1.0 | 0.5 | 0.0 | 0.0 | 1.0 | 0.4 |

## XI) Conclusions

Robust and efficient finite difference WENO schemes are very valuable in science and engineering because they allow us to carry out high order accurate multidimensional simulations with a speed and efficiency that cannot be matched by finite volume WENO or DG schemes. Because of this reason, they occupy a unique place in computational science. Robust FD-WENO schemes have been presented in the literature by Jiang and Shu [40] and Balsara and Shu [3] and these classical FD-WENO schemes have been refined by many others since. It is also noteworthy that closed form expressions for the smoothness indicators that are expressed as a sum of perfect squares have been available for FD-WENO in Balsara, Garain and Shu [6]. That work has also been extended to the smoothness indicators for finite volume WENO schemes in the Supplement to Balsara *et al*. [8]. While such advances bring a modicum of completeness to the WENO algorithms available to us, a very substantial deficiency has been the absence (to the best of our knowledge) of FD-WENO algorithms for hyperbolic systems with non-conservative products. The fact that many recently discovered PDE systems have such non-conservative products, as well as having stiff source terms, makes this a problem for which a good resolution is sorely desired.

In this paper we first analyze the issues associated with designing a FD-WENO scheme that could accommodate hyperbolic PDEs with non-conservative products (Section II). By focusing on the HLLI Riemann solver of Dumbser and Balsara [26] (Section III), we present the final FD-WENO scheme in Section IV. The use of the HLLI Riemann solver enables us to obtain schemes that have a good supersonic limit. We have also been able to obtain good results with the LLF variant of the same Riemann solver. The availability of an anti-diffusive flux ensures that both HLL and LLF variants of the FD-WENO schemes presented here can accurately capture stationary linearly degenerate discontinuities on the mesh. Section V provides the analysis that enables us to understand why this scheme is a natural higher order extension of the work in Dumbser and Balsara [26]. Section VI describes the pointwise implementation of the scheme presented here. Because of their use of pointwise primal variables, the treatment of stiff source terms is also very easy for FD-WENO schemes, as shown in Section X.

Our new FD-WENO scheme uses the same reconstruction strategy that has been used in classical FD-WENO schemes, making it very easy for practitioners to transition over to the new formulation. Our present FD-WENO formulation also uses only one reconstruction step as opposed to two reconstruction steps (for each of the right-going and left-going fluxes) that were needed in LLF versions of the classical FD-WENO. Along with the closed form expressions for WENO reconstruction in Balsara, Garain and Shu [6], the reduced reconstruction steps in our new formulation should make for substantially more efficient FD-WENO schemes.

The new FD-WENO is shown to perform as well as the classical version of FD-WENO, with two major advantages:- 1) It can capture jumps in stationary linearly degenerate wave families exactly. 2) It only requires the WENO reconstruction to be applied once. To highlight the versatility of the new FD-WENO method, we have focused on three major hyperbolic systems with non-conservative products: a) The Baer-Nunziato model for compressible multiphase flow, b) The multiphase debris flow model of Pitman and Le and c) The two-layer shallow water equations. To show that the method is perfectly general, we also present several stringent examples from Euler flow. A very extensive list of one-dimensional and multidimensional test problems are also presented here. Because these test problems could also help others, some of them are concatenated in detail in the Appendices of this paper. Section VII shows numerous non-trivial, multidimensional accuracy studies involving conservative systems as well as many hyperbolic systems with non-conservative products. We show that our FD-WENO methods meet their design accuracies for numerous smooth multidimensional test problems. Section VIII shows several one-dimensional test problems. It also demonstrates the ability of the method to capture isolated stationary discontinuities and shows that the method performs well on test problems that require well-balancing. Section IX presents numerous multidimensional test problems, many of which are genuinely novel to this field. We anticipate that they will serve as useful multidimensional test problems in future work. Section X shows several test problems involving stiff source terms.

## Acknowledgements

DSB acknowledges support via NSF grants NSF-19-04774, NSF-AST-2009776, NASA-2020-1241 and NASA grant 80NSSC22K0628. DSB and HK acknowledge support from a Vajra award, VJR/2018/00129 and also a travel grant from Notre Dame International. CWS acknowledges support via AFOSR grant FA9550-20-1-0055 and NSF grant DMS-2010107.

## Ethical Statement

**i. Compliance with Ethical Standards** : This manuscript complies with all ethical standards for scientific publishing.

**ii. (in case of Funding) Funding** : The funding has been acknowledged. DSB acknowledges support via NSF grants NSF-19-04774, NSF-AST-2009776, NASA-2020-1241 and NASA-80NSSC22K0628. DSB and HK acknowledge support from a Vajra award, VJR/2018/00129. CWS acknowledges support via AFOSR grant FA9550-20-1-0055 and NSF grant DMS-2010107.

**iii. Conflict of Interest** : On behalf of all authors, the corresponding author states that there is no conflict of interest.

**iv. Ethical approval** : N/A

**v. Informed consent** : N/A

**Appendix A) Description of the 2D Vortex for the Baer-Nunziato equations**

In this Appendix, we give description of the 2D Vortex problem for the Baer-Nunziato equations which has been proposed in Dumbser *et al*. [22]. We consider the same set up here by considering a smooth unsteady vortex. The exact solution of this smooth unsteady vortex is achieved in two steps. First, an exact stationary vortex is constructed using the polar form of the Baer-Nunziato equations. Second, the two-dimensional vortex attains unsteady motion by superimposing constant velocity field $\bar{\mathbf{v}} = (\bar{u}, \bar{v})$. The exact solution is then given by the advection of the initial conditions, advected with the speed $\bar{\mathbf{v}}$.

In the polar coordinates ($r-\theta$), we denote the angular velocities by $v_k^\theta$ and the radial velocities by $v_k^r$, where the subscript $k=1$ describes the solid phase velocities and $k=2$ describes the gas phase velocities. The radial velocities are set to zero. Following Dumbser *et al*. [22], particular expressions for the densities $\rho_k$, pressures $p_k$ and solid volume fraction $\phi_1$ are given by

$$
\begin{aligned}
&\rho_1 = 1, \quad \rho_2 = 2, \\
&p_k = p_{k0}\left(1 - \frac{1}{4}\mathrm{EXP}\left[1 - r'^2 / s_k^2\right]\right), \qquad (k = 1, 2), \\
&\phi_1 = \frac{1}{3} + \frac{1}{2\sqrt{2\pi}}\mathrm{EXP}\left[-r'^2 / 2\right],
\end{aligned}
$$

where $r' \equiv r / r_s$. Here $r_s$ gives us the length scale of the vortex. For the vortex that is used in Section VII.b, the constants $r_s$, $p_{k0}$ and $s_k$ are specified by the following values

$$
r_s = 1, \quad p_{10} = 1, \quad p_{20} = \frac{3}{2}, \quad s_1 = \frac{3}{2}, \quad s_2 = \frac{7}{5}.
$$

We have used the implementation friendly notation “EXP” to denote the exponential function. Once the particular values for the densities $\rho_k$, pressures $p_k$ and solid volume fraction $\phi_1$ are known, the angular velocities are derived from the following momentum equations

$$\frac{\partial}{\partial r}(\phi_1 p_1) = p_2 \frac{\partial}{\partial r}\phi_1 + \frac{1}{r}(v_1^\theta)^2 \phi_1 \rho_1,$$
$$\frac{\partial}{\partial r}(\phi_2 p_2) = p_2 \frac{\partial}{\partial r}\phi_2 + \frac{1}{r}(v_2^\theta)^2 \phi_2 \rho_2. \tag{A.1}$$

The system (A.1) is a simple algebraic system, which can be solved to obtained the following expressions for the steady state angular velocities $v_k^\theta$

$$v_1^\theta = \sqrt{\frac{r}{\phi_1 \rho_1}\left(\phi_1 A_1 + (p_1 - p_2) B_1\right)},$$
$$v_2^\theta = \sqrt{\frac{r}{\rho_2} A_2},$$
where $A_k \equiv \frac{\partial p_k}{\partial r} = (p_{k0} - p_k)\frac{2r'}{r_s s_k^2}$ and $B_1 \equiv \frac{\partial \phi_1}{\partial r} = \frac{r'}{r_s}\left(\frac{1}{3} - \phi_1\right)$.

For the vortex that is used in Section VII.b, the particular values for the remaining parameters are given by

$$\bar{u} = 2, \quad \bar{v} = 2, \quad \gamma_1 = 1.4, \quad \pi_1 = 0, \quad \gamma_2 = 1.35, \quad \pi_2 = 0.$$

Note that the transformation $x = r\cos(\theta),\ y = r\sin(\theta)$, where $r = \sqrt{x^2 + y^2}$ and $\theta := \theta(x, y)$ is the angular coordinate, can be used to convert the polar coordinates in to the Cartesian coordinates. Indeed the initial velocities for the unsteady case in the Cartesian coordinate system are given by

$$u_1 = \bar{u} - v_1^\theta \sin(\theta), \quad v_1 = \bar{v} + v_1^\theta \cos(\theta),$$
$$u_2 = \bar{u} - v_2^\theta \sin(\theta), \quad v_2 = \bar{v} + v_2^\theta \cos(\theta).$$

For the vortex that is used in Section VII.b, after one-advection period $(t = 5)$, the exact solution is given by the initial condition. The test is run on a two-dimensional computational domain $\Omega = [-5,5]\times[-5,5]$ with the periodic boundary conditions at all four boundaries. The stopping time (one-advection period) is $t = 5$. For the fifth, seventh and ninth order schemes we double the computational domain and stopping time to minimize the effect of small jumps in the velocity field at the periodic boundaries.

**Appendix B) Description of the 2D Vortex for the multiphase debris flow equations**

In this Appendix, we give description of the 2D Vortex problem for the multiphase debris flow equations which has been described in Dumbser *et al.* [21]. We consider the same initial data. The initial heights for the solid and fluid flows are given by the following expressions

$$h_s = \frac{1}{4}\frac{v_{10}^2 s_2 \mathrm{EXP}\left[2s_1\right](1-\mathrm{EXP}\left[-2s_1 r'^2\right]) - v_{20}^2 s_1 \mathrm{EXP}\left[2s_2\right](1-\mathrm{EXP}\left[-2s_2 r'^2\right]) + 4h_{10} g s_1 s_2 (1-\rho)}{g(1-\rho)s_1 s_2},$$
$$h_f = \frac{1}{4}\frac{v_{20}^2 s_1 \mathrm{EXP}\left[2s_2\right](1-\mathrm{EXP}\left[-2s_2 r'^2\right]) - \rho v_{10}^2 s_2 \mathrm{EXP}\left[2s_1\right](1-\mathrm{EXP}\left[-2s_1 r'^2\right]) + 4h_{20} g s_1 s_2 (1-\rho)}{g(1-\rho)s_1 s_2},$$

where $r' \equiv r / r_s$. For the vortex that is used in Section VII.c, we have used $g = 10$, $\rho = \frac{9}{10}$, $s_1 = \frac{1}{2}$, $s_2 = 1$, $v_{10} = \frac{3}{4}$, $v_{20} = \frac{1}{10}$, $h_{10} = 1$ and $h_{20} = 1$. Also for the vortex used in Section VII.c, we set the scaling factor $r_s = 1$. The initial bottom topography is set to zero ($b = 0$). We have used the implementation friendly notation "EXP" to denote the exponential function.

In the polar coordinates ($r - \theta$), the angular velocities are denoted by $v_k^\theta$ and the radial velocities are denoted by $v_k^r$, where the subscript $k = s$ refers to the solid flow velocities and $k = f$ refers to the fluid flow velocities. The radial velocities are set to zero. Following Dumbser *et al.* [21], solid and fluid angular velocities for the steady vortex are given by

$$v_s^\theta = \sqrt{r' g \frac{\partial}{\partial r'}\left(h_s + \frac{1+\rho}{2} h_f\right) + r' \frac{1-\rho}{2} g \frac{h_f}{h_s} \frac{\partial}{\partial r'} h_s},$$
$$v_f^\theta = r' v_{10} \mathrm{EXP}\left[s_1 (1 - r'^2)\right].$$

Using the polar transformation $x = r\cos(\theta)$, $y = r\sin(\theta)$, where $r = \sqrt{x^2 + y^2}$ and $\theta := \theta(x, y)$ is the angular coordinate, the initial solid and fluid velocities, after applying the Galilean transformation are given by

$$u_s = \bar{u} - v_s^\theta \sin(\theta), \quad v_s = \bar{v} + v_s^\theta \cos(\theta),$$
$$u_f = \bar{u} - v_f^\theta \sin(\theta), \quad v_f = \bar{v} + v_f^\theta \cos(\theta).$$

For the vortex that is used in Section VII.c, we take $\bar{u} = \bar{v} = 5$. After one-advection period, the exact solution is given by the initial condition. For the vortex that is used in Section VII.c, the test is run on a two-dimensional computational domain $\Omega = [-5,5]\times[-5,5]$ with the periodic boundary conditions at all four boundaries. The stopping time (one-advection period) is $t = 2.$ To minimize the effect of small jumps in the velocity field at the periodic boundaries we double the computational domain and stopping time for the fifth, seventh and ninth order schemes.

**Appendix C) Description of the 2D Vortex for the two-layer shallow water equations.**

In this last Appendix, we give description of the 2D Vortex problem for the two-layer shallow water equations which has been proposed in Dumbser *et al*. [21]. We consider the same set up here. The initial heights for the upper and lower fluid, in terms of the polar coordinates, are given by the following expressions

$$h_1 = \frac{1}{4}\frac{v_{10}^2 s_2 \mathrm{EXP}\left[2s_1\right](1-\mathrm{EXP}\left[-2s_1 r'^2\right]) - v_{20}^2 s_1 \mathrm{EXP}\left[2s_2\right](1-\mathrm{EXP}\left[-2s_2 r'^2\right]) + 4h_{10} g s_1 s_2 (1-\rho)}{g(1-\rho)s_1 s_2},$$
$$h_2 = \frac{1}{4}\frac{v_{20}^2 s_1 \mathrm{EXP}\left[2s_2\right](1-\mathrm{EXP}\left[-2s_2 r'^2\right]) - \rho v_{10}^2 s_2 \mathrm{EXP}\left[2s_1\right](1-\mathrm{EXP}\left[-2s_1 r'^2\right]) + 4h_{20} g s_1 s_2 (1-\rho)}{g(1-\rho)s_1 s_2},$$

where $r' \equiv r / r_s$ For the vortex used in Section VII.d, we have used $g = 10$, $\rho = \frac{9}{10}, s_1 = \frac{1}{2}$, $s_2 = 1, v_{10} = \frac{3}{4}, v_{20} = \frac{1}{10}, h_{10} = 1,$ $h_{20} = 1$ and we set the scaling factor $r_s = 1$. We remark that the above expressions are exactly the same expressions which are used in Appendix B to define the solid and fluid heights for the multi-phase debris flow equations. The initial bottom topography is set to zero ($b = 0$).

Following Dumbser *et al*. [21], the radial velocities are set to zero and the angular velocities for the steady vortex are given by

$$v_1^\theta = r' v_{10} \mathrm{EXP}\left[s_1(1-r'^2)\right],$$
$$v_2^\theta = r' v_{20} \mathrm{EXP}\left[s_2(1-r'^2)\right].$$

The steady vortex is superimposed with a uniform velocity field $(\bar{u},\bar{v})$ to obtain the following initial velocities for the unsteady two-dimensional vortex:

$$u_1 = \bar{u} - v_1^{\theta}\sin(\theta), \quad v_1 = \bar{v} + v_1^{\theta}\cos(\theta),$$
$$u_2 = \bar{u} - v_2^{\theta}\sin(\theta), \quad v_2 = \bar{v} + v_2^{\theta}\cos(\theta).$$

For the vortex that is used in Section VII.d, we take mean velocities $\bar{u} = \bar{v} = 5$. After one advection period, the exact solution is given by the initial condition. For the vortex used in Section VII.d, the test is run on a two-dimensional computational domain $\Omega = [-5,5]\times[-5,5]$ with the periodic boundary conditions at all four boundaries. The stopping time is $t = 2$. We remark that at the periodic boundaries the velocity field may not exactly vanish to zero. Therefore, to reduce the effect of small fluctuations in the velocity field we double the computational domain and stopping time for the higher order schemes (fifth, seventh and ninth order schemes in our case).

**Appendix D) Algorithmic Touch-ups to the WENO-AO and Multiresolution WENO Algorithms**

WENO-AO is discussed in Balsara Garain and Shu [6] where we write their eqns. (3.6) and (3.7a) as

$$\tau = \frac{1}{3}\left(\left|\beta_3^{r5} - \beta_1^{r3}\right| + \left|\beta_3^{r5} - \beta_2^{r3}\right| + \left|\beta_3^{r5} - \beta_3^{r3}\right|\right) \tag{D.1}$$

and

$$\begin{aligned} w_3^{r5} &= \gamma_3^{r5}\left(1+\tau^2\Big/\left(\beta_3^{r5}+\varepsilon\right)^2\right) \quad ; \quad w_1^{r3} = \gamma_1^{r3}\left(1+\tau^2\Big/\left(\beta_1^{r3}+\varepsilon\right)^2\right) \quad ; \\ w_2^{r3} &= \gamma_2^{r3}\left(1+\tau^2\Big/\left(\beta_2^{r3}+\varepsilon\right)^2\right) \quad ; \quad w_3^{r3} = \gamma_3^{r3}\left(1+\tau^2\Big/\left(\beta_3^{r3}+\varepsilon\right)^2\right) \end{aligned} \tag{D.2}$$

For WENO-AO(5,3), the above two equations are entirely correct. For smooth flow, the non-linear weights differ from the linear weights by terms that are proportional to $O\left(\Delta x^5\right)$. However, the above formulae cannot be extended as-is to seventh and ninth orders. Instead, for the seventh order WENO-AO(7,3) we should write

$$\tau = \frac{1}{3}\left(\left|\beta_4^{r7} - \beta_1^{r3}\right| + \left|\beta_4^{r7} - \beta_2^{r3}\right| + \left|\beta_4^{r7} - \beta_3^{r3}\right|\right) \tag{D.3}$$

and

$$w_4^{r7} = \gamma_4^{r7}\left(1+\tau^3\Big/\left(\beta_4^{r7}+\varepsilon\right)^2\right) \quad ; \quad w_1^{r3} = \gamma_1^{r3}\left(1+\tau^3\Big/\left(\beta_1^{r3}+\varepsilon\right)^2\right) \quad ;$$
$$w_2^{r3} = \gamma_2^{r3}\left(1+\tau^3\Big/\left(\beta_2^{r3}+\varepsilon\right)^2\right) \quad ; \quad w_3^{r3} = \gamma_3^{r3}\left(1+\tau^3\Big/\left(\beta_3^{r3}+\varepsilon\right)^2\right) \qquad \text{(D.4)}$$

With the change in eqn. (D.4), for very smooth flow, the non-linear weights differ from the linear weights by $O\left(\Delta x^7\right)$. As a result, a good order property is preserved. Likewise, for the ninth order WENO-AO(9,3) we should write

$$\tau = \frac{1}{3}\left(\left|\beta_5^{r9} - \beta_1^{r3}\right| + \left|\beta_5^{r9} - \beta_2^{r3}\right| + \left|\beta_5^{r9} - \beta_3^{r3}\right|\right) \qquad \text{(D.5)}$$

and

$$w_5^{r9} = \gamma_5^{r9}\left(1+\tau^4\Big/\left(\beta_5^{r9}+\varepsilon\right)^2\right) \quad ; \quad w_1^{r3} = \gamma_1^{r3}\left(1+\tau^4\Big/\left(\beta_1^{r3}+\varepsilon\right)^2\right) \quad ;$$
$$w_2^{r3} = \gamma_2^{r3}\left(1+\tau^4\Big/\left(\beta_2^{r3}+\varepsilon\right)^2\right) \quad ; \quad w_3^{r3} = \gamma_3^{r3}\left(1+\tau^4\Big/\left(\beta_3^{r3}+\varepsilon\right)^2\right) \qquad \text{(D.6)}$$

With the change in eqn. (D.6), for very smooth flow, the non-linear weights differ from the linear weights by $O\left(\Delta x^9\right)$. Notice too from comparing eqns. (D.2), (D.4) and (D.6) that the power to which "$\tau$" is raised increases with the order of accuracy. That is the essential change from Balsara Garain and Shu [6].

Zhu and Shu [64] also design a scheme which can selectively reduce the order of accuracy of its interpolation. They present an approach that can graciously lower its order of accuracy on a sequence of nested, centered stencils if a loss of smoothness is detected in one of the stencils. This process can continue until it reaches a first order polynomial $q_1\left(x\right)$ in their eqn. (2.5) which is indeed first order accurate. When we tried the third order multiresolution scheme, we found that on somewhat large meshed it degrades to less than second order of accuracy. The problem stems from the fact that the scheme can use either a third order reconstruction or it can degrade to first order, with nothing in between. Eqns. (2.16) to (2.19) of Zhu and Shu [64] provide a recipe for order reduction to first order and we leave them unchanged. However, we recommend using either

$$q_1\left(x\right) = \overline{w}_i + MC_\beta\left(\overline{w}_{i+1} - \overline{w}_i, \overline{w}_i - \overline{w}_{i-1}\right) \qquad \text{(D.7)}$$

or

$$q_1(x) = \bar{w}_i + vanAlbada_\varepsilon\left(\bar{w}_{i+1} - \bar{w}_i, \bar{w}_i - \bar{w}_{i-1}\right) \ . \tag{D.8}$$

The limiters used in the above two equations are explicitly given as the well-known MC limiter (with a minor modification stemming from the use of "$\beta$"):-

$$MC_\beta(a,b) \equiv 0.5\left(\text{sgn}(a) + \text{sgn}(b)\right)\min\left(|a+b|/2, \beta|a|, \beta|b|\right) \quad \text{with} \quad 1 \le \beta \le 2 \tag{D.9}$$

And the van Albada limiter (with a minor modification stemming from the use of "$\varepsilon$"):-

$$vanAlbada_\varepsilon(a,b) \equiv \frac{\left(b^2+\varepsilon\right)a + \left(a^2+\varepsilon\right)b}{a^2+b^2+2\varepsilon} \ . \tag{D.10}$$

Eqn. (D.10) could prove to be useful in implicit problems because of its continuous variation with respect to the left and right slopes. In eqn. (D.10), "$\varepsilon$" is a very small user-defined number. The use of a TVD limited solution in eqns. (D.7) or (D.8) allows the third order Multiresolution WENO reconstruction to indeed converge with third order accuracy when the solution is smooth enough on the mesh. When the solution becomes slightly non-smooth, the reconstruction has the option to gradually transition from third to second order. With this change, the solution is only reconstructed with first order accuracy when it is completely non-smooth.

## Figure Captions

*Fig. 1 shows part of the mesh around zone "i". The mesh functions are collocated at the zone centers, as shown by the thick dots. The zone boundaries are shown by the vertical lines. The figure also shows the stencils associated with the zone "i" for the fifth order WENO-AO reconstruction. We have three smaller third order stencils and a large fifth order stencil. The reconstructed variables at the zone boundaries are shown with a caret. The variables with a superscript star are resolved states obtained by the pointwise application of a simple HLL or LLF Riemann solver at the zone boundaries.*

*Fig. 2) Debris Flow: Jump in the linearly degenerate fields using the $3^{rd}$ order accurate HLL-based FD-WENO scheme with 200 zones. The $3^{rd}$ order LLF-based FD-WENO also shows identical results. Identical results were also obtained for $5^{th}$ and $7^{th}$ order HLL-based and LLF-based FD-WENO schemes.*

*Fig. 3a) Debris Flow: Riemann problem-2 (top) and -3 (bottom) using the $3^{rd}$ order accurate HLL-based FD-WENO scheme with 200 zones.*

*Fig. 3b) Debris Flow: Riemann problem-2 (top) and -3 (bottom) using the $5^{th}$ order accurate HLL-based FD-WENO scheme with 200 zones.*

*Fig. 3c) Debris Flow: Riemann problem-2 (top) and -3 (bottom) using the $7^{th}$ order accurate HLL-based FD-WENO scheme with 200 zones.*

*Fig. 4) Debris Flow: Non-constant Bathymetry test using the $3^{rd}$ order (left column), $5^{th}$ order (middle column) and $7^{th}$ order (right column) accurate HLL-based FD-WENO scheme with 400 zones. The top three panels show the solution at time of $t=0.5$; the lower three panels show the solution at a time of $t=1.25$.*

*Fig. 5a) Baer-Nunziato: Abgrall problem using the 3$^{rd}$ order (left column), 5$^{th}$ order (middle column) and 7$^{th}$ order (right column) accurate HLL-based FD-WENO scheme with 200 zones.*

*Fig. 5b) Baer-Nunziato: Abgrall problem using the 3$^{rd}$ , 5$^{th}$ and 7$^{th}$ order accurate LLF-based FD-WENO scheme with 200 zones.*

*Fig. 6a) Baer-Nunziato: RP-1 (left), RP-2 (middle) and RP-3 (right) using the 3$^{rd}$ order accurate HLL-based FD-WENO scheme with 200 zones.*

*Fig. 6b) Baer-Nunziato: RP-4 (left), RP-5 (middle) and RP-6 (right) using the 3$^{rd}$ order accurate HLL-based FD-WENO scheme with 200 zones.*

*Fig. 6c) Baer-Nunziato: RP-1 (left), RP-2 (middle) and RP-3 (right) using the 5$^{th}$ order accurate HLL-based FD-WENO scheme with 200 zones.*

*Fig. 6d) Baer-Nunziato: RP-4 (left), RP-5 (middle) and RP-6 (right) using the 5$^{th}$ order accurate HLL-based FD-WENO scheme with 200 zones.*

*Fig. 6e) Baer-Nunziato: RP-1 (left), RP-2 (middle) and RP-3 (right) using the 7$^{th}$ order accurate HLL-based FD-WENO scheme with 200 zones.*

*Fig. 6f) Baer-Nunziato: RP-4 (left), RP-5 (middle) and RP-6 (right) using the 7$^{th}$ order accurate HLL-based FD-WENO scheme with 200 zones.*

*Fig. 6g) Verification of mass conservation in the solid and gas masses for the Riemann problem-2 (left panel), Riemann problem-4 (middle panel) and Riemann problem-5 (right panel) using the 3$^{rd}$, 5$^{th}$ and 7$^{th}$ order HLL-based FD-WENO schemes. Conservation of solid and gas masses is respected to at least one part in 10$^{12}$ .*

*Fig. 7) Two-Layer Shallow water: Riemann problem-1 using the 3$^{rd}$ order accurate HLL-based FD-WENO scheme with 200 zones. The 3$^{rd}$ order LLF-based FD-WENO also shows identical results. Identical results were also obtained for 5$^{th}$ and 7$^{th}$ order HLL-based FD-WENO and LLF-based FD-WENO schemes.*

*Fig. 8a) Two-Layer Shallow water: Riemann problem-2(left) and -3(right) using the $3^{rd}$ order accurate HLL-based FD-WENO scheme with 200 zones.*

*Fig. 8b) Two-Layer Shallow water: Riemann problem-2(left) and -3(right) using the $5^{th}$ order accurate HLL-based FD-WENO scheme with 200 zones.*

*Fig. 8c) Two-Layer Shallow water: Riemann problem-2(left) and -3(right) using the $7^{th}$ order accurate HLL-based FD-WENO scheme with 200 zones.*

*Fig. 9a) Two-Layer Shallow Water Model: Non-constant Bathymetry test. Left: Large perturbation, Right: Small perturbation using the $3^{rd}$ order accurate HLL-based FD-WENO scheme with 400 zones. The $3^{rd}$ order LLF-based FD-WENO scheme produces identical results.*

*Fig. 9b) Two-Layer Shallow Water Model: Non-constant Bathymetry test. Left: Large perturbation, Right: Small perturbation using the $5^{th}$ order accurate HLL-based FD-WENO scheme with 400 zones. The $5^{th}$ order LLF-based FD-WENO scheme produces identical results.*

*Fig. 9c) Two-Layer Shallow Water Model: Non-constant Bathymetry test. Left: Large perturbation, Right: Small perturbation using the $7^{th}$ order accurate HLL-based FD-WENO scheme with 400 zones. The $7^{th}$ order LLF-based FD-WENO scheme produces identical results.*

*Fig. 10) EULER: Blast wave interaction using the $3^{rd}$ , $5^{th}$ and $7^{th}$ order accurate HLL-based FD-WENO scheme with 1000 zones. Figs. 10a, 10b and 10c individually show the density at $3^{rd}$ , $5^{th}$ and $7^{th}$ orders.*

*Fig. 11a, b, c) Baer Nunziato: Shock-bubble interaction problem using the $3^{rd}$, $5^{th}$ and $7^{th}$ order accurate HLL-based FD-WENO schemes with 700 $\times$300 zones. The solid density profiles has been shown.*

*Fig. 12) Baer-Nunziato: Shock-Vortex Interaction using the $7^{th}$ order accurate HLL-based FD-WENO scheme with 600 $\times$600 zones at time levels t=0.0, 0.23 and 0.84. Figs. 12a, 12b and 12c show solid volume fraction at times t=0.0, 0.23, 0.84. Figs. 12d, 12e, 12f show solid x-velocity at times t=0.0, 0.23, 0.84. For the solid volume fraction, 30 contours were fit between a range of 0.33 and 0.530. For the solid x-velocity, 30 contours were fit between a range of -0.5 and 1.95.*

*Fig. 13a, b, c) Two-Layer Shallow Water Model : Shock-bubble interaction problem using the $3^{rd}$, $5^{th}$ and $7^{th}$ order accurate HLL-based FD-WENO schemes with 700 $\times$300 zones. Height of the lower fluid has been shown.*

*Fig. 14) Two-Layer Shallow Water Model: Shock-Vortex Interaction using the $7^{th}$ order accurate HLL-based FD-WENO scheme with 600 $\times$600 zones at time levels t=0.0, 0.06 and 0.24. Figs. 14a, 14b and 14c show height of the upper fluid at times t=0.0, 0.06, 0.24. Figs. 14d, 14e, 14f show x-velocity of the upper fluid at times t=0.0, 0.06, 0.24. For the height, 40 contours were fit between a range of 1.0 and 3.5. For the velocity, 40 contours were fit between a range of 1.8 and 5.6.*

*Fig. 15a, b, c) Debris Flow : Shock-bubble interaction problem using the $3^{rd}$, $5^{th}$ and $7^{th}$ order accurate HLL-based FD-WENO schemes with 700 $\times$500 zones. Height of the fluid phase has been shown.*

*Fig. 16) Debris Flow: Shock-Vortex Interaction using the $7^{th}$ order accurate HLL-based FD-WENO scheme with 600 $\times$600 zones at time levels t=0.0, 0.06 and 0.24. Figs. 16a, 16b and 16c show the solid height at times t=0.0, 0.06, 0.24. Figs. 16d, 16e, 16f show the solid x-velocity at times t=0.0, 0.06, 0.24. For the height, 40 contours were fit between a range of 1.0 and 3.0. For the velocity, 40 contours were fit between a range of 1.5 and 5.6.*

*Fig. 17) EULER: Forward facing step problem using the 7<sup>th</sup> order accurate HLL-based FD-WENO scheme with 1440 × 480 zones. 30 contours were fit between a range of 0.105 and 6.699.*

*Fig. 18) EULER: Double Mach reflection problem using the 7<sup>th</sup> order accurate HLL-based FD-WENO scheme with 1920 × 480 zones. Fig. 18a shows the density profile at the 7<sup>th</sup> orders. Fig. 18b shows the detailed view of the density profile. 20 contours were fit between a range of 1.0 and 20.97.*

*Fig. 19) Baer-Nunziato with stiff source: Results for the one-dimensional Riemann Problem using the 7<sup>th</sup> order accurate HLL-based FD-WENO scheme with 400 zones. Figs. 19a, 19b and 19c show the solid density, solid x-velocity and solid pressure.*

*Fig. 20) Baer-Nunziato with stiff source: Results for the two-dimensional Riemann Problem using the 7<sup>th</sup> order accurate HLL-based FD-WENO scheme with 400 × 400 zones. Fig. 20a shows the solid density, Fig. 20b shows the gas density and Fig. 20c shows the solid volume fraction. 30 equidistant contour lines are shown over the color plots.*

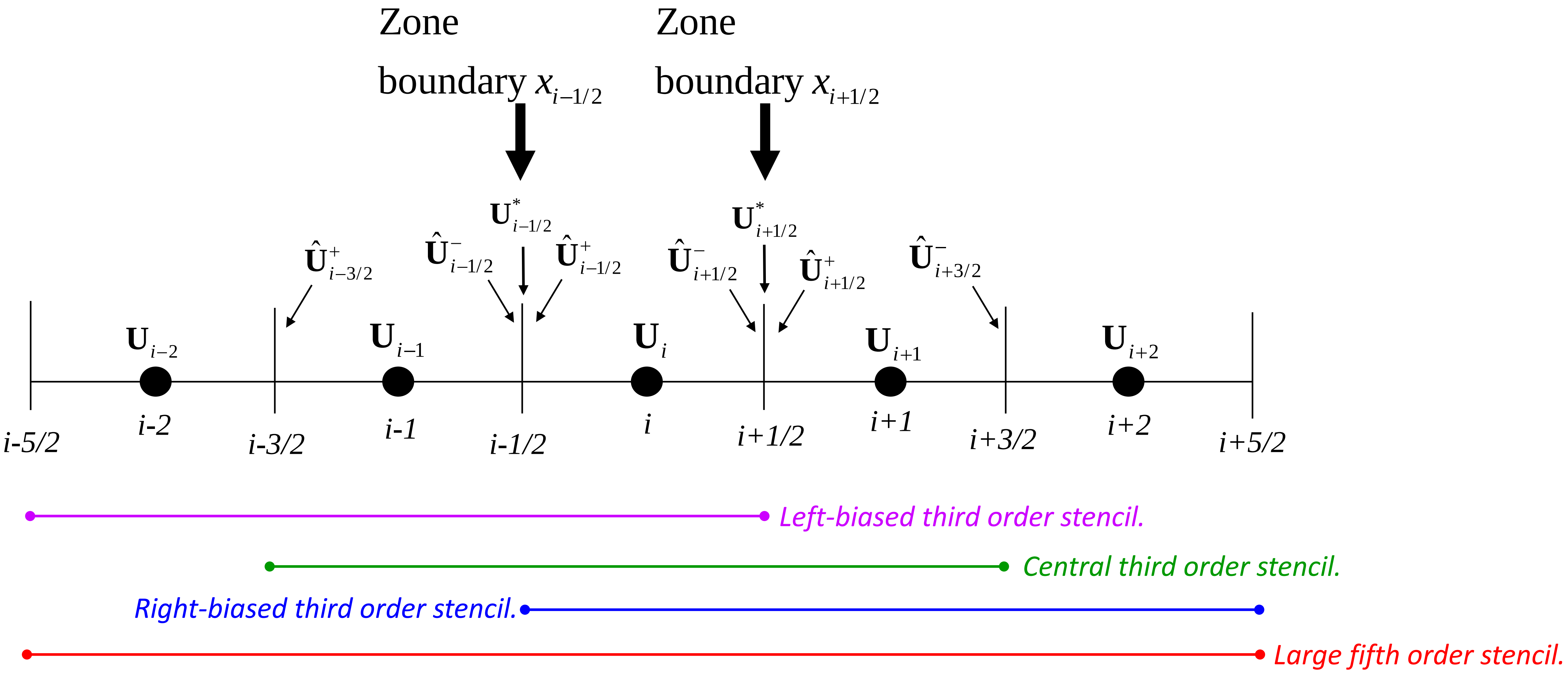

*Fig. 1 shows part of the mesh around zone "i". The mesh functions are collocated at the zone centers, as shown by the thick dots. The zone boundaries are shown by the vertical lines. The figure also shows the stencils associated with the zone "i" for the fifth order WENO-AO reconstruction. We have three smaller third order stencils and a large fifth order stencil. The reconstructed variables at the zone boundaries are shown with a caret. The variables with a superscript star are resolved states obtained by the pointwise application of a simple HLL or LLF Riemann solver at the zone boundaries.*

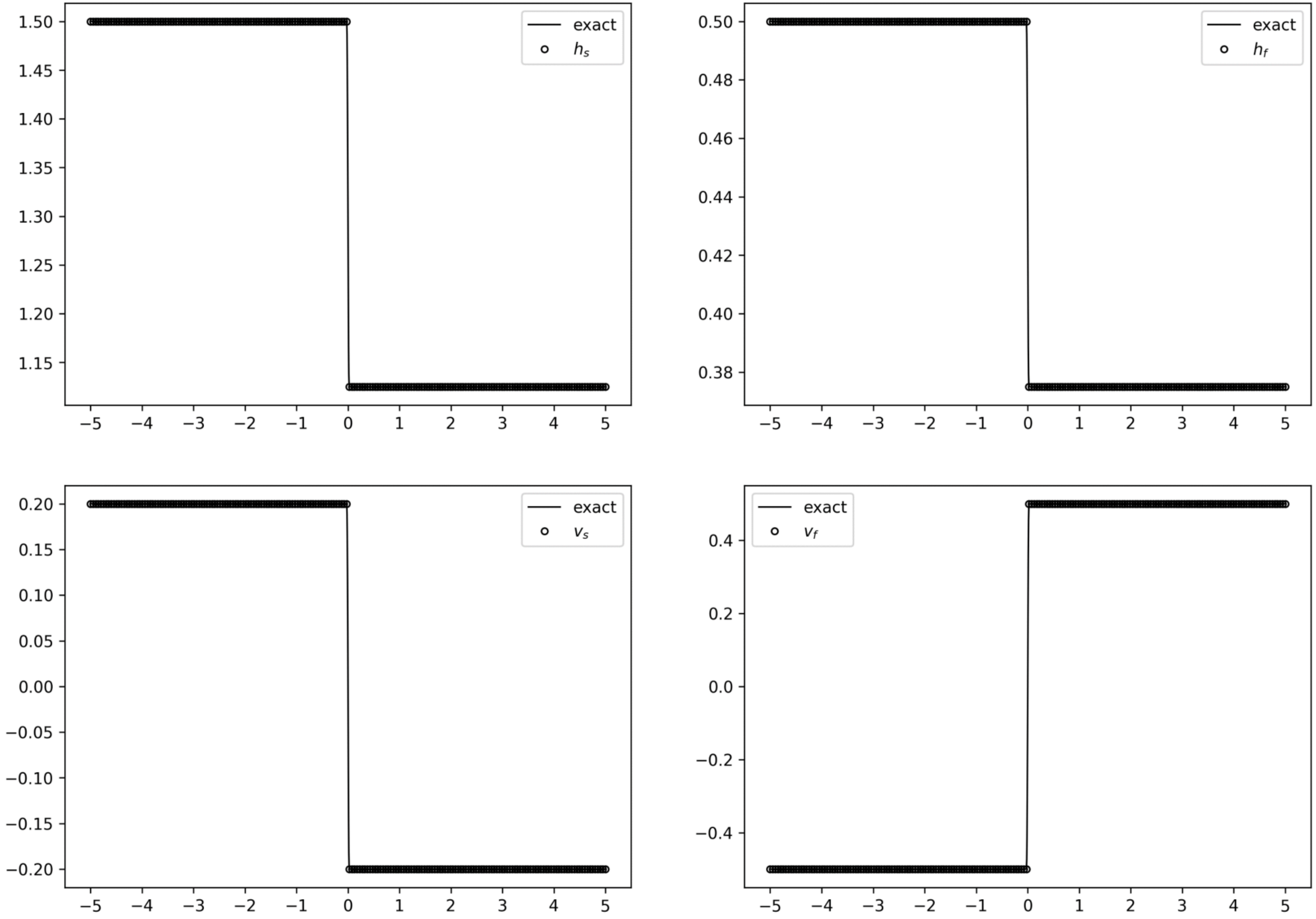

*Fig. 2) Debris Flow: Jump in the linearly degenerate fields using the 3rd order accurate HLL-based FD-WENO scheme with 200 zones. The 3rd order LLF-based FD-WENO also shows identical results. Identical results were also obtained for 5th and 7th order HLL-based and LLF-based FD-WENO schemes.*

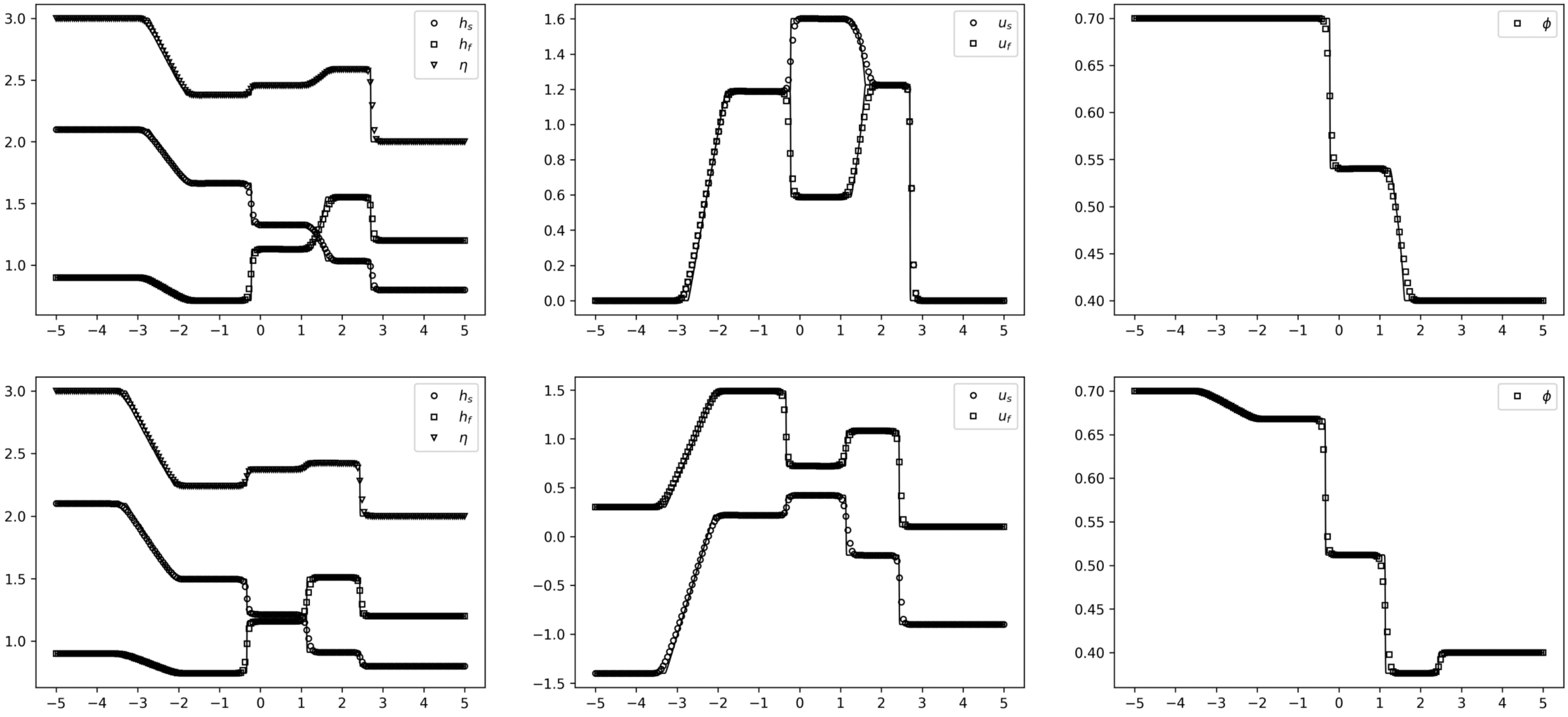

*Fig. 3a) Debris Flow: Riemann problem-2 (top) and -3 (bottom) using the 3rd order accurate HLL-based FD-WENO scheme with 200 zones.*

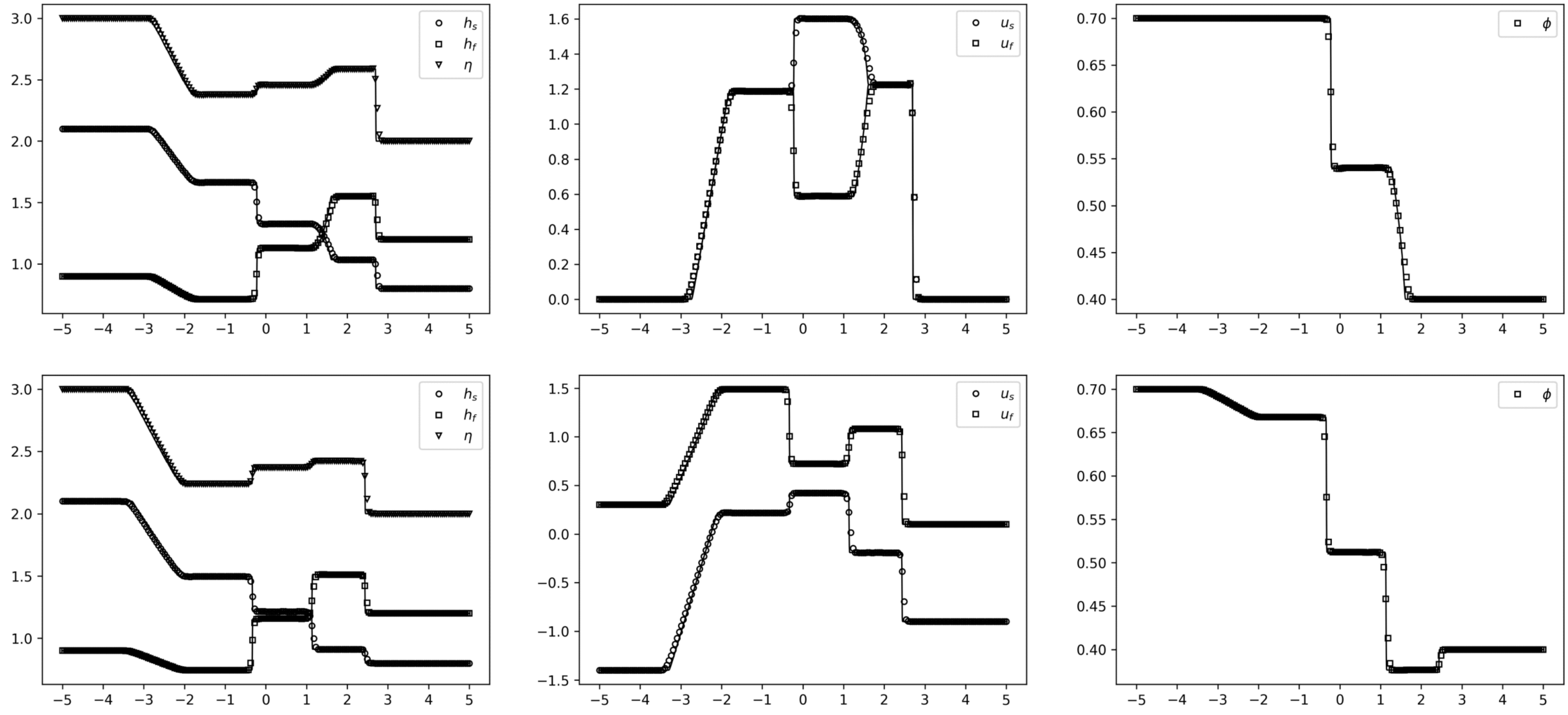

*Fig. 3b) Debris Flow: Riemann problem-2 (top) and -3 (bottom) using the 5th order accurate HLL-based FD-WENO scheme with 200 zones.*

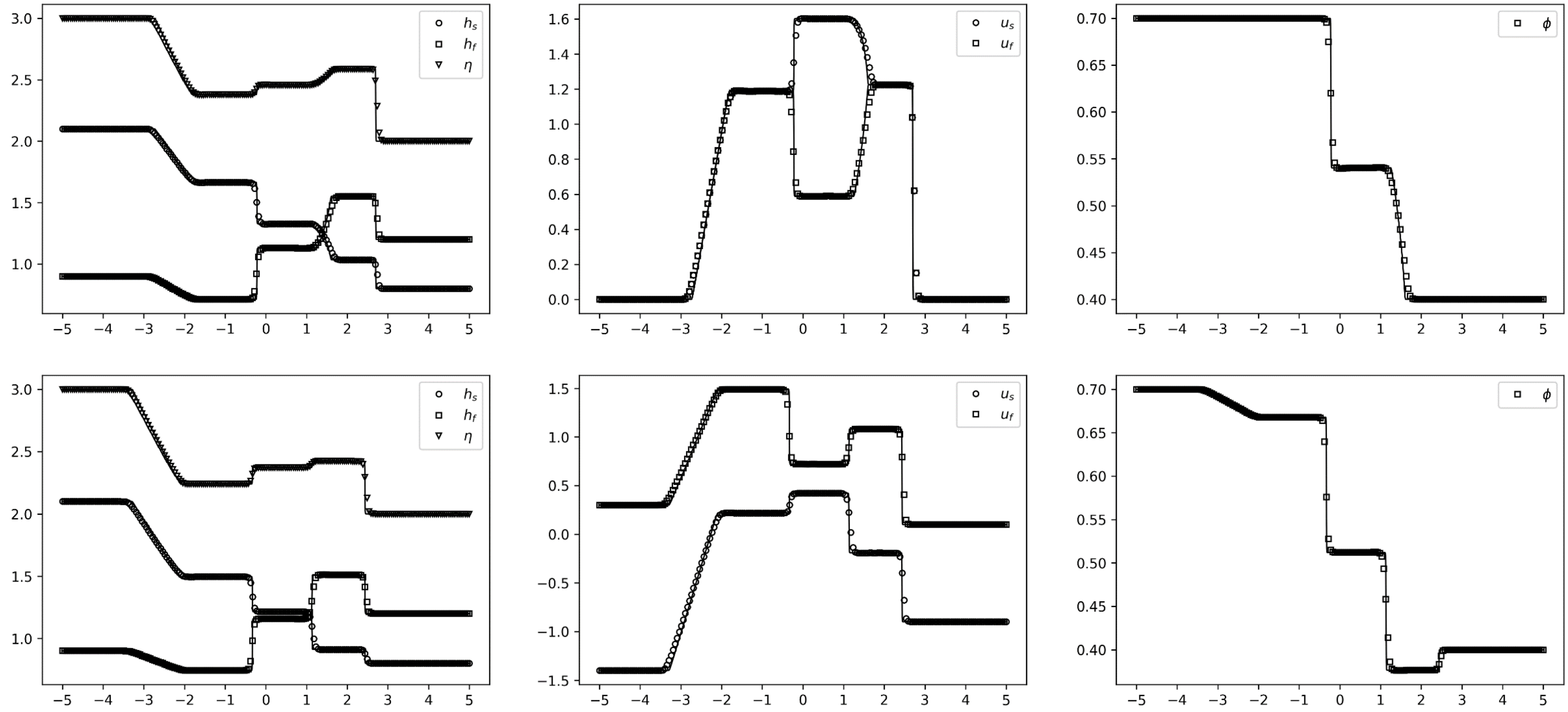

*Fig. 3c) Debris Flow: Riemann problem-2 (top) and -3 (bottom) using the 7$^{th}$ order accurate HLL-based FD-WENO scheme with 200 zones.*

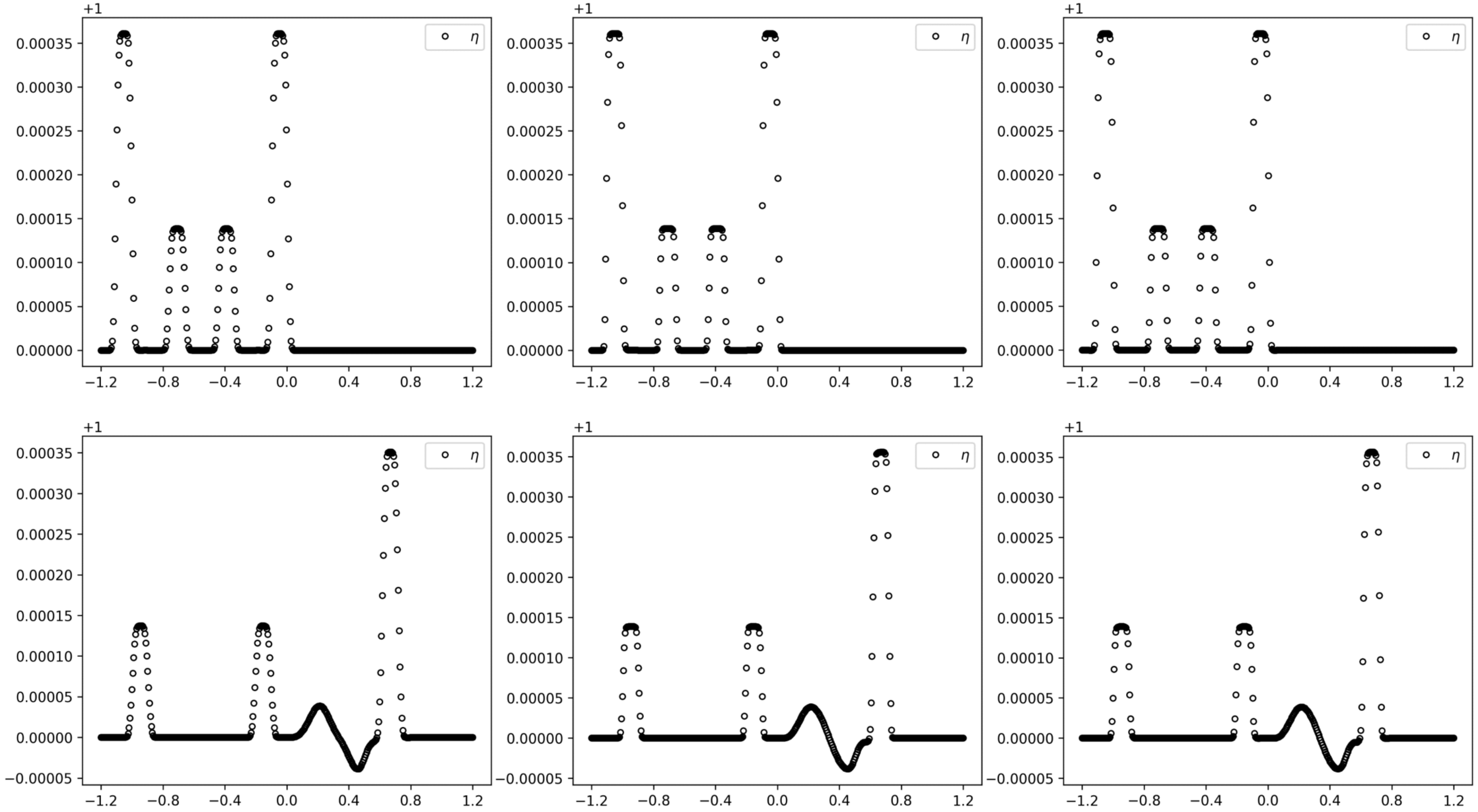

*Fig. 4) Debris Flow: Non-constant Bathymetry test using the 3rd order (left column), 5th order (middle column) and 7th order (right column) accurate HLL-based FD-WENO scheme with 400 zones. The top three panels show the solution at time of t=0.5; the lower three panels show the solution at a time of t=1.25.*

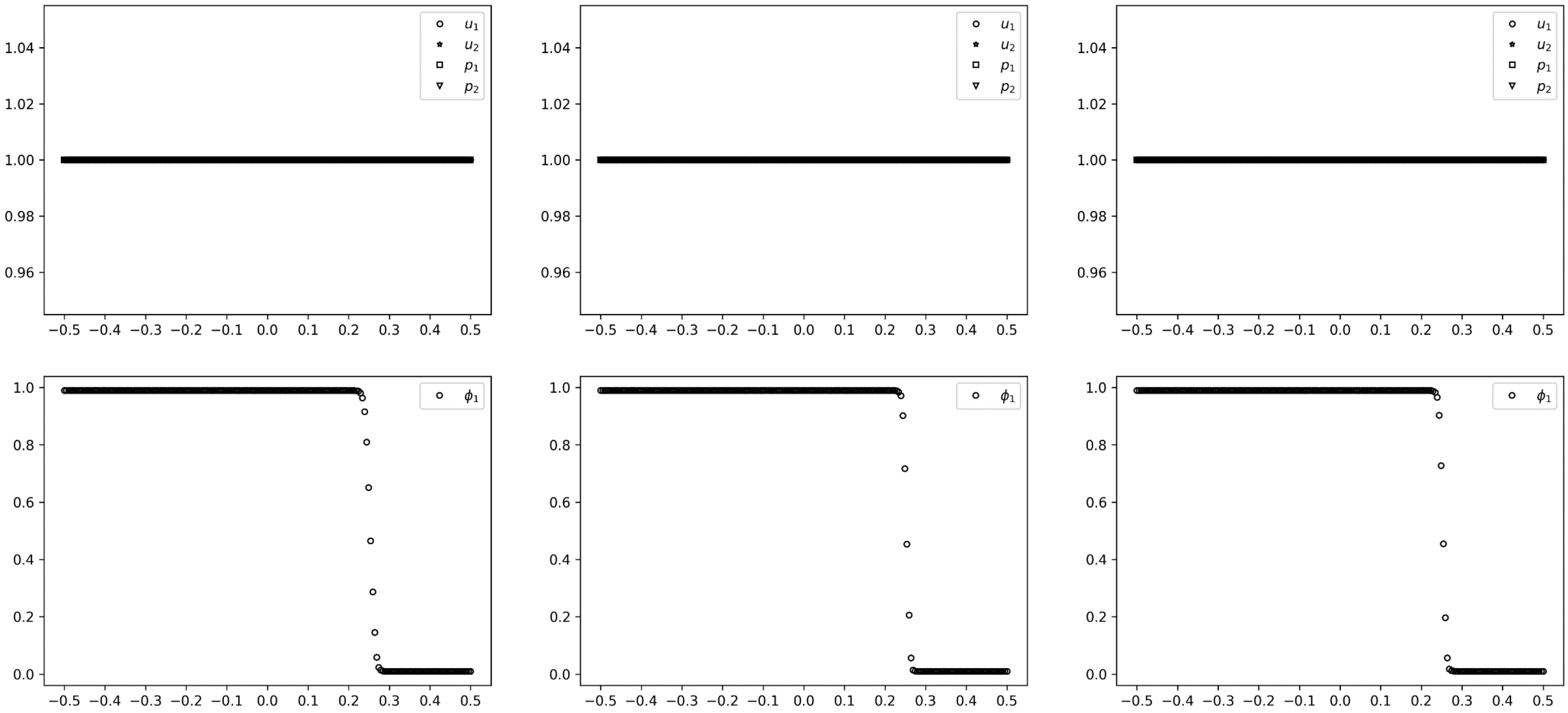

*Fig. 5a) Baer-Nunziato: Abgrall problem using the 3rd order (left column), 5th order (middle column) and 7th order (right column) accurate HLL-based FD-WENO scheme with 200 zones.*

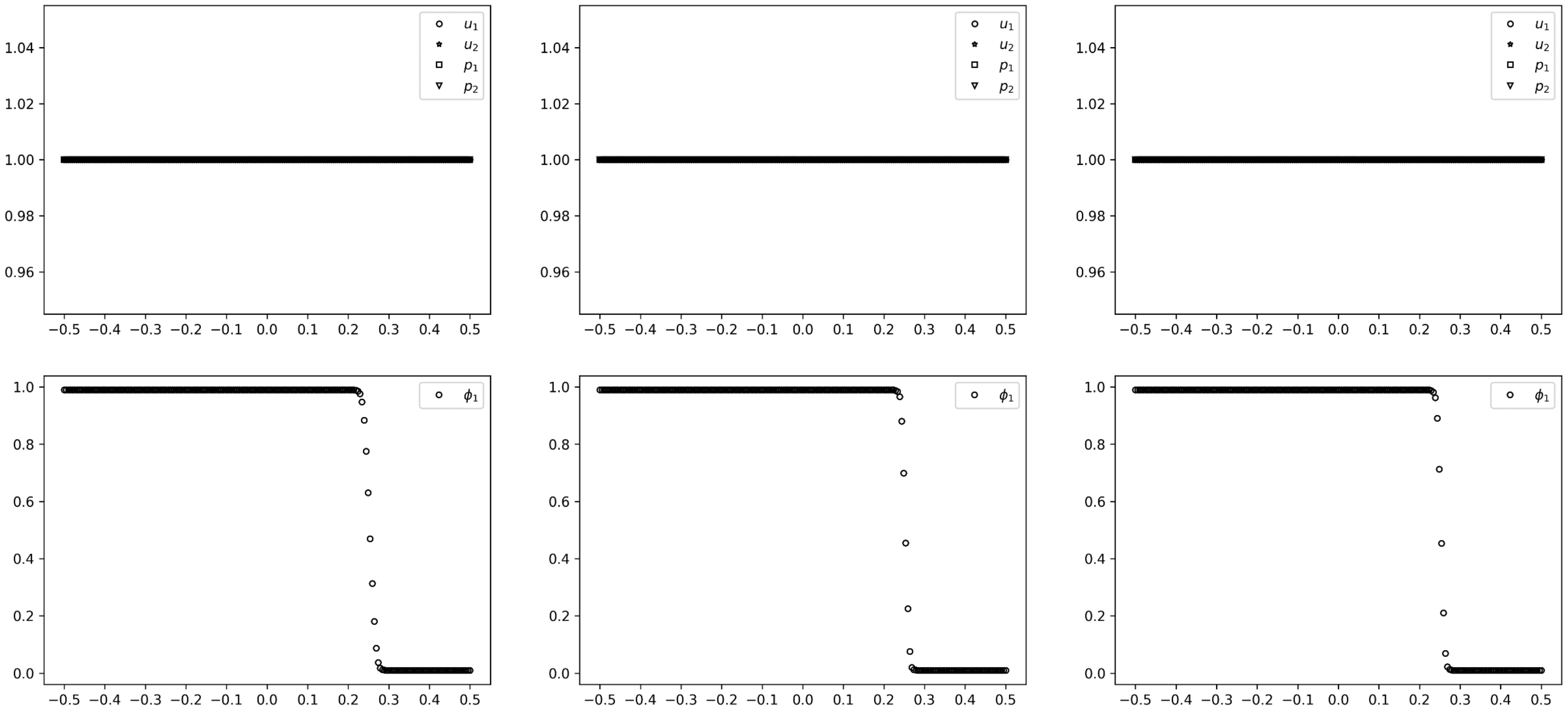

*Fig. 5b) Baer-Nunziato: Abgrall problem using the 3rd , 5th  and 7th order accurate LLF-based FD-WENO scheme with 200 zones.*

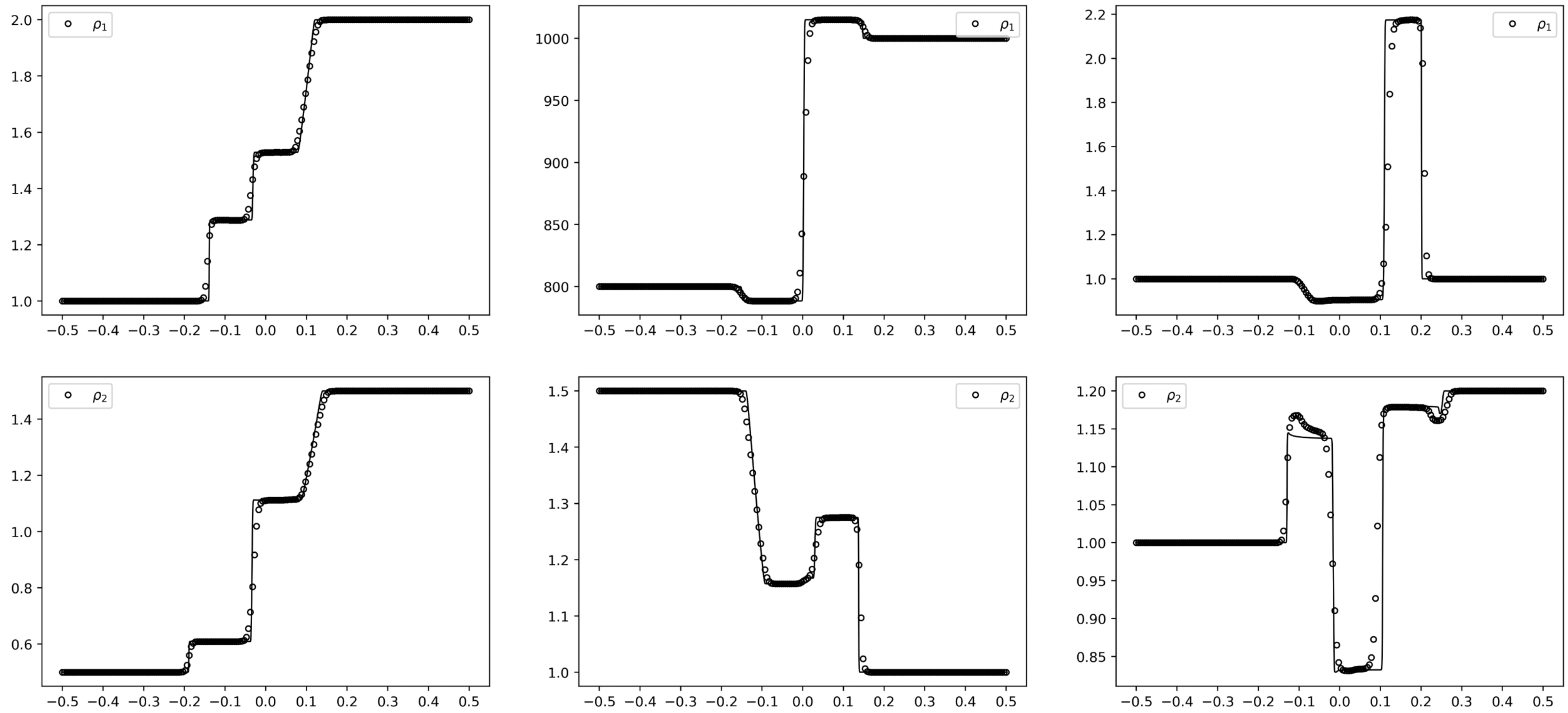

*Fig. 6a) Baer-Nunziato: RP-1 (left), RP-2 (middle) and RP-3 (right) using the 3rd order accurate HLL-based FD-WENO scheme with 200 zones.*

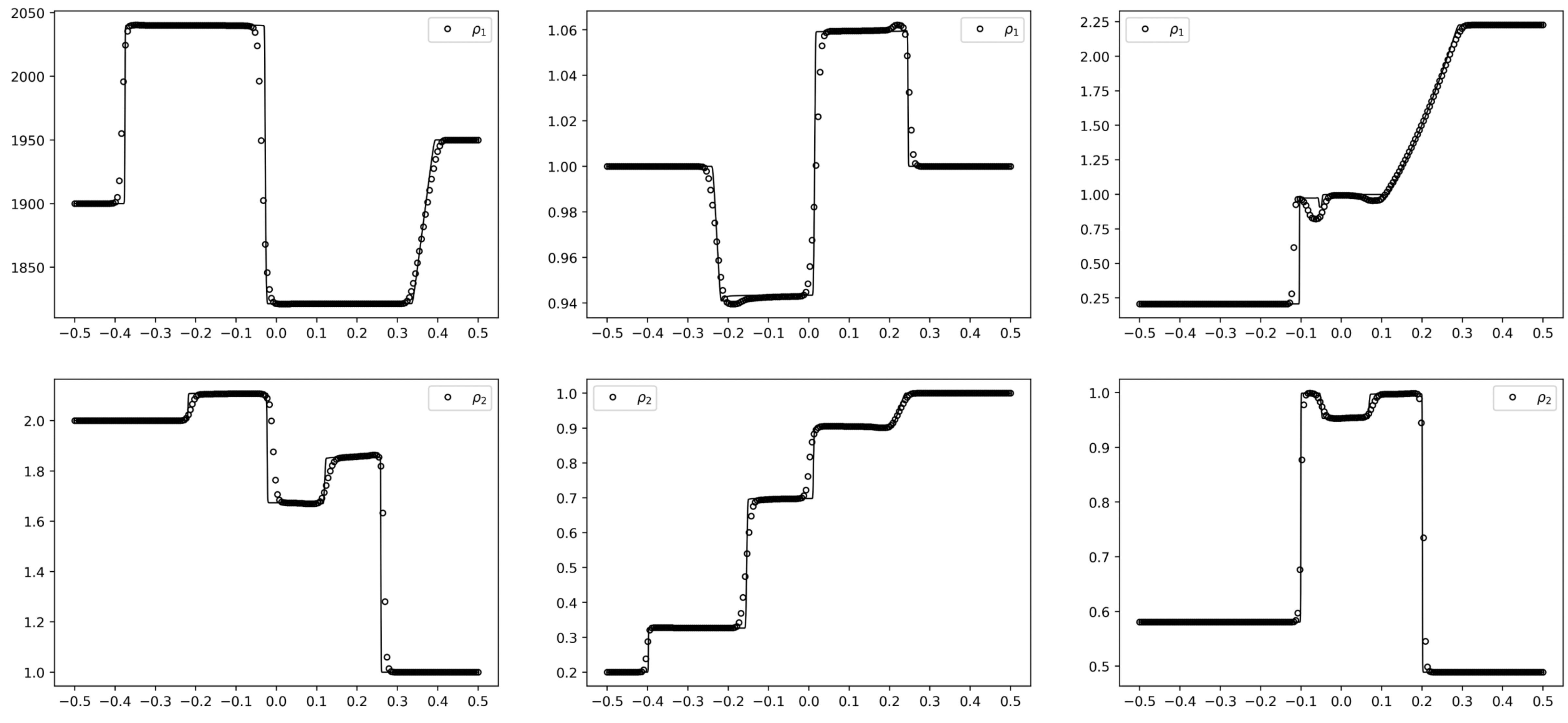

*Fig. 6b) Baer-Nunziato: RP-4 (left), RP-5 (middle) and RP-6 (right) using the 3rd order accurate HLL-based FD-WENO scheme with 200 zones.*

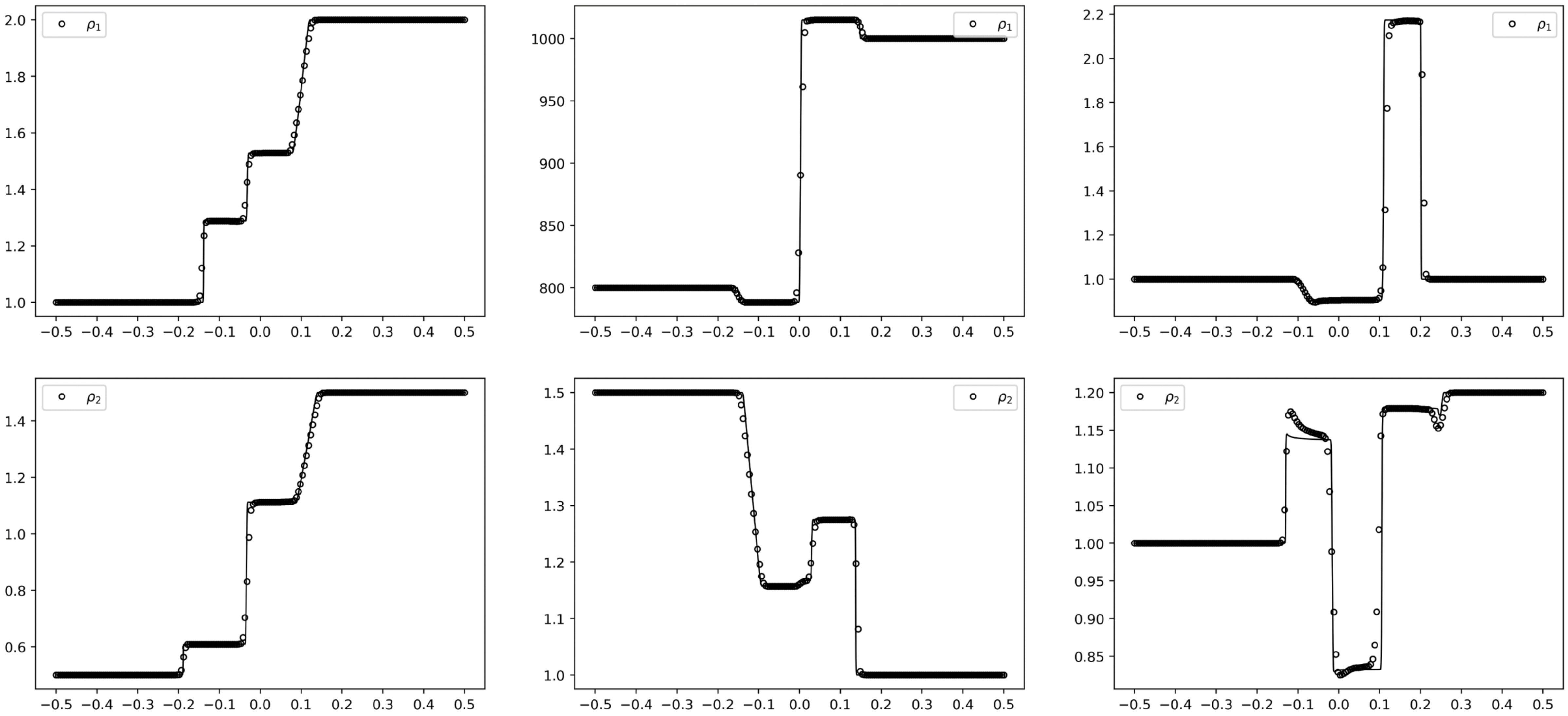

*Fig. 6c) Baer-Nunziato: RP-1 (left), RP-2 (middle) and RP-3 (right) using the $5^{th}$ order accurate HLL-based FD-WENO scheme with 200 zones.*

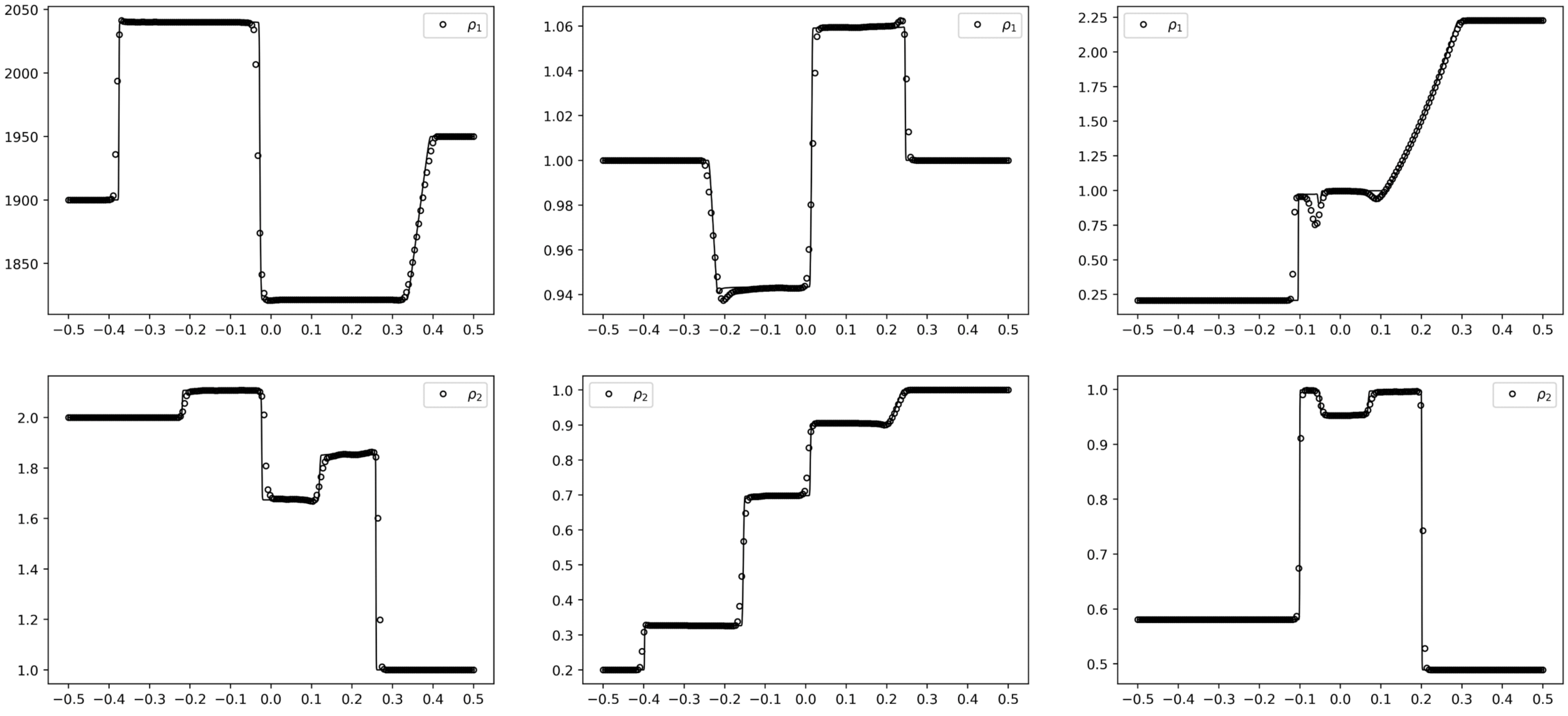

*Fig. 6d) Baer-Nunziato: RP-4 (left), RP-5 (middle) and RP-6 (right) using the 5th order accurate HLL-based FD-WENO scheme with 200 zones.*

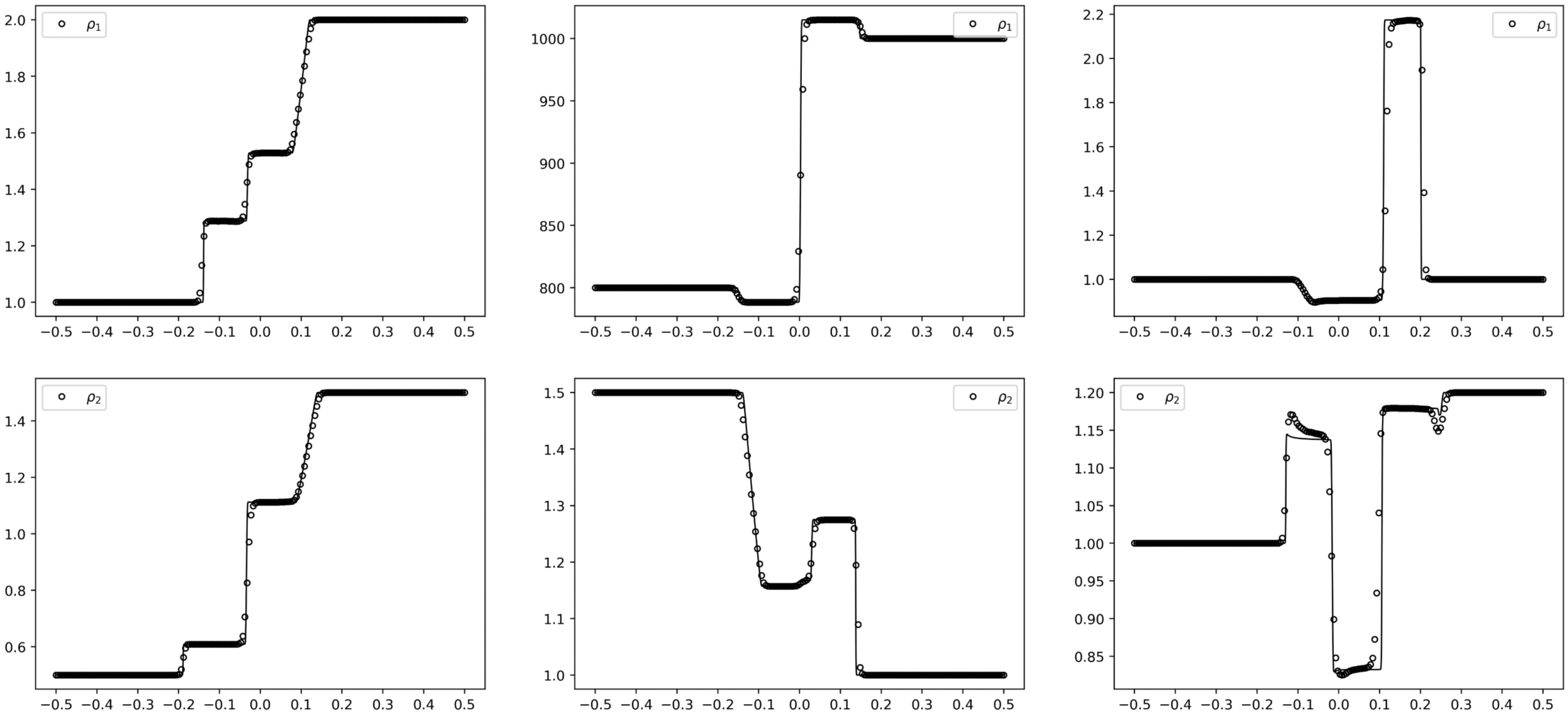

*Fig. 6e) Baer-Nunziato: RP-1 (left), RP-2 (middle) and RP-3 (right) using the 7th order accurate HLL-based FD-WENO scheme with 200 zones.*

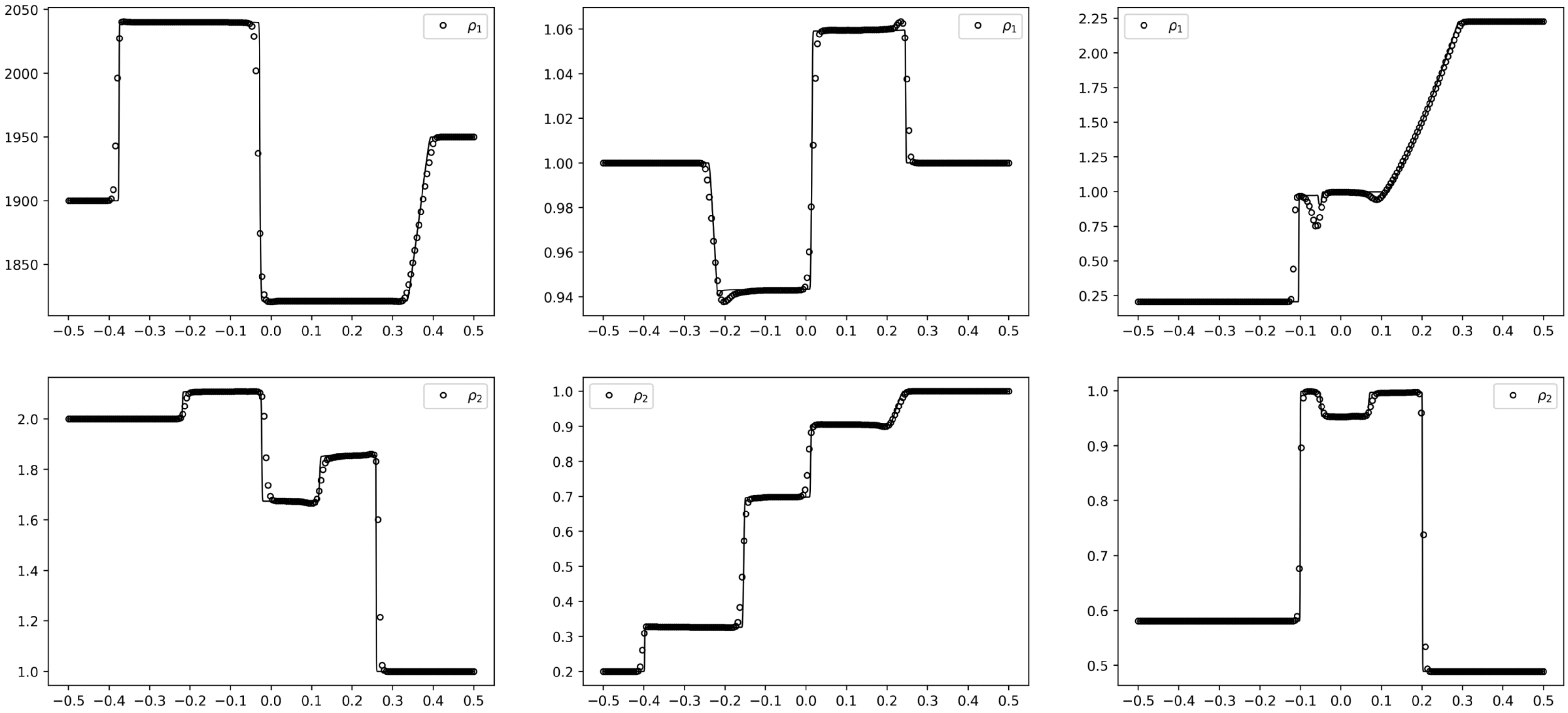

*Fig. 6f) Baer-Nunziato: RP-4 (left), RP-5 (middle) and RP-6 (right) using the 7th order accurate HLL-based FD-WENO scheme with 200 zones*

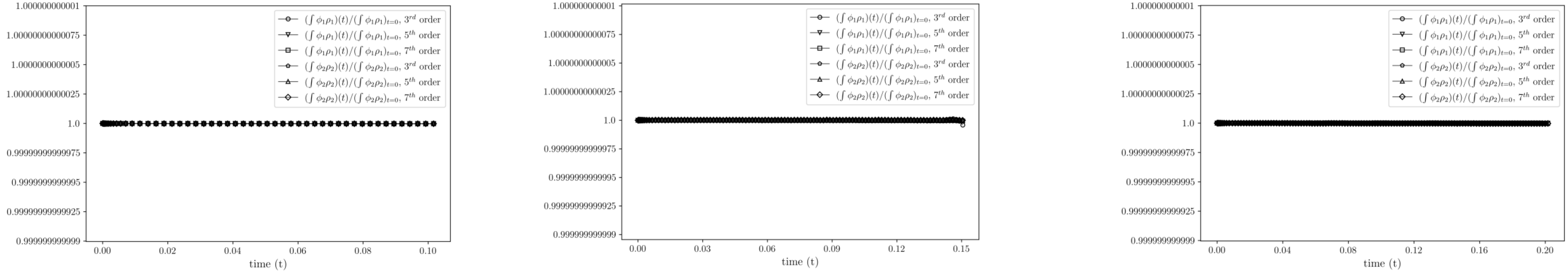

*Fig. 6g) Verification of mass conservation in the solid and gas masses for the Riemann problem-2 (left panel), Riemann problem-4 (middle panel) and Riemann problem-5 (right panel) using the $3^{rd}$, $5^{th}$ and $7^{th}$ order HLL-based FD-WENO schemes. Conservation of solid and gas masses is respected to at least one part in $10^{12}$ .*

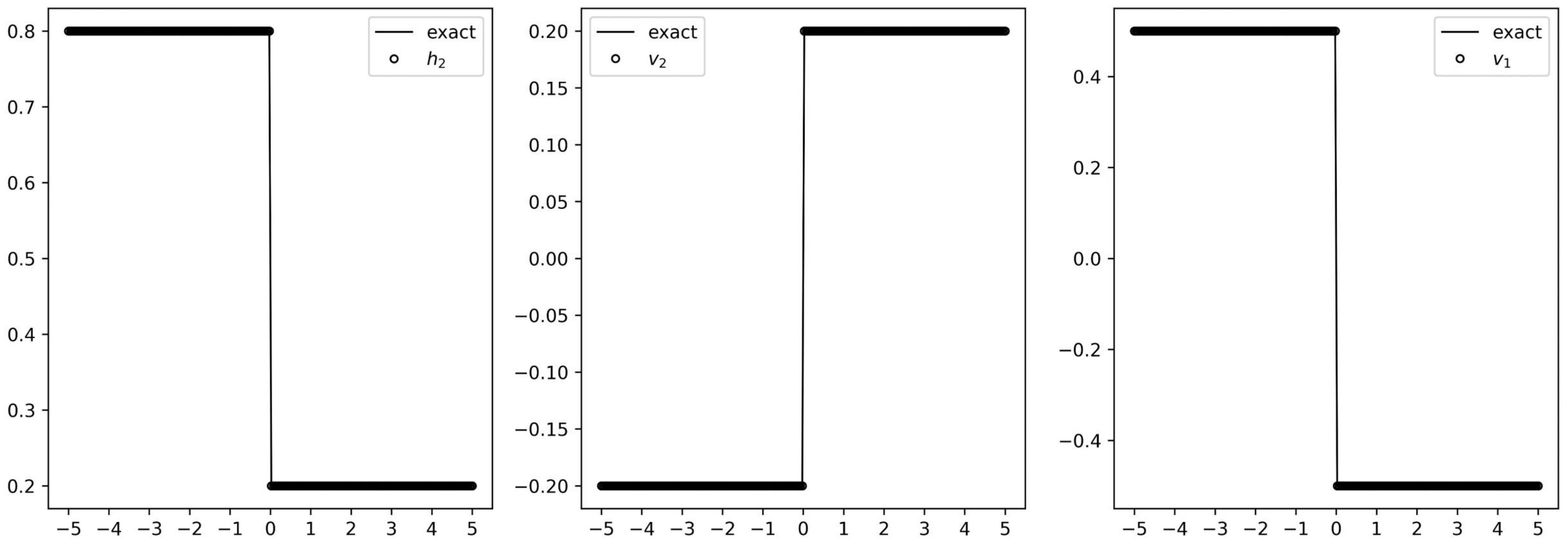

*Fig. 7) Two-Layer Shallow water: Riemann problem-1 using the 3rd order accurate HLL-based FD-WENO scheme with 200 zones. The 3rd order LLF-based FD-WENO also shows identical results. Identical results were also obtained for 5th and 7th order HLL-based FD-WENO and LLF-based FD-WENO schemes.*

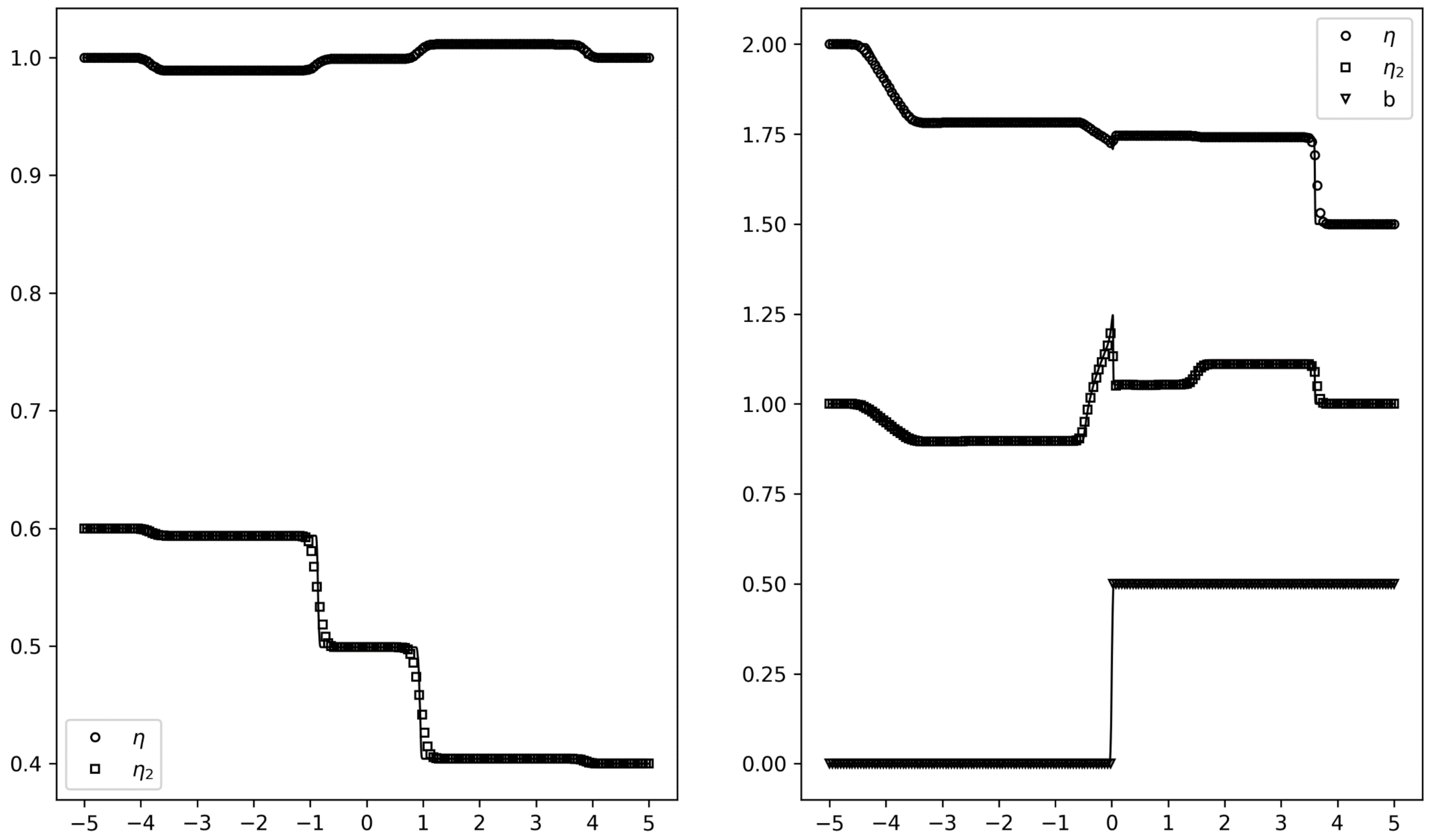

*Fig. 8a) Two-Layer Shallow water: Riemann problem-2(left) and -3(right) using the 3rd order accurate HLL-based FD-WENO scheme with 200 zones.*

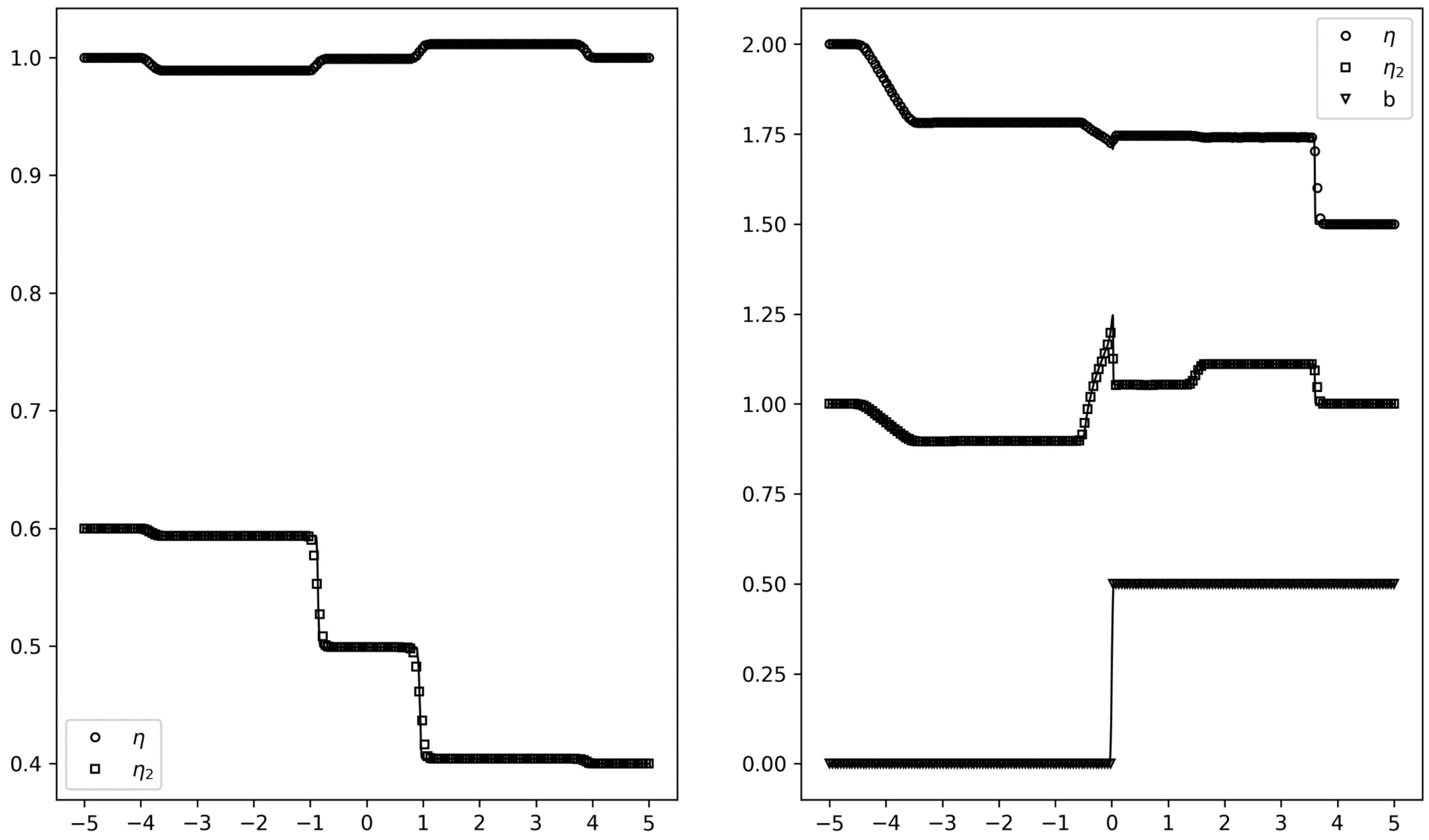

*Fig. 8b) Two-Layer Shallow water: Riemann problem-2(left) and -3(right) using the 5th order accurate HLL-based FD-WENO scheme with 200 zones.*

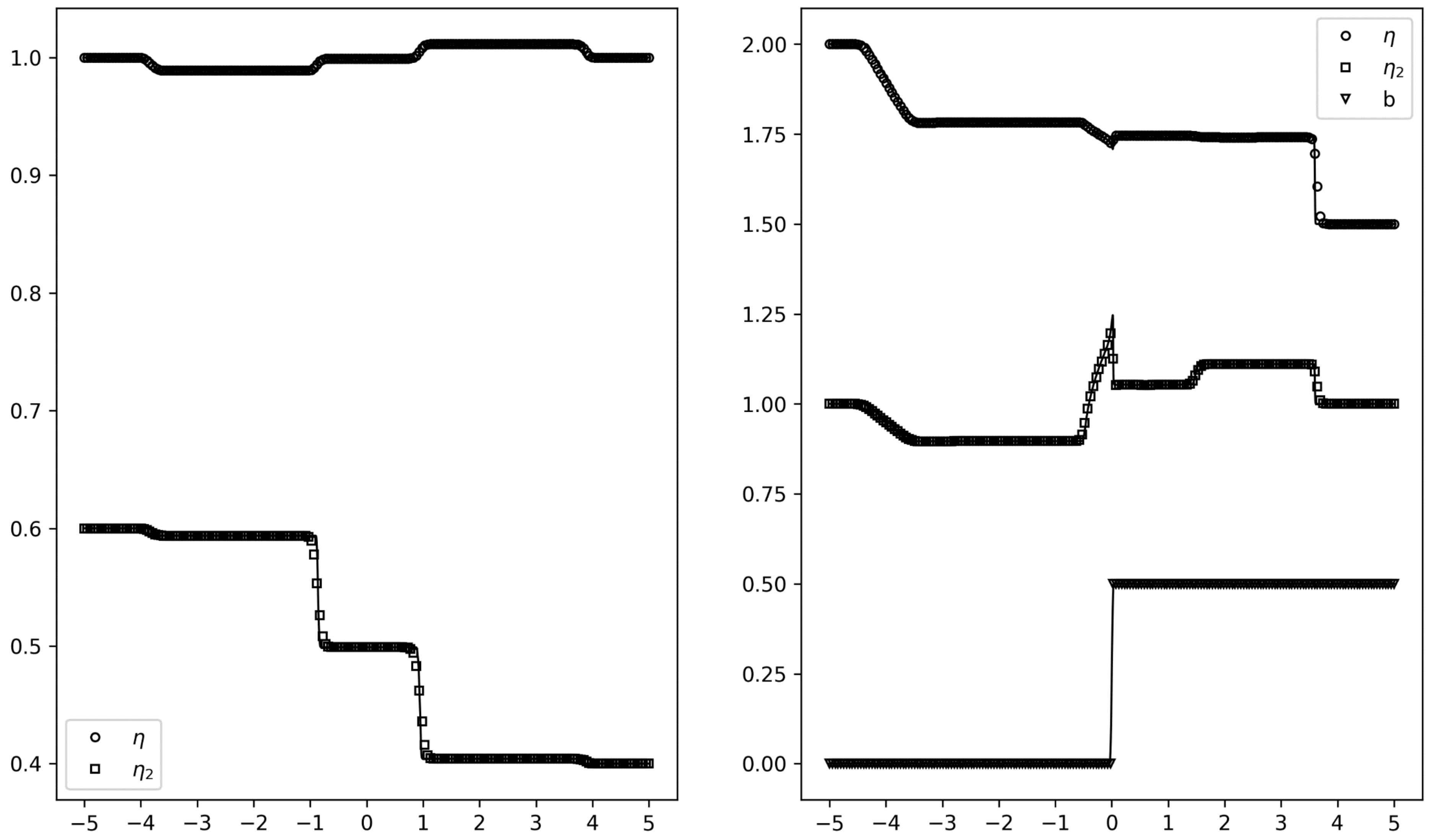

*Fig. 8c) Two-Layer Shallow water: Riemann problem-2(left) and -3(right) using the 7th order accurate HLL-based FD-WENO scheme with 200 zones.*

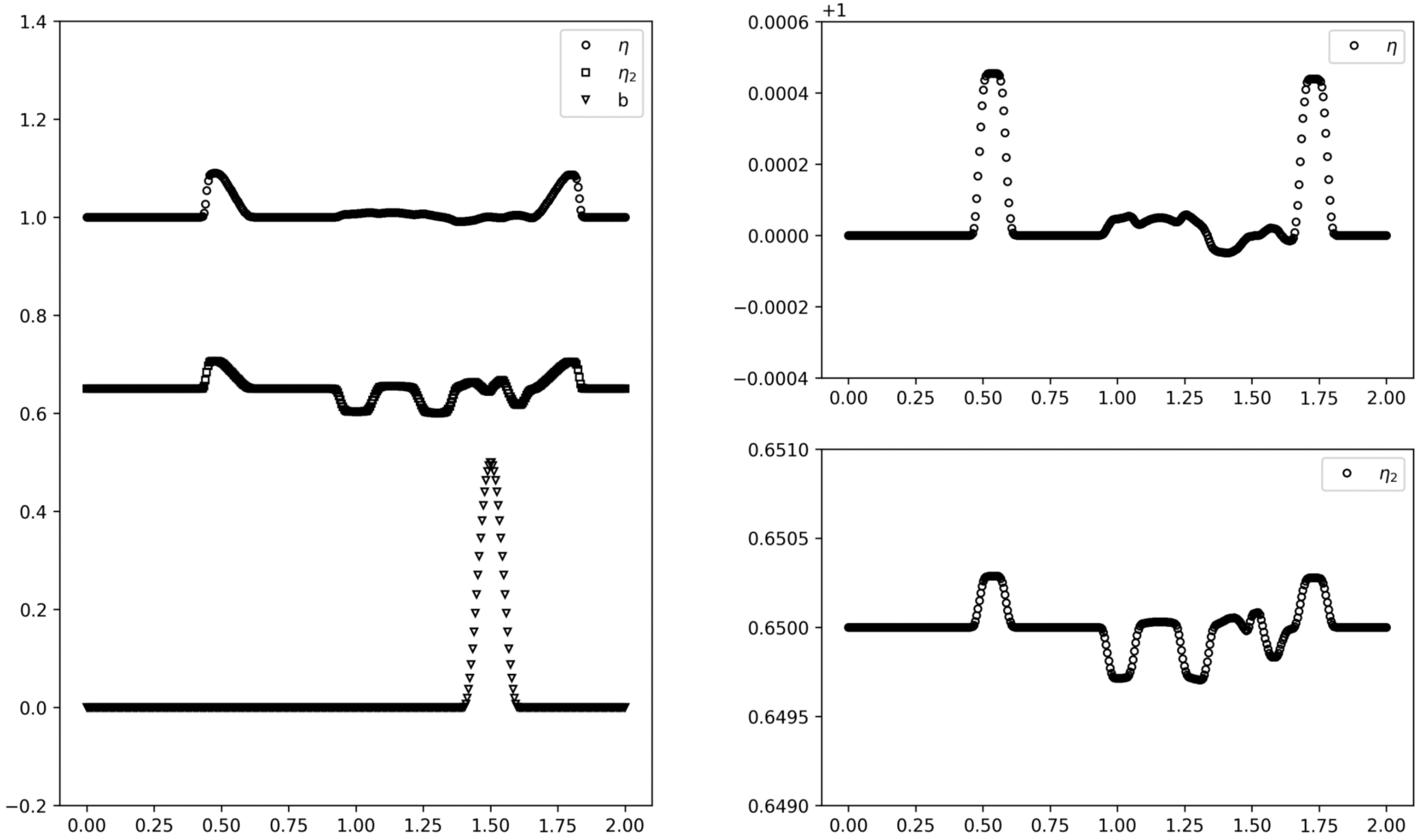

*Fig. 9a) Two-Layer Shallow Water Model: Non-constant Bathymetry test. Left: Large perturbation, Right: Small perturbation using the 3rd order accurate HLL-based FD-WENO scheme with 400 zones. The 3rd order LLF-based FD-WENO scheme produces identical results.*

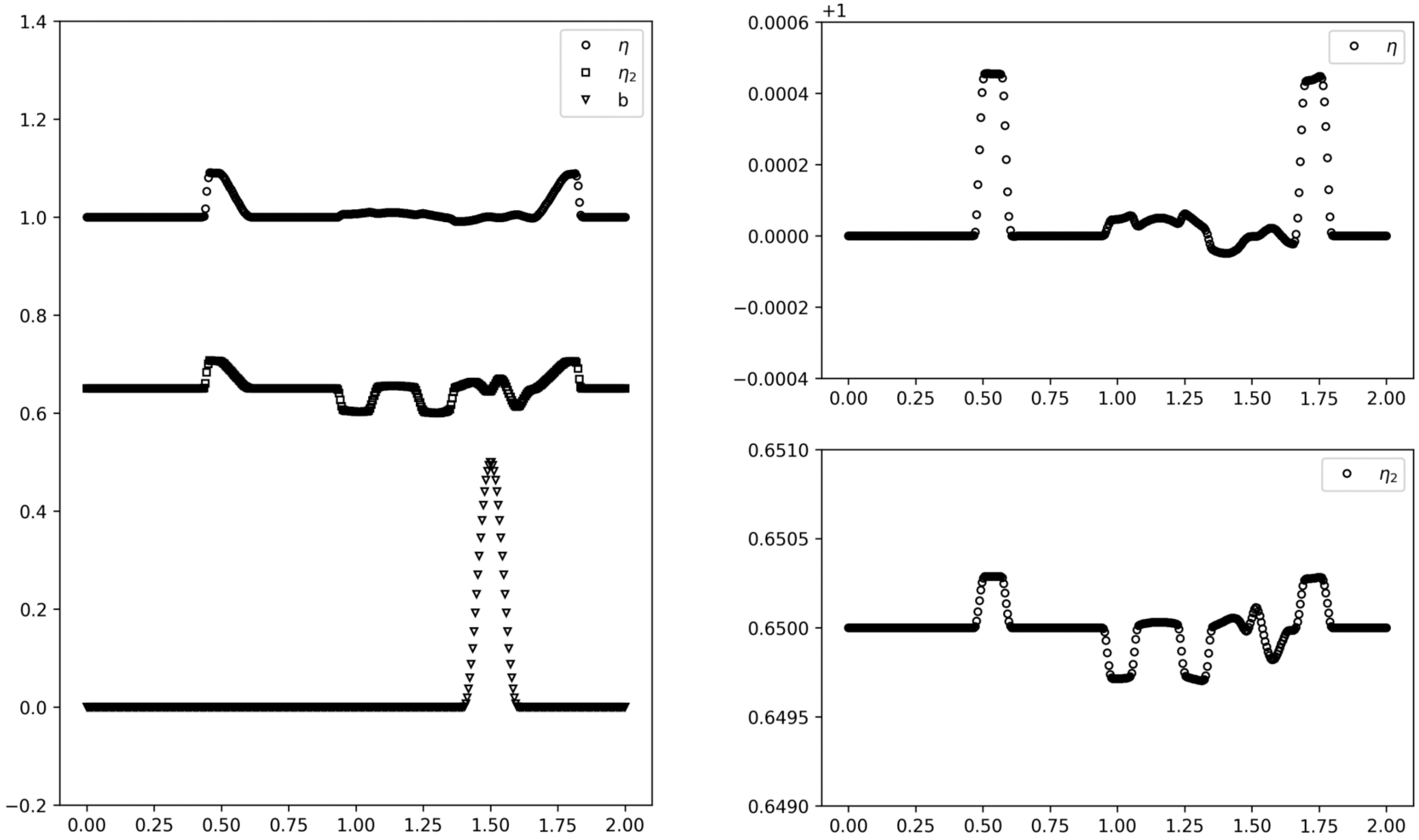

*Fig. 9b) Two-Layer Shallow Water Model: Non-constant Bathymetry test. Left: Large perturbation, Right: Small perturbation using the 5$^{th}$ order accurate HLL-based FD-WENO scheme with 400 zones. The 5$^{th}$ order LLF-based FD-WENO scheme produces identical results.*

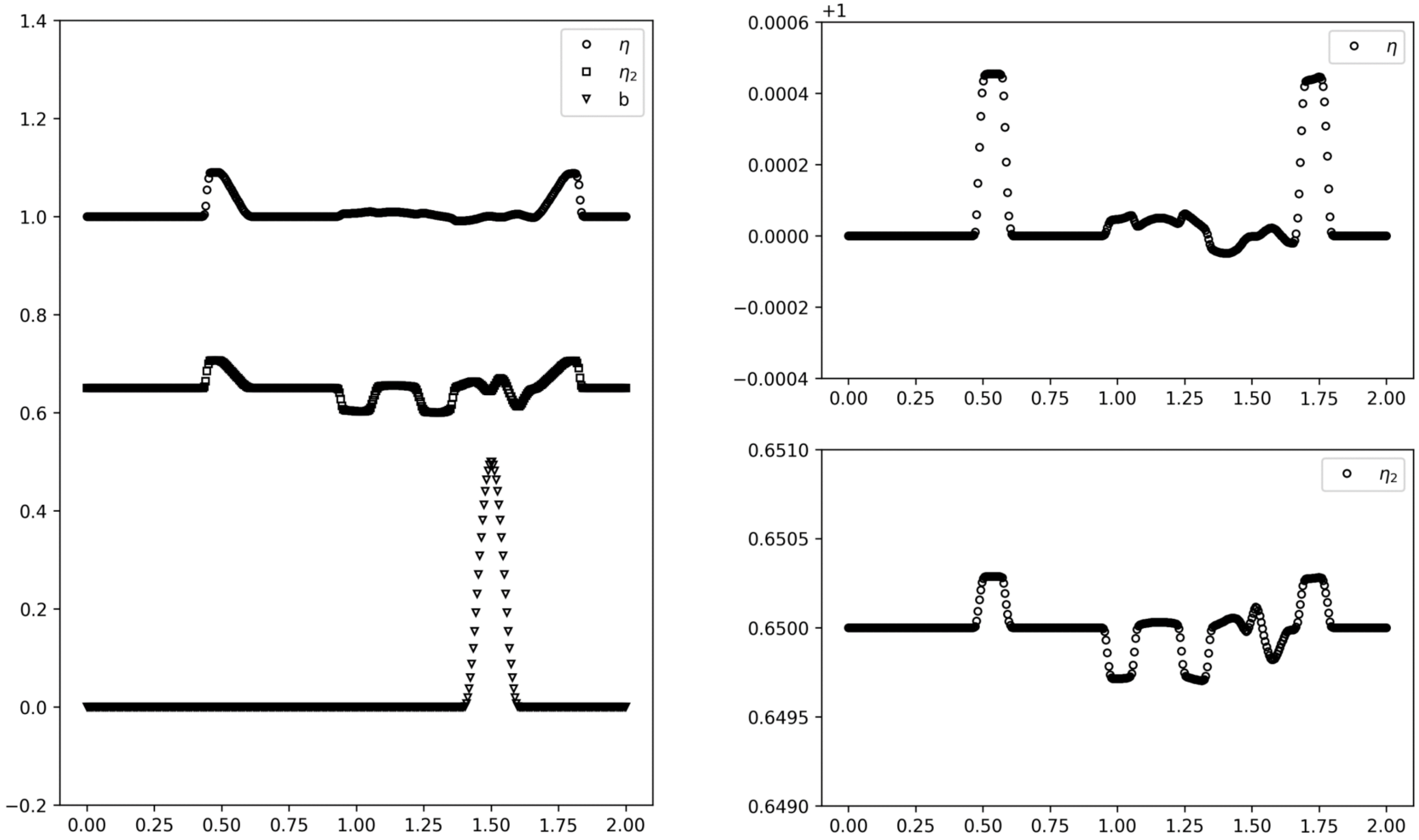

*Fig. 9c) Two-Layer Shallow Water Model: Non-constant Bathymetry test. Left: Large perturbation, Right: Small perturbation using the 7th order accurate HLL-based FD-WENO scheme with 400 zones. The 7th order LLF-based FD-WENO scheme produces identical results.*

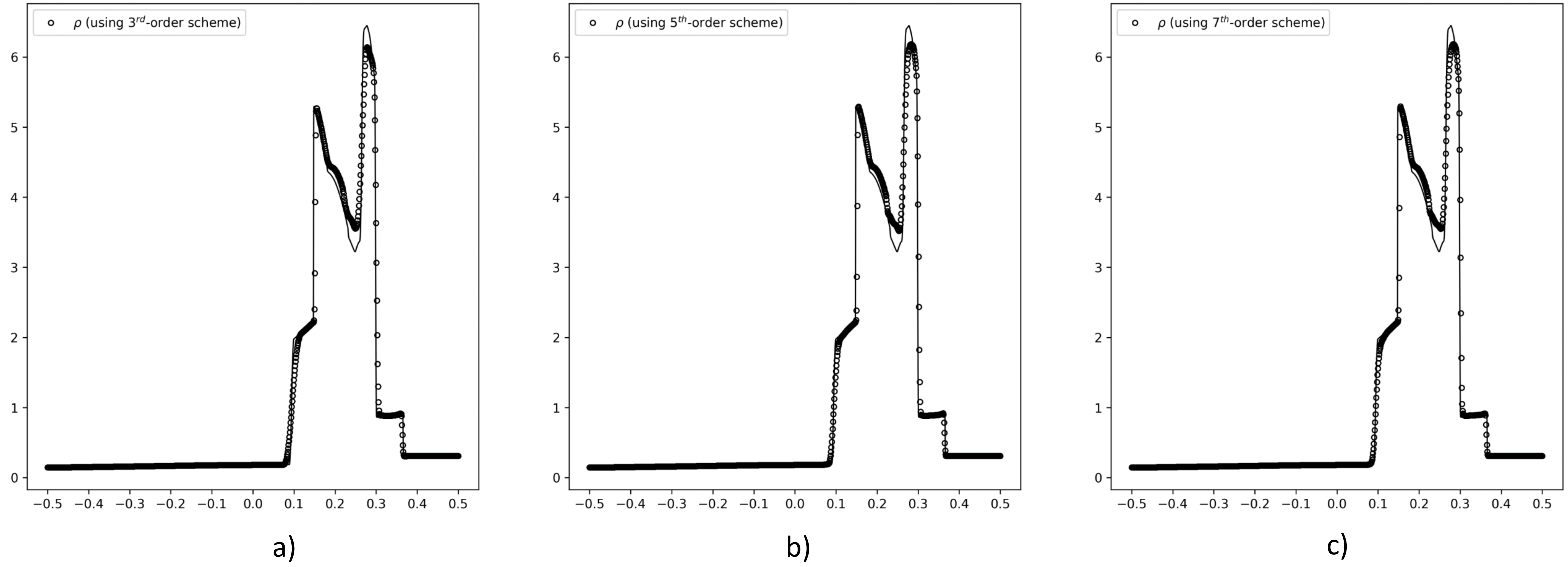

*Fig. 10) EULER: Blast wave interaction using the 3rd , 5th and 7th order accurate HLL-based FD-WENO scheme with 1000 zones. Figs. 10a, 10b and 10c individually show the density at 3rd , 5th and 7th orders.*

a)

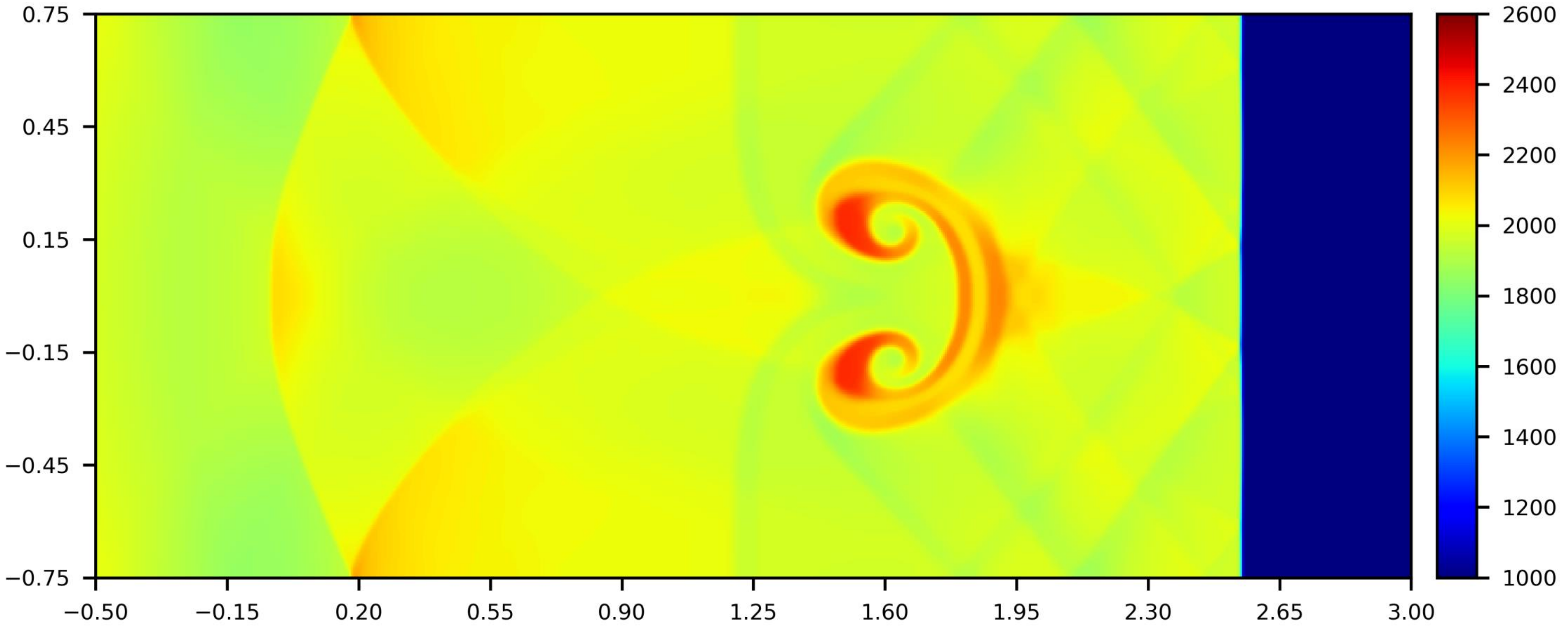

b)

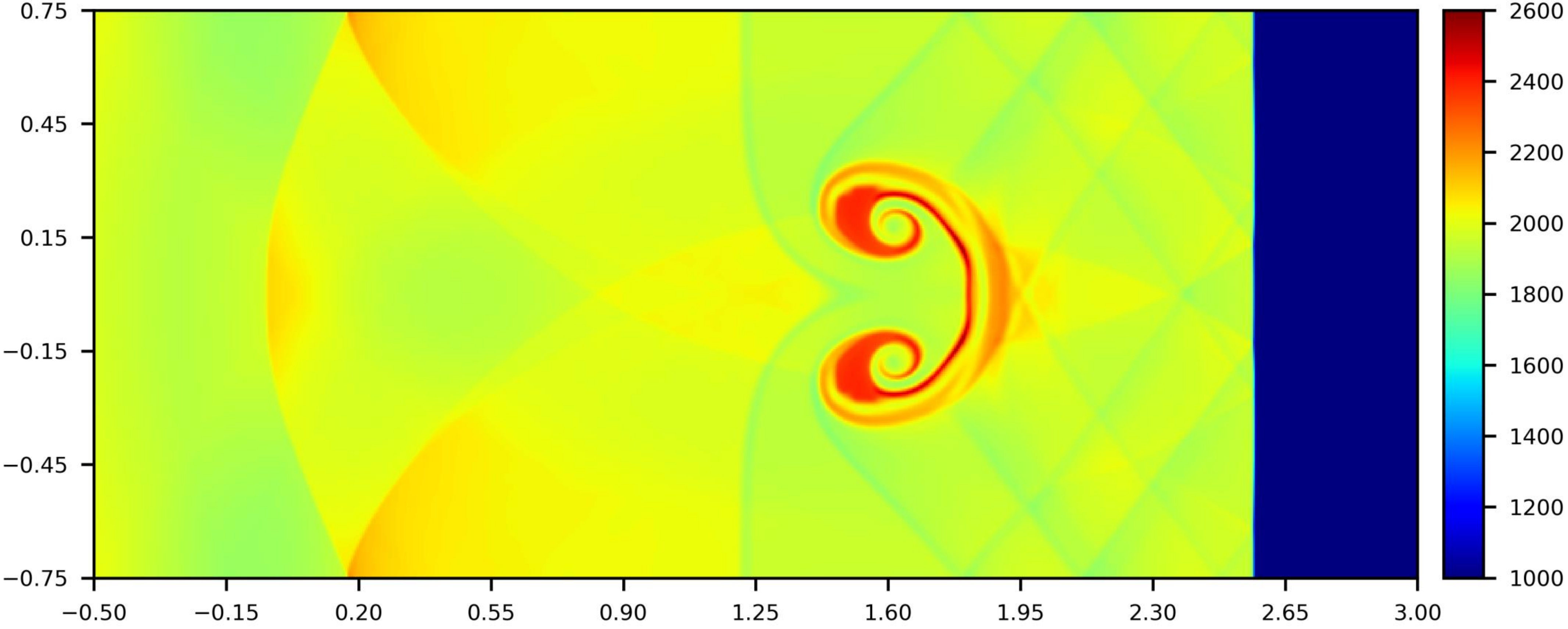

c)

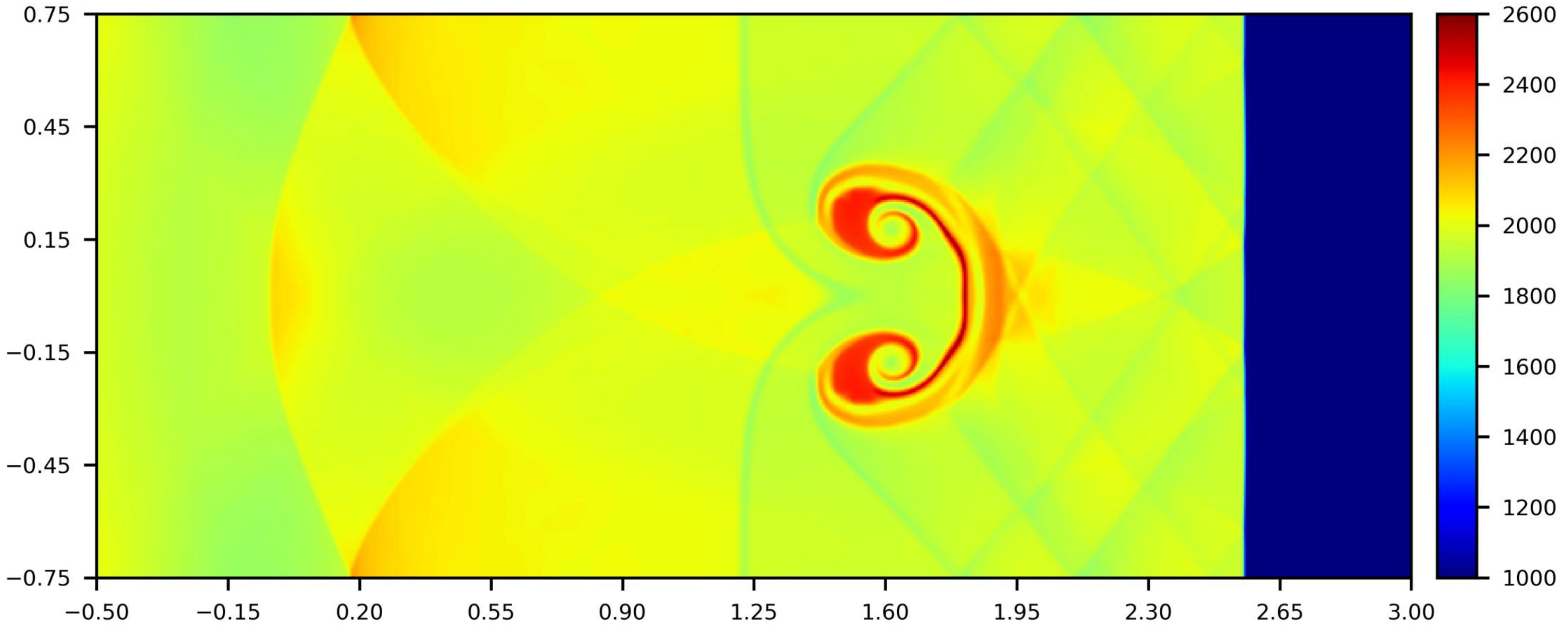

*Fig. 11a, b, c) Baer Nunziato: Shock-bubble interaction problem using the 3rd, 5th and 7th order accurate HLL-based FD-WENO schemes with 700×300 zones. The solid density profiles has been shown.*

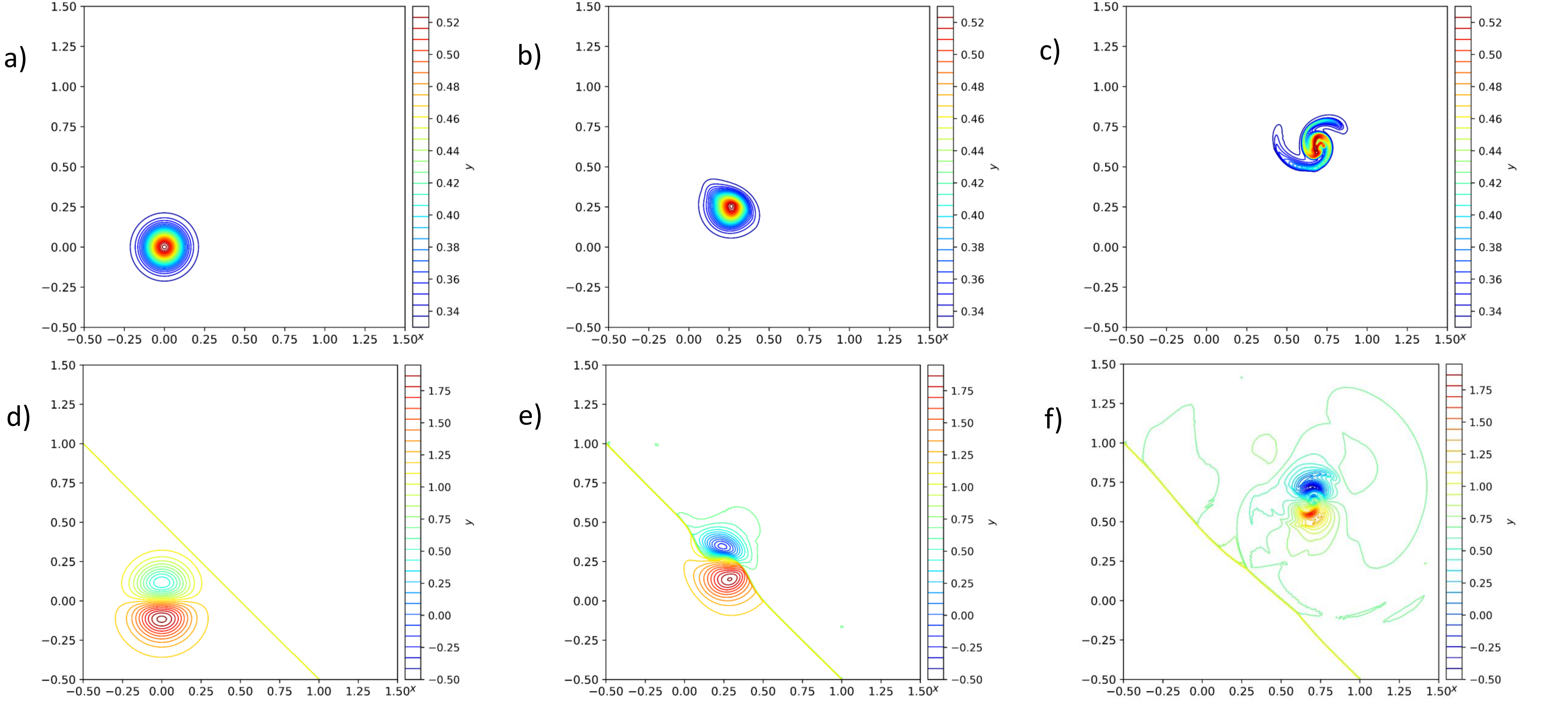

*Fig. 12) Baer-Nunziato: Shock-Vortex Interaction using the 7[th] order accurate HLL-based FD-WENO scheme with 600 ×600 zones at time levels t=0.0, 0.23 and 0.84. Figs. 12a, 12b and 12c show solid volume fraction at times t=0.0, 0.23, 0.84. Figs. 12d, 12e, 12f show solid x-velocity at times t=0.0, 0.23, 0.84. For the solid volume fraction, 30 contours were fit between a range of 0.33 and 0.530. For the solid x-velocity, 30 contours were fit between a range of -0.5 and 1.95.*

a)

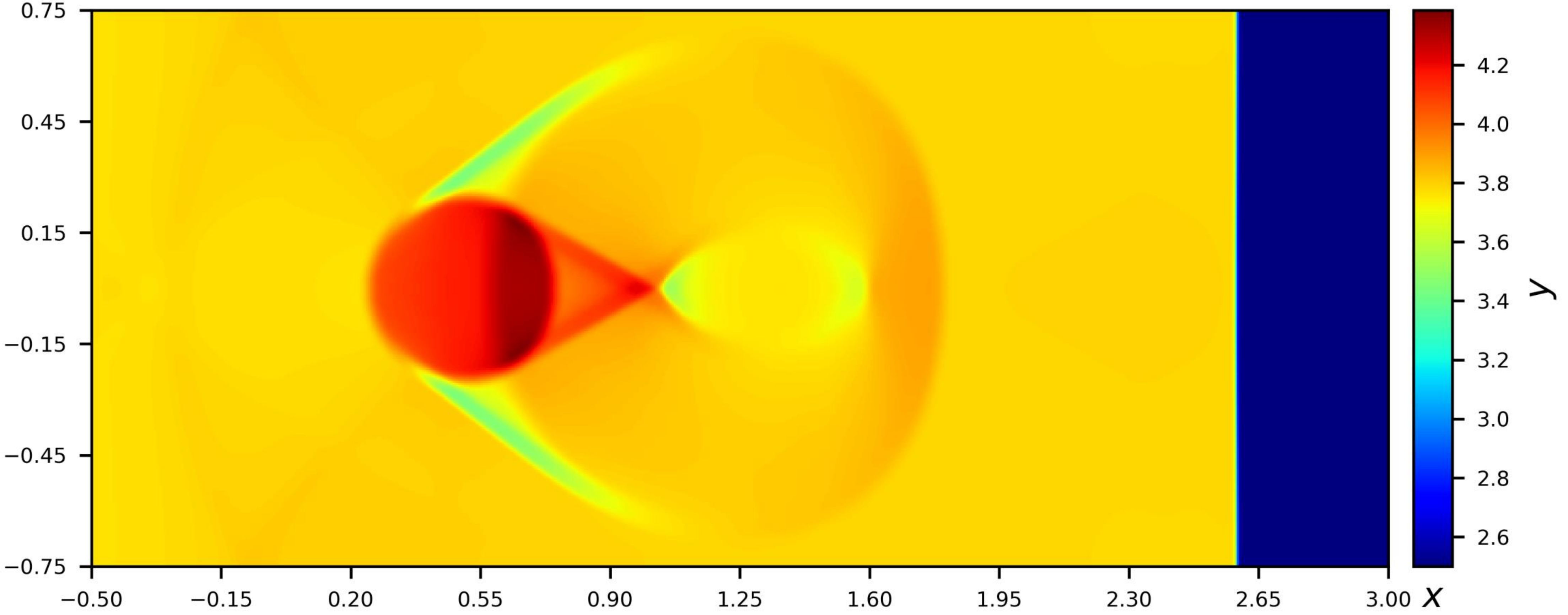

b)

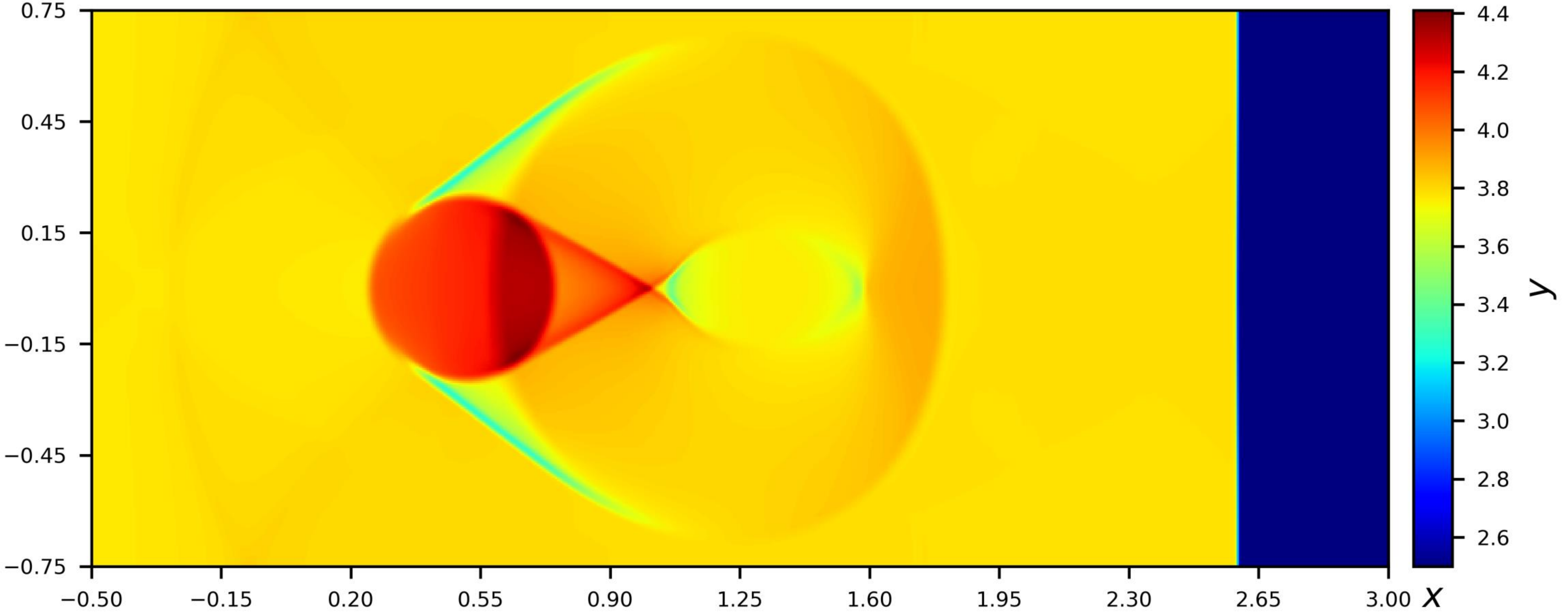

c)

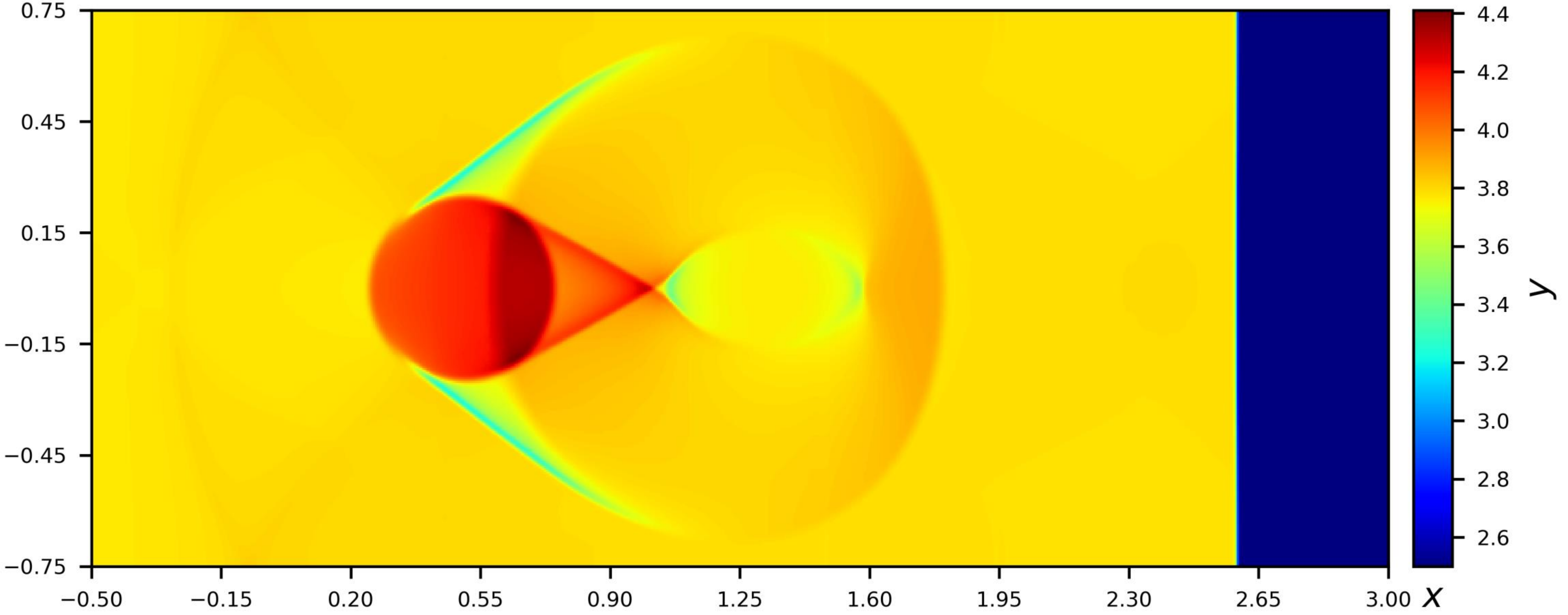

*Fig. 13a, b, c) Two-Layer Shallow Water Model : Shock-bubble interaction problem using the 3rd, 5th and 7th order accurate HLL-based FD-WENO schemes with 700×300 zones. Height of the lower fluid has been shown.*

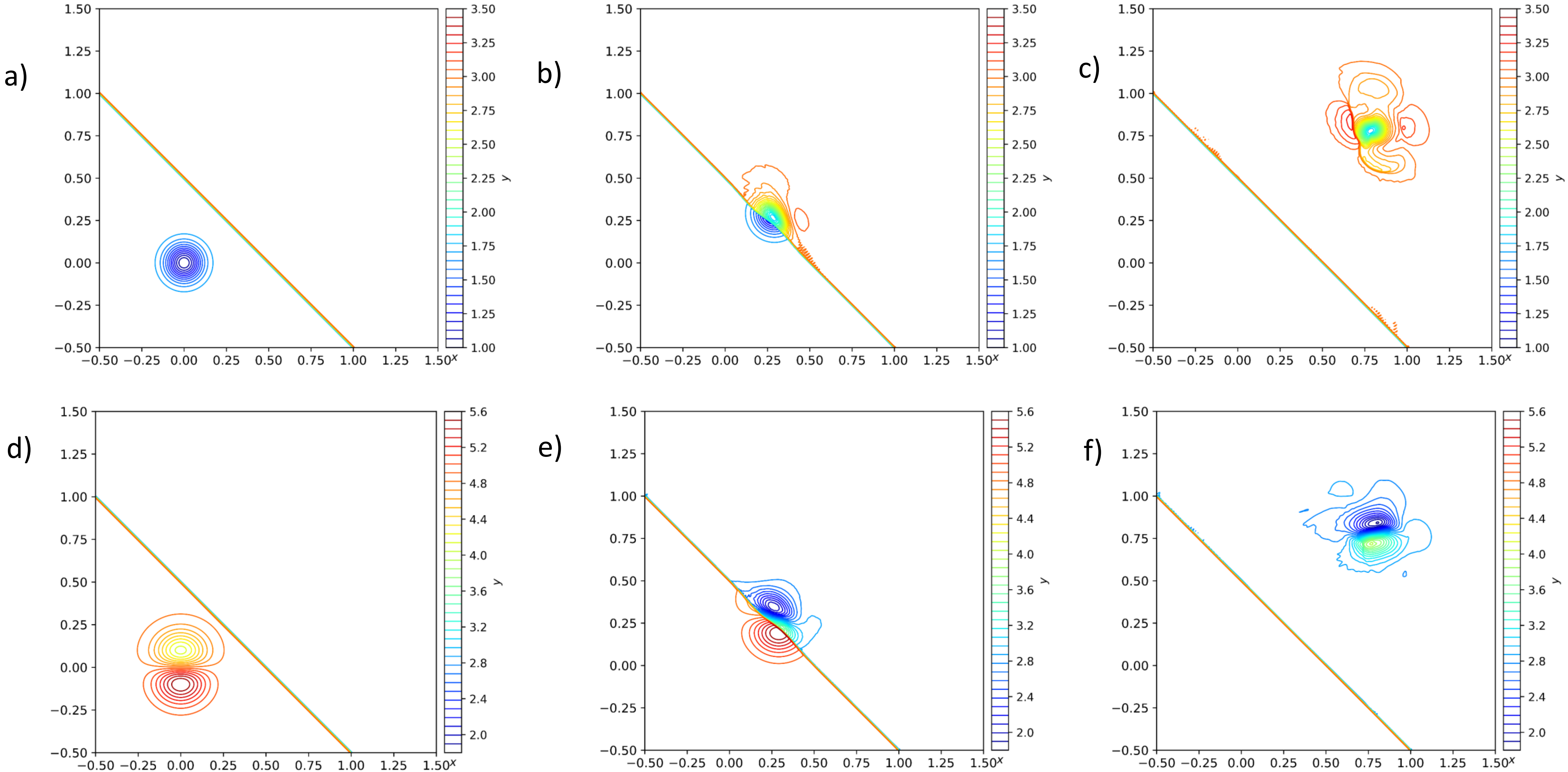

*Fig. 14) Two-Layer Shallow Water Model: Shock-Vortex Interaction using the 7th order accurate HLL-based FD-WENO scheme with 600 ×600 zones at time levels t=0.0, 0.06 and 0.24. Figs. 14a, 14b and 14c show height of the upper fluid at times t=0.0, 0.06, 0.24. Figs. 14d, 14e, 14f show x-velocity of the upper fluid at times t=0.0, 0.06, 0.24. For the height, 40 contours were fit between a range of 1.0 and 3.5. For the velocity, 40 contours were fit between a range of 1.8 and 5.6.*

a)

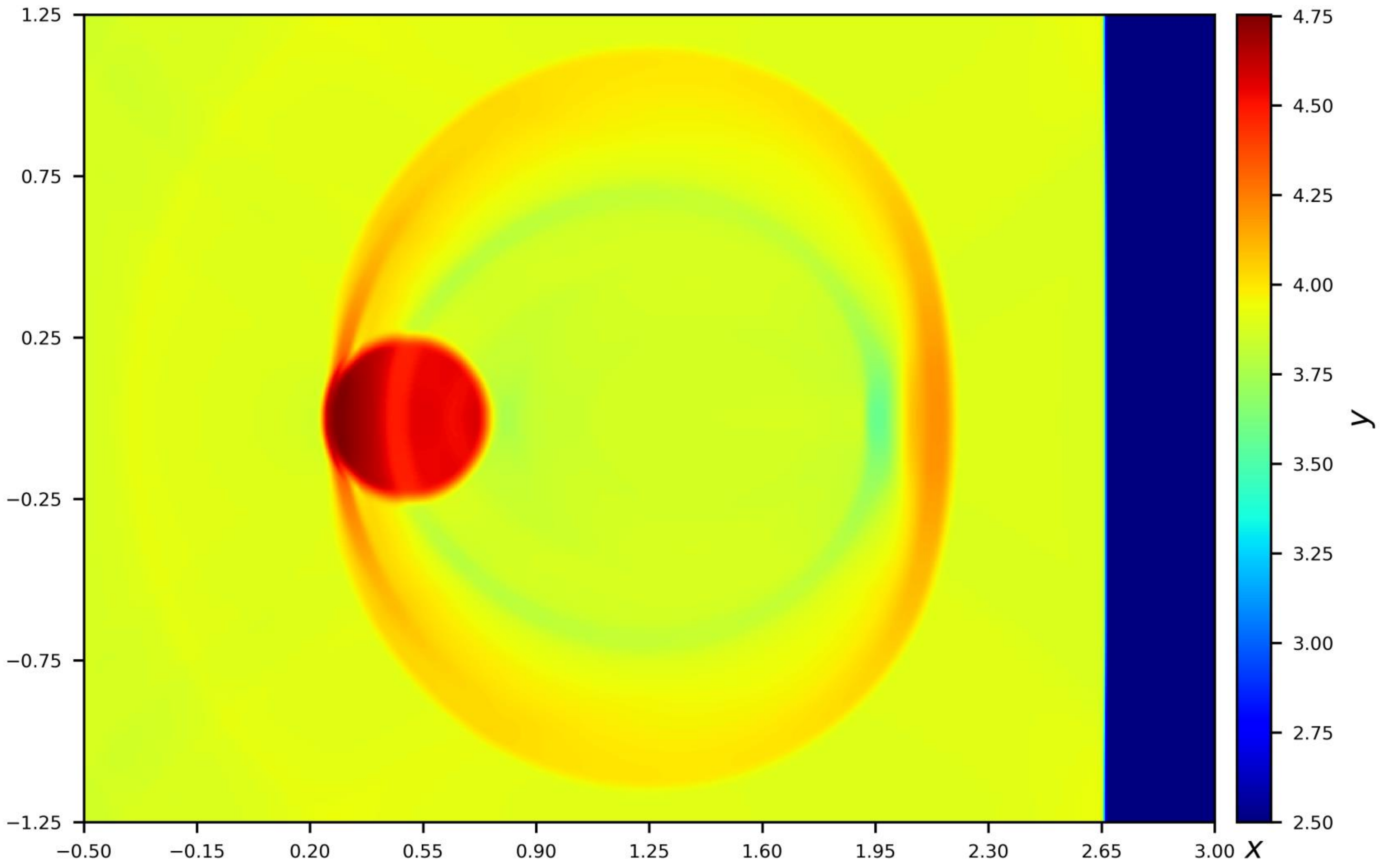

b)

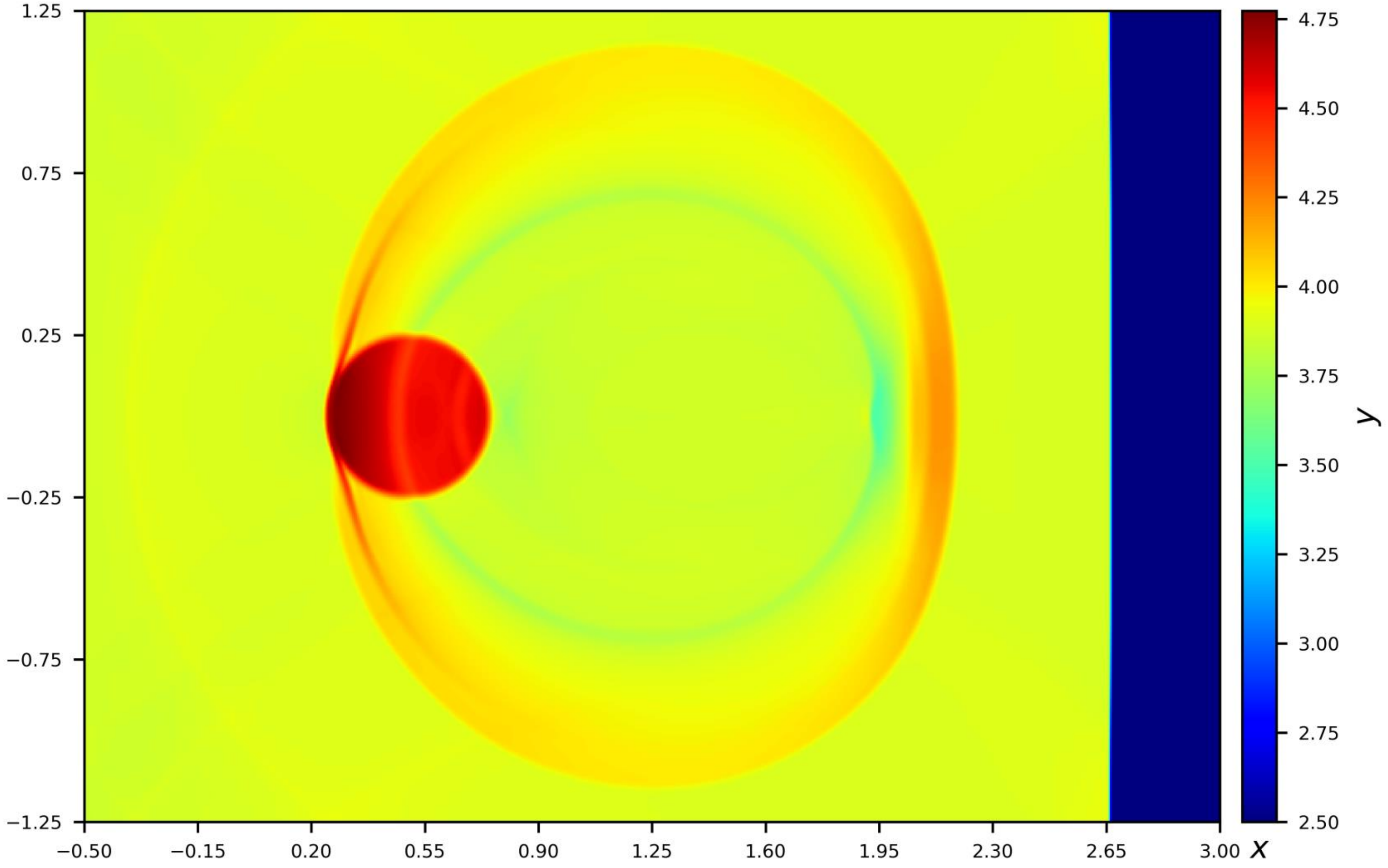

c)

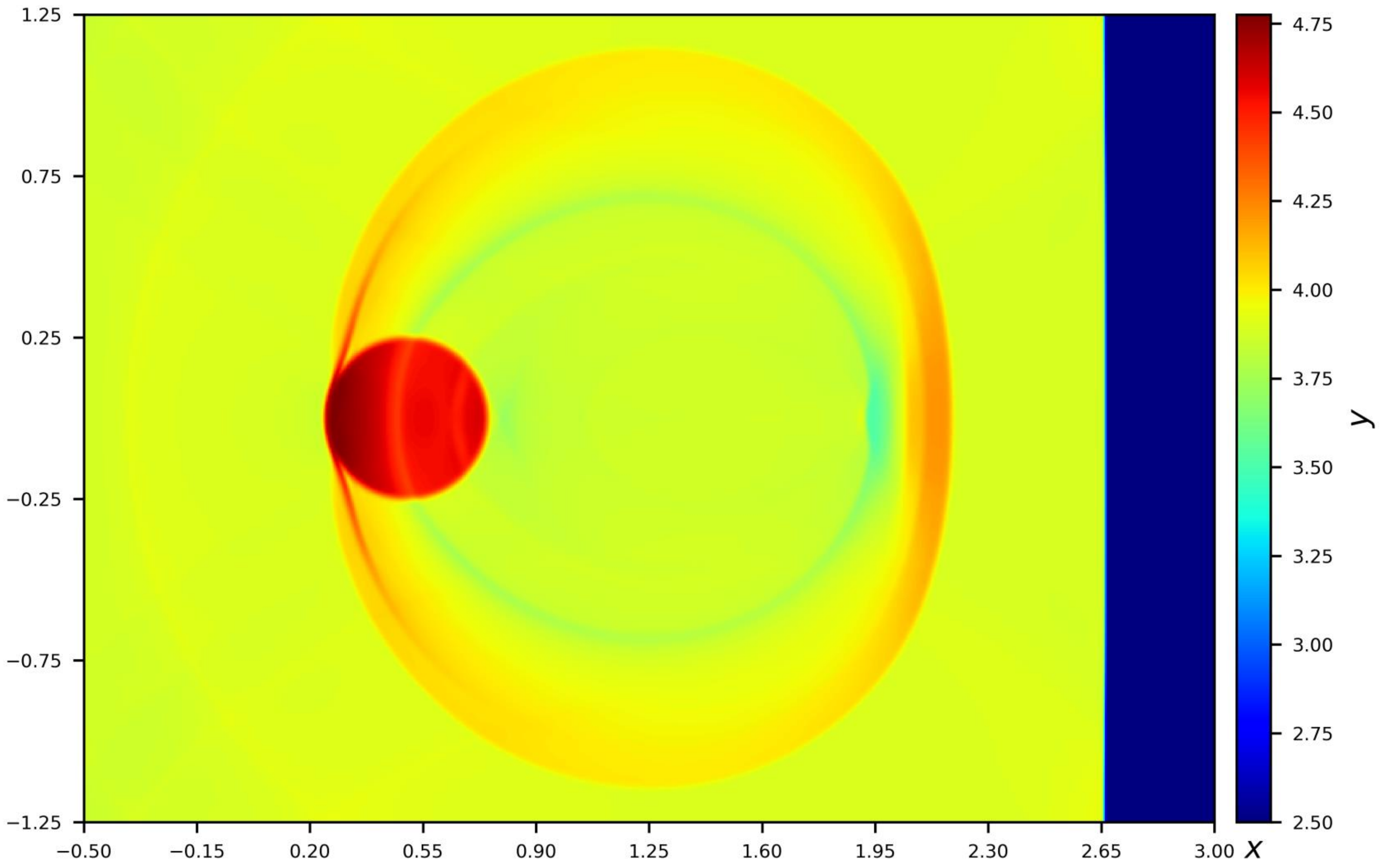

*Fig. 15a, b, c) Debris Flow : Shock-bubble interaction problem using the 3rd, 5th and 7th order accurate HLL-based FD-WENO schemes with 700×500 zones. Height of the fluid phase has been shown.*

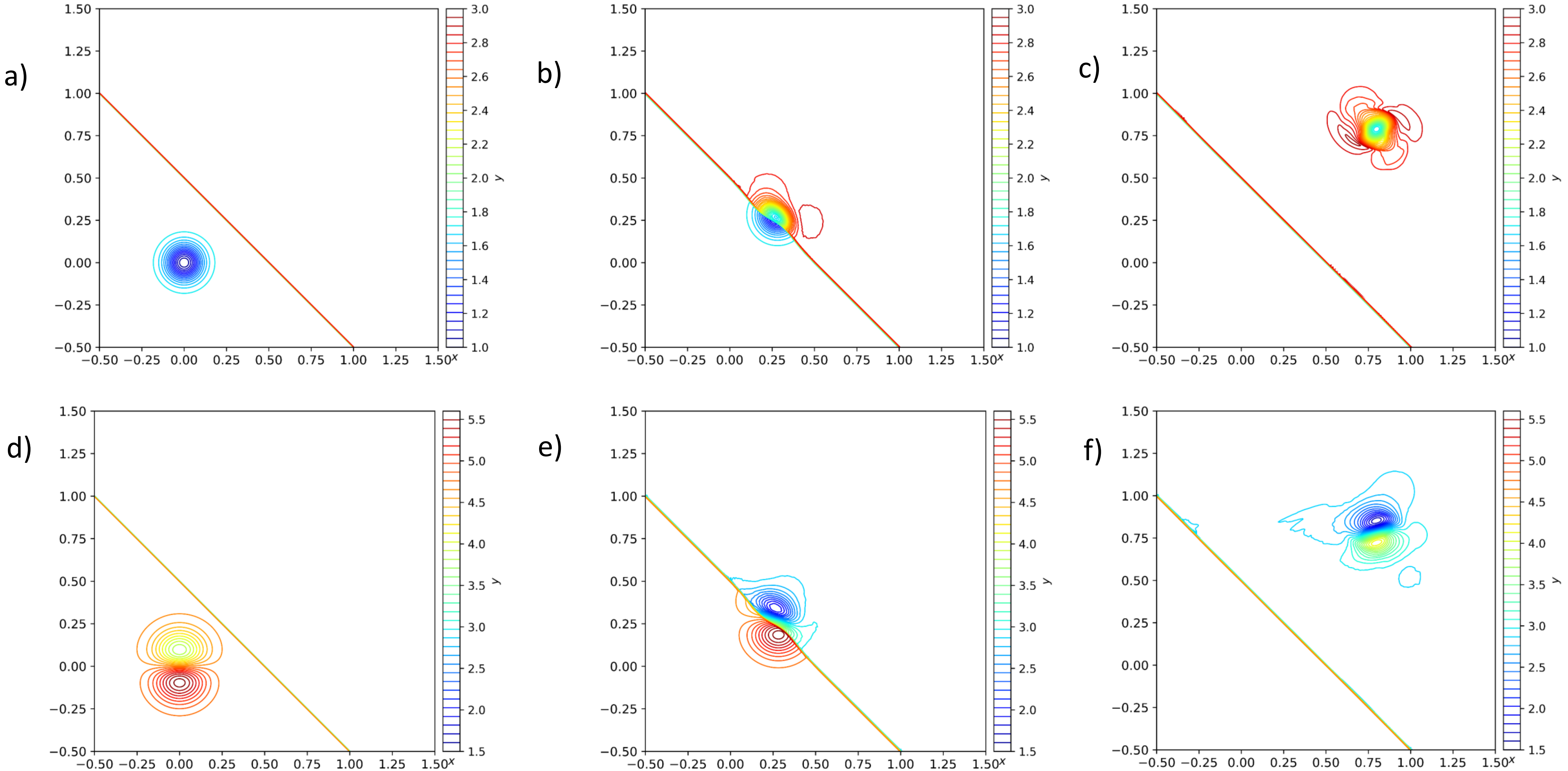

*Fig. 16) Debris Flow: Shock-Vortex Interaction using the 7th order accurate HLL-based FD-WENO scheme with 600×600 zones at time levels t=0.0, 0.06 and 0.24. Figs. 16a, 16b and 16c show the solid height at times t=0.0, 0.06, 0.24. Figs. 16d, 16e, 16f show the solid x-velocity at times t=0.0, 0.06, 0.24. For the height, 40 contours were fit between a range of 1.0 and 3.0. For the velocity, 40 contours were fit between a range of 1.5 and 5.6.*

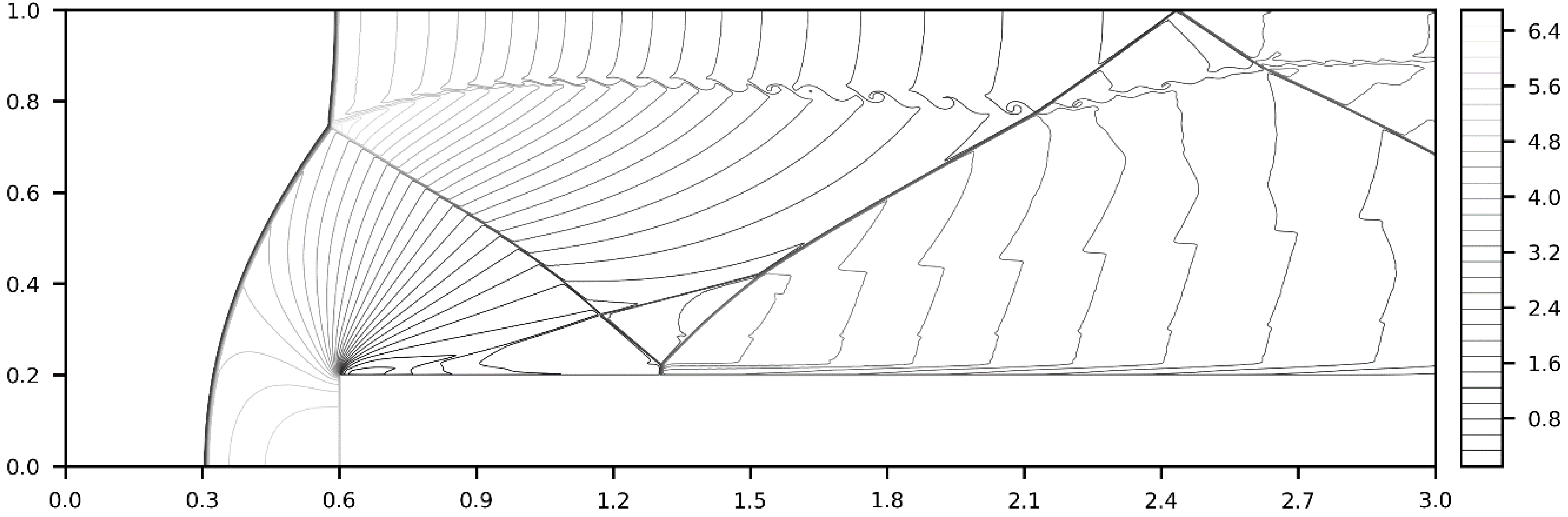

*Fig. 17) EULER: Forward facing step problem using the 7th order accurate HLL-based FD-WENO scheme with 1440×480 zones. 30 contours were fit between a range of 0.105 and 6.699.*

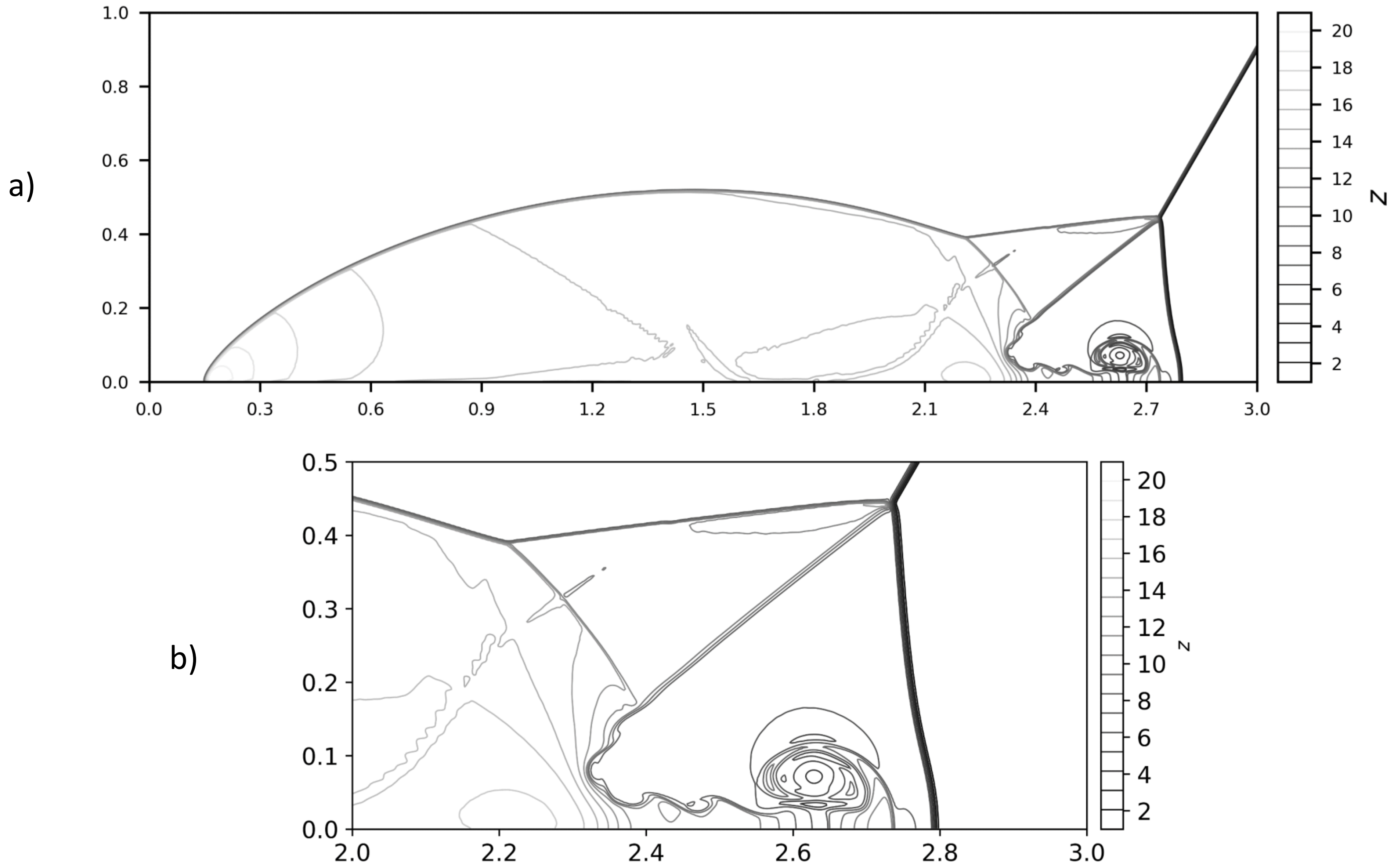

*Fig. 18) EULER: Double Mach reflection problem using the 7th order accurate HLL-based FD-WENO scheme with 1920×480 zones. Fig. 18a shows the density profile at the 7th orders. Fig. 18b shows the detailed view of the density profile. 20 contours were fit between a range of 1.0 and 20.97.*

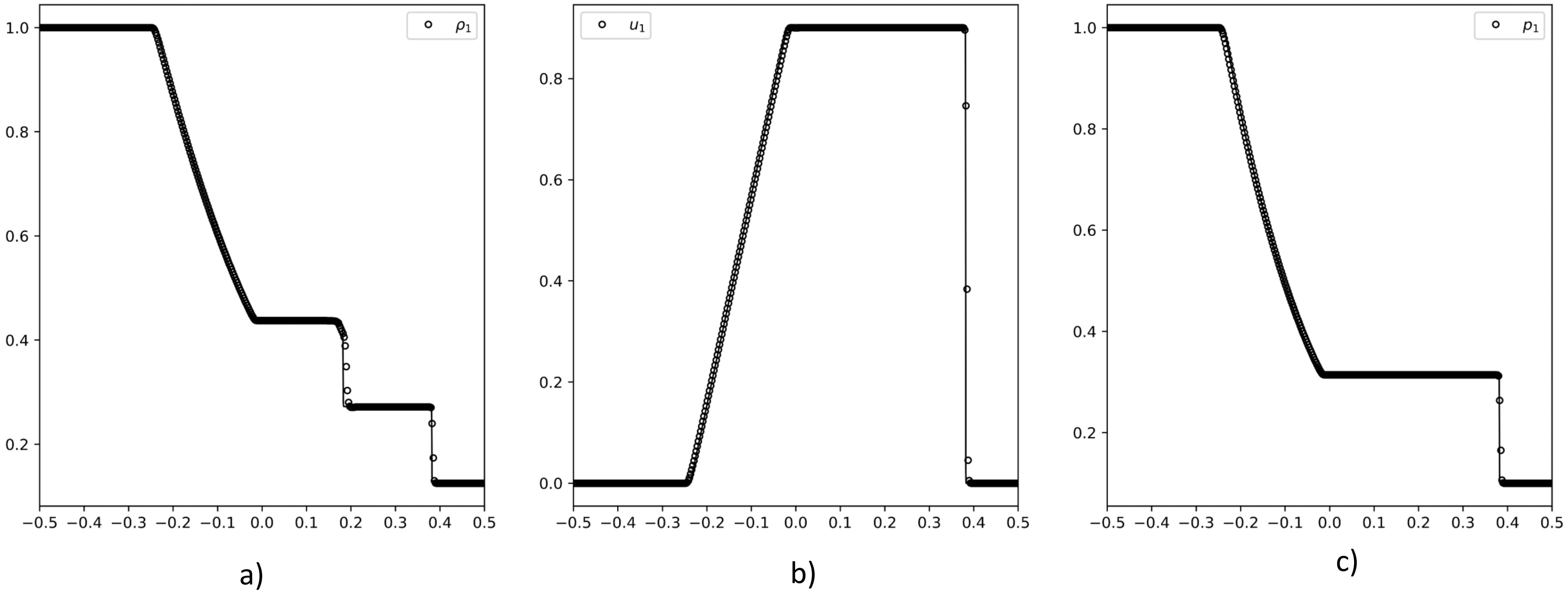

*Fig. 19) Baer-Nunziato with stiff source: Results for the one-dimensional Riemann Problem using the 7th order accurate HLL-based FD-WENO scheme with 400 zones. Figs. 19a, 19b and 19c show the solid density, solid x-velocity and solid pressure.*

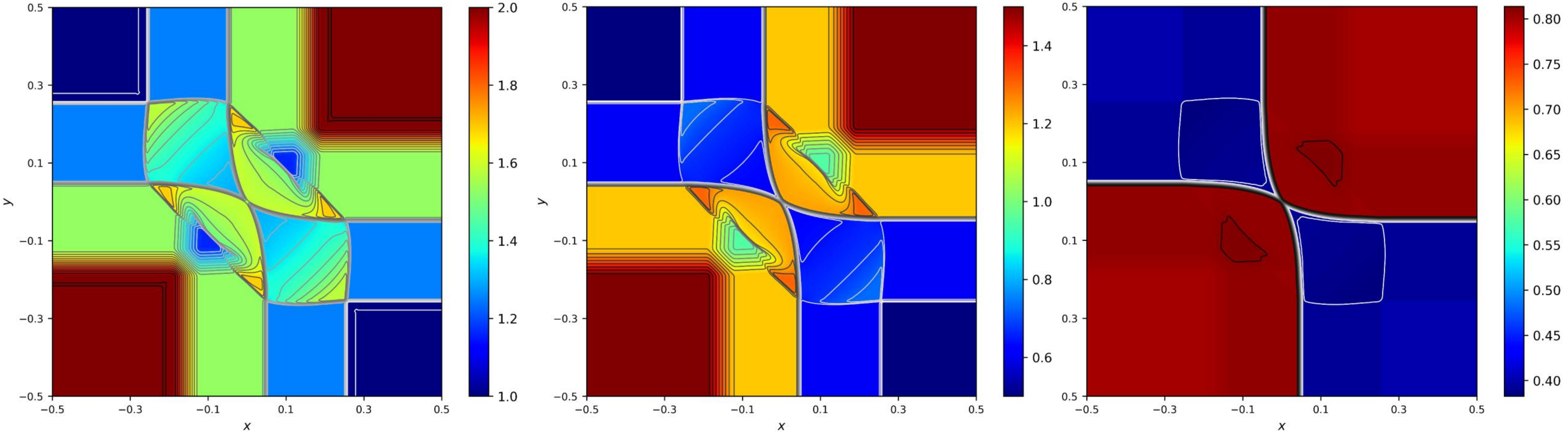

a) b) c)

*Fig. 20) Baer-Nunziato with stiff source: Results for the two-dimensional Riemann Problem using the 7th order accurate HLL-based FD-WENO scheme with 400 ×400 zones. Fig. 20a shows the solid density, Fig. 20b shows the gas density and Fig. 20c shows the solid volume fraction. 30 equidistant contour lines are shown over the color plots.*